\documentclass{article}

\usepackage[english,USenglish,UKenglish]{babel}
\usepackage{amssymb}
\usepackage{amstext}
\usepackage{amsmath}

\usepackage{amsmath}
\usepackage{txfonts}
\usepackage{theorem}
\usepackage{graphicx}
\usepackage[all]{xy}
\usepackage{latexsym}
\usepackage{enumitem}

\newtheorem{defn}{Definition}[section]
\newtheorem{thm}[defn]{Theorem}
\newtheorem{prop}[defn]{Proposition}
\newtheorem{lem}[defn]{Lemma}
\newtheorem{claim}[defn]{Claim}

\newtheorem{cor}[defn]{Corollary}
\newtheorem{prob}[defn]{Problem}
{\theorembodyfont{\upshape} \newtheorem{rem}[defn]{Remark}}
{\theorembodyfont{\upshape} \newtheorem{ex}[defn]{Example}}
{\theorembodyfont{\upshape} \newtheorem{proof}{Proof}}

\newcommand{\ora}{\overrightarrow}

\allowdisplaybreaks[4]

\begin{document}

\title{Fukaya categories in Koszul duality theory}
\author{Satoshi Sugiyama\thanks{sugi3@ms.u-tokyo.ac.jp}}
\date{\today}
\maketitle

\begin{abstract}

In this paper, we define $A_{\infty}$-Koszul duals for directed $A_{\infty}$-categories
in terms of twists in their $A_{\infty}$-derived categories. Then, we
compute a concrete formula of $A_{\infty}$-Koszul duals for path algebras with directed $A_n$-type
Gabriel quivers. To compute an $A_\infty$-Koszul dual of such an algebra $A$, we construct a directed subcategory of a Fukaya category which are $A_\infty$-derived equivalent to the category of $A$-modules and compute Dehn twists as twists.
The formula unveils  all the ext groups of simple modules of the parh algebras
and their higher composition structures. 

\end{abstract}

\setcounter{tocdepth}{2}
\tableofcontents

\section{Introduction}\label{sec:intro}

The purpose of this paper is to give a new expression of $A_\infty$-Koszul duals of certain path algebras with relations (Theorem \ref{thm:MAIN1TheFormulaOfAInftyKoszulDual}). We use the technique of the Fukaya categories and Dehn twists to compute $A_\infty$-Koszul duals.
Our approach does not contain anything new in the standpoint of the abstract
theory of Koszul duality. However, we show that the technique of the Fukaya categories can be used for a concrete computation of an algebraic problem. Moreover, our description computed via the Fukaya categories provides a new way of understanding of Koszul duality as a duality between higher products and relations.

The Fukaya categories are  $A_\infty$-categories associated to symplectic
manifolds defined by using the technique of Floer theory \cite{FOOO10}, \cite{Se08}.
The Fukaya categories are mainly studied in the context of homological mirror
symmetry \cite{Ko94}.
The concept of  Fukaya categories emerges in the context of Koszul duality in the paper
of A. J. Blumberg, R. L. Cohen, and C. Teleman \cite{BCT09} and the paper
of T. Etg\"{u} and Y. Lekili \cite{EL16}.
These papers state that  End $A_\infty$-algebras of two certain objects
in some Fukaya categories are Koszul dual to each other.
Therefore, they say that the Koszul duality patterns emerge in the context
of
Fukaya categories.
In our paper, the direction is opposite.
We use the Fukaya categories to compute  $A_\infty$-Koszul duals of path
algebras with relations.
Therefore we can say that Fukaya categories emerge in the context
of Koszul duality theory.

Before we see the main theorem of this paper, let us review the fundamental results about Koszul duality in \cite{Lo86}. (The results presented here is a simplified version.) Let $A_0 = k$ be a field, $A_1$
be a finite dimensional vector space and $I$ be a subspace of $A_1 \otimes
A_1$.
Define $A \coloneqq T(A_1)/I$ as the quotient algebra of the tensor algebra
of $A_1$ over $A_0 = k$.
Then, we have $E \coloneqq \textrm{Ext}_A(k, k) \cong T(A_1^*)/I^\perp$, where $(-)^*$ is the linear dual over $k$ and $I^\perp \subset A_1^* \otimes A_1^*$
is the annihilating submodule of $I \subset A_1 \otimes A_1$ (we use the natural
isomorphism between $A_1^* \otimes A_1^*$ and $(A_1 \otimes A_1)^*$).
Let us  fix an isomorphism between $A_1$ and $A_1^*$. Then, $I$ and $I^\perp$
are mutually complemental.
Hence, we can say that the products and relations interchange between $A$ and $E$.
By the above computation, $\textrm{Ext}_E(k, k)$ is naturally isomorphic
to $A$.
This is what we call Koszul duality and we can say that Koszul duality is
a duality between products and relations represented by the Yoneda Ext algebra.
Moreover, certain derived categories of $A$ and $E$ are equivalent \cite{BGS96}.
(In that paper, the setting above is generalized to the case of that $A_0$ is a finite dimensional semi-simple algebra.)

Nowadays, many phenomena related to the Koszul duality are widely observed,
for example, the Koszul duality for Koszul algebras \cite{Pr70}, \cite{Lo86},
 \cite{BGS96}, its generalisation to augumented-$A_\infty$ algebras \cite{LPWZ04},
a generalisation to Koszul operads \cite{GK94}, \cite{Va07}, \cite{LV12},
and its relation to the study of symplectic geometry \cite{EL16} and mirror
symmetry \cite{AKO08}.

In this paper, we are interested in the case that there exist higher degree (homogenous) relation, i.e. for the algebra $A = T(A_1)/I$ with $I (\not \subset A_1^{\otimes 2})
\subset \bigoplus_{d \geq 2}A_1^{\otimes d}$.
In general, there is no easy description of $E$.
Moreover, the ext algebra $\textrm{Ext}_E(A_0, A_0)$ and $A$ are no longer isomorphic.
However, we can overcome this difficulty by referring the results in \cite{LPWZ04}.
They generalise the concept of Koszul dual to the augmented $A_\infty$-algebras.
After that, they prove that the twice dual is quasi-isomorphic to the original
augmented $A_\infty$-algebra and their derived categories are equivalent (under some finiteness
condition).
The above algebra $A$ is an example of an augmented $A_\infty$-algebra,
so we have its dual. But the description is too complicated and we can not
interpret the Koszul dual as the duality between products and relations.

In this paper, we define the notion of $A_\infty$-Koszul dual for directed $A_\infty$-categories (Definition \ref{def:AInftyKoszulDual}) and  present an explicit description of $A_\infty$-Koszul dual
of certain class of path algebras with relations (Theorem \ref{thm:MAIN1TheFormulaOfAInftyKoszulDual}) which enable us to understand the Koszul duality as a duality between higher products and relations. The notion of $A_\infty$-Koszul dual is a natural generalisation. This is supported by the following two corollaries: the $A_\infty$-Koszul
dual $\mathcal{C}$ of $\mathcal{B}$ is naturally quasi-isomorphic to $\mathcal{A}$
(Corollary \ref{cor:KoszulDualIsDuality}); $\mathcal{A}$ and its Koszul dual
$\mathcal{B}$ are $A_\infty$-derived equivalent, i.e. $Tw\mathcal{A} \cong
Tw\mathcal{B}$ (hence, in particular, they are derived equivalent , i.e.
$D\mathcal{A}
\cong D\mathcal{B}$) (Corollary \ref{cor:DerEquivOfKoszulDuals}).

The computation of the $A_\infty$-Koszul dual takes place in the Fukaya categories
of exact Riemann surfaces.
The rough sketch of the computation is as follows.
In general, the Koszul dual can be computed by the operation in the derived
category called twist.
First we ``embed" our directed $A_\infty$-category $\mathcal{A} = \mathcal{A}(R)$
into the Fukaya category $\mathcal{F} = Fuk(M)$ of an exact Riemann
surface $M$ constructed by using the data of relations of $R$.
Seidel proved in \cite{Se08} that the twists are ``quasi-isomorphic" to the
Dehn twists in the Fukaya category.
Thus, we compute the Dehn twists of the objects which are lying in the image
of the ``embedding" $\mathcal{A} \hookrightarrow Fuk(M)$.
Finally, we investigate how the resulting curves intersect and encircle polygons
to compute the morphism spaces and their higher compositions.
After that, we find that there is a $(d+1)$-gon in $M$ corresponding to a
degree
$d$ relation, and the $(d+1)$-gon generates the $d$-th higher composition
$\mu^d$.
This is our geometric explanation of the duality between higher products and relations.
Some typical example is presented in Corollary \ref{cor:typicalEx} and Subsection \ref{subsec:SomeExamples}.

Here, we fix some notations we often use.
In this paper, 
$k$ is a fixed field;
all categories are of over $k$;
all graded vector spaces are assumed to have the property that theie total
dimensions are finite;
for a graded vector space $V = (V^d)_{d \in \mathbb{Z}}$, $V[r] \coloneqq
(V^{d+r})_{d \in \mathbb{Z}}$ is the $r$-th shift of $V$;
all modules are always right modules;
all manifolds are oriented;
all the additional structures on manifolds are assumed to be compatible with
their orientations;
the character $\mathcal{F}$ always stands for the Fukaya category $Fuk(M)$
of $M$ where $M$ is ``the" exact symplectic manifold we consider in each
paragraph; 
if $M$ has some subscripts like $M_1$ then $\mathcal{F}_1$ stands for the
Fukaya category of $M_1$,
unless otherwise stated.

The structure of this paper is as follows.
In section 2, we prepare the algebraic notions and define the $A_\infty$-Koszul
dual.
In section 3, we prepare the geometric notions, e.g. exact symplectic
manifolds and their Fukaya categories. At the last part of the section, we present
the key theorem proved by Seidel which states the equivalence of algebraic
twists and Dehn twists.
In section 4, we state the main theorem. In section 5 and 6, we construct
exact Riemann surfaces whose Fukaya categories are the targets of the ``embedding"
from directed $A_\infty$-categories.
In section 7, we do the computation of $A_\infty$-Koszul duals, i.e. the computation
of Dehn twists.
The computation and the formula of $A_\infty$-Koszul duals are the main ingredients
of this paper.

{\bf Acknowledgement}

I would like to thank my supervisor Toshitake Kohno for giving me great advice and navigating me and this study to an appropriate direction.
I also want to thank A. Ikeda for teaching me about Koszul duality theory and to F. Sanda, M. Kawasaki, T. Kuwagaki, J. Yoshida, and R. Sato for fruitful discussion.
Finally, I am deeply greatful to my friend I. Hoshimiya, A. Kiriya, R. Shibuki, and Y. Todo for supporting me when I was in difficult situations.

This work was supported by the Program for Leading Graduate 
Schools, MEXT, Japan.

\section{Algebraic preliminaries}\label{sec:AlgDef}

In this section, we review the definitions of algebraic objects we use in
this paper, and define the $A_{\infty}$-Koszul dual, the key concept in
this paper.
For the notation of signs, we follow Seidel's notation in \cite{Se08}.
The definition of the Koszul dual for $A_{\infty}$-algebras with some properties already
exists \cite{EL16}, \cite{LPWZ04}.
Our construction is a generalisation to directed $A_{\infty}$-categories.

\subsection{Basic definitions and properties of $A_{\infty}$-categories}

\begin{defn}[$A_\infty$-category]\label{def:AInfCat}
An $A_\infty$-category $\mathcal{A}$, consists
of the following data:
\begin{enumerate}
\item a set $Ob(\mathcal{A})$,
\item a $\mathbb{Z}$-graded vector space $\displaystyle \hom_{\mathcal{A}}(X,
Y) = \bigoplus_{i \in \mathbb{Z}}\hom_\mathcal{A}^i(X, Y)$ for each $X, Y
\in Ob(\mathcal{A})$,
\item maps called higher composition maps
\begin{multline*}
\mu^d \colon \textrm{hom}_\mathcal{A}(X_{d-1}, X_d) \otimes \hom_\mathcal{A}(X_{d-2},
X_{d-1}) \otimes \cdots \otimes \hom_\mathcal{A}(X_0, X_1) 
\\ \to \hom_\mathcal{A}(X_0, X_d)[2-d],
\end{multline*}
for $d \geq 1$ and $X_0, X_1, \dots , X_d \in Ob(\mathcal{A})$.
\end{enumerate}
We impose that the \(\mu\)'s satisfy the \(A_\infty\)-associativity relation:
\begin{equation*}
\sum _{i, j,l} (-1)^{\bigstar _i} \mu ^l (a_d, \dots , a_{i+j+1}, \mu ^j
(a_{i+j}, \dots , a_{i+1}), a_i, \dots a_1) = 0
\end{equation*}
for $d \geq 1$, where $\displaystyle \bigstar _i = \sum _{1 \leq l \leq i}
( |a_l| - 1)$
, $\, (|a_i| = deg(a_i))$.
\end{defn}

Let us see the first few $A_\infty$-relations.
The $A_\infty$-relation of $d = 1$ is $\mu^1(\mu^1(a_1)) = 0$ and $\deg (\mu^1)
= 2-1 = 1$. Hence, $(\hom_\mathcal{A}(X_0, X_1), \mu^1)$ forms a cochain
complex.
The second case, the relation is $\mu^1(\mu^2(a_2, a_1)) + \mu^2(a_2, \mu^1(a_1))
-(-1)^{|a_1|} \mu^2(\mu^1(a_2), a_1) = 0$.
When we write $da = (-1)^{|a|}\mu^1(a)$ and $a_2 \circ a_1 = (-1)^{|a_1|}\mu^2(a_2,
a_1)$, the relation is written by $d(a_2 \circ a_1) = da_2 \circ a_1 + (-1)^{|a_2|}a_2
\circ da_1$. Thus, the second relation expresses the graded Leibniz' rule.
If all the higher composition maps are zero, i.e. $\mu^d = 0$ for $d \geq
3$, then the $A_\infty$-category is nothing but a dg category by the above
$d$ and $- \circ -$. Therefore, the notion of $A_\infty$-categories is a generalisation
of dg categories.

The third relation is somewhat complicated:
\[ a_3 \circ (a_2 \circ a_1) - (a_3 \circ a_2) \circ a_1 = \pm d(\mu^3(a_3,
a_2, a_1)) \pm \mu^3(da_3, a_2, a_1) +(\text{other two terms}). \]
In general, the right hand side does not vanish, so the composition defined
by $\mu^2$ is not associative.
However, $\mu^3$ forms a homotopy between $a_3 \circ (a_2
\circ a_1)$ and $(a_3 \circ a_2) \circ a_1$, hence $\mu^2$ defines an associative
composition on cohomology level.
We define the \textit{cohomology category} $H(\mathcal{A})$ by
$Ob(H(\mathcal{A})) \coloneqq Ob(\mathcal{A})$,
$\hom_{H(\mathcal{A})}(X_0, X_1) \coloneqq H(\hom_\mathcal{A}(X_0, X_1),
\mu^1)$,
and $[a_2] \circ [a_1] \coloneqq (-1)^{|a_1|}[a_2 \circ a_1]$. The resulting
category $H(\mathcal{A})$ has an associative composition.
Thus, we say that $\mu^2$ is homotopy associative.
We also define $H^0(\mathcal{A})$ in the obvious way.

We don't assume that the $A_\infty$-category admits identity morphisms,
so $H(\mathcal{A})$ and $H^0(\mathcal{A})$ may not have identity morphisms.
If $H(\mathcal{A})$ admits identity morphisms for each object, then we say
that $\mathcal{A}$ is \textit{cohomologically unital} or \textit{c-unital}.
In this paper, all the $A_\infty$-categories are of c-unital unless otherwise
stated.
We say that two objects $X_0$ and $X_1$ in an $A_\infty$-category are \textit{quasi-isomorphic}
if they are isomorphic in $H^0(\mathcal{A})$.

We do not present the definitions of $A_\infty$-functors, quasi-equivalences,
and quasi-isomorphisms of $A_\infty$-categories here. These are generalisations
in the case of dg categories. For precise definition and properties, please
refer Section 1 and 2 in \cite{Se08}.

\subsection{Directed $A_\infty$-categories and $A_\infty$-Koszul duals}\label{subsec:AInftyKoszulDual}

In this paper, we mainly consider the directed $A_\infty$-categories.

\begin{defn}

An $A_\infty$-category $\mathcal{A}$ is said to be directed when
\begin{enumerate}
\item the set $Ob(\mathcal{A})$ is finite,
\item $\hom_\mathcal{A}(X, X) = k \cdot 1_X$, and
\item there exists a total order on $Ob(\mathcal{A})$ such that the hom space
$\hom_{\mathcal{A}}(X, Y) \neq 0$ only when $X\leq Y$.
\end{enumerate}

\end{defn}

For a totally ordered finite set $A$, we have a canonical isomorphism $A
\cong \{0<2< \cdots <n\}$. Therefore we write the objects of an directed $A_\infty$-category
as $0<1< \cdots <n$, $X_0 < X_1 < \dots < X_n$, and so on.

\begin{defn}\label{def:DirSUbcat}

Let $\mathcal{A}$ be an $A_\infty$-category and $\boldsymbol{Y} = (Y_1, Y_2,
\dots , Y_n)$ be a collection of objects in $\mathcal{A}$.
Then, we define the associated directed subcategory $\mathcal{A}^\to(\boldsymbol{Y})$
of $\mathcal{A}$ by setting $Ob(\mathcal{A}) = \{ Y_1, Y_2, \dots , Y_n \}$,
\begin{align*}
\hom_{\mathcal{A}^\to(\boldsymbol{Y})}(Y_i, Y_j) = 
\begin{cases}
\hom_\mathcal{A}(Y_i, Y_j) & (i < j) \\
k \cdot e_i & (i = j) \\
0 & (i > j),
\end{cases}
\end{align*}
and $\mu$'s of $\mathcal{A}^\to(\boldsymbol{Y})$ are canonically induced
from those of $\mathcal{A}$.

\end{defn}

Now, we begin the definition of $A_\infty$-Koszul duals.
For an $A_\infty$-category $\mathcal{A}$, we call an $A_\infty$-functor $\mathcal{M}$
from $\mathcal{A}^{op}$ to $C(k) \coloneqq C^b_{dg}(k)$ a (right) $A_\infty$-module,
where $\mathcal{A}^{op}$ is the opposite $A_\infty$-category of $\mathcal{A}$
and $C^b_{dg}(k)$ is the dg category of bounded cochain complexes of finite dimensional
$k$ vector spaces considered as an $A_\infty$-category.
It is known that such $A_\infty$-modules form a triangulated dg category
$\mathcal{Q} \coloneqq mod(\mathcal{A})$.
(Note that all the hom spaces of this category are finite dimensional iff
$\# Ob(\mathcal{A}) < \infty$.)
Let $\mathcal{A}$ be a directed $A_\infty$-category with its object set $\{ 0
< 1 < \cdots < n \}$.
We define an $A_\infty$-module $\mathcal{S}(j)$ for $j \in Ob(\mathcal{A})$
determined by the data
\begin{equation*}
\mathcal{S}(j)(i) =
\begin{cases}
k & \text{if }i = j, \\
0 & \text{if }i \neq j
\end{cases}
\end{equation*}
 and call it a simple $\mathcal{A}$-module corresponds to $j$, where we consider
$k$ as a one-dimensional cochain complex concentrated in the degree zero part.
Then, it is known that the full sub $A_\infty$-category $\mathcal{A}^!_\infty$
of $\mathcal{Q}$ with object set $\{ \mathcal{S}(n) < \mathcal{S}(n-1) <
\dots < \mathcal{S}(0)\}$ forms a directed $A_\infty$-category.
Hence, $\mathcal{A}^!_\infty$ is naturally isomorphic to $\mathcal{Q}^\to(\boldsymbol{\mathcal{S}})$,
where $\boldsymbol{\mathcal{S}} = (\mathcal{S}(n), \mathcal{S}(n-1), \dots
, \mathcal{S}(0))$ is a collection of objects in $\mathcal{Q}$.
The details can be found in (5j) and (5o) in \cite{Se08}.

\begin{defn}\label{def:AInftyKoszulDual}

Let $\mathcal{A}$ be a directed $A_\infty$-category  with object set $\{
0 < 1 < \cdots < n\}$. A directed $A_\infty$-category $\mathcal{B}$ quasi-isomorphic
to $\mathcal{A}^!_\infty$ is called an $A_\infty$\textit{-Koszul dual} of
$\mathcal{A}$.

\end{defn}

\begin{rem}

The above definition is an analogy or a category version of the definition in \cite{BGS96}.
In that paper, they deal with Koszul ring $A$ and give a different definition
of its Koszul dual $A^!$.
However, Theorem 2.10.1 in that paper states that $\textrm{Ext}^\bullet_A(k, k) \cong (A^!)^\textrm{opp}$
canonically.
Even though there exist many different notations, we can translate from one to the other.

Also, our definition is an analogy of the definition in \cite{LPWZ04}.
In that paper, they define Koszul dual $E(A)$ (in their notation) for Adams connected $A_\infty$-algebra $A$ by $E(A) \coloneqq \text{RHom}_{A^\circ}(k, k)$.
The right hand side of the definition is a straightforward generalisation of the definition in \cite{BGS96}, hence the definitions in that paper and in our paper shares the common origin.

In this paper, we treat with $A_\infty$-categories, not $A_\infty$-algebras and we focus on the very special case, directed $A_\infty$-categories.

\end{rem}

\begin{ex}\label{ex:DirAlgAsAInftyCat}

Let $R = k(\ora{\Delta}, \rho)$ be a path algebra with relations over a finite directed quiver $\overrightarrow{\Delta}$.
Here, a finite quiver is a quiver with a finite set of vertices (which we write $\Delta_0$) and a finite set of arrows (which we write $\Delta_1$);
a directed quiver is a quiver without oriented cycles.
We can see $R$ as an $A_\infty$-category $\mathcal{A} = \mathcal{A}(R)$ by
setting $Ob(\mathcal{A}) = \Delta_0$, $\hom_\mathcal{A}^0(i, j) = e_j A e_i$,
$\hom_\mathcal{A}^{d}(i, j) = 0$ for $d \neq 0$, $\mu^2$ is induced from
the product structure of $A$, and $\mu^d = 0$ for $d \neq 2$.
(We write the product of two paths $\alpha$ from $i$ to $j$
and $\beta$ from $j$ to $l$ as $\beta \alpha$ later on.)
Now, the dimension $\dim_k R$ as a $k$ vector space is finite since its quiver
$\overrightarrow{\Delta}$ has no oriented cycles. Thus, we can deduce that
$mod(\mathcal{A})$ and $C(R) \coloneqq C^b_{dg}(R)$ are naturally isomorphic
as triangulated dg categories, where $C^b_{dg}(R)$ is the dg category of
finitely generated $R$-modules.
(Recall that a functor from $\mathcal{A}(R)^{op}$ (considered as $k$-linear
category) to the category of finite dimensional $k$ vector spaces $vect(k)$
can be naturally  considered as a right $R$-module.)
The natural isomorphism maps $\mathcal{S}(j)$ in $mod(\mathcal{A})$ for $j
\in \Delta_0$ into  the simple module $S(j)$ in $C^b_{dg}(R)$ corresponds to  $j \in
\Delta_0$.
Set a graded algebra $R^!_{dg} \coloneqq \hom^*_{C^b_{dg}(R)} \left( \bigoplus
\widetilde{S}^\bullet (j), \bigoplus \widetilde{S}^\bullet(j) \right)$, where
$\widetilde{S}^\bullet (j)$ is a projective resolution of $S(j)$ and the
direct sum is taken over $\Delta_0$  . We call it the \textit{dg Koszul dual}
of $R$.
Then we can compute $\mathcal{A}(R)^!_\infty$ by $Ob(\mathcal{A}(R)^!_\infty)
= \Delta_0$, $\hom_{\mathcal{A}(R)^!_\infty}^d(i, j) = \hom^d_{C^b_{dg}(R)}
\left( \bigoplus
\widetilde{S}^\bullet (i), \bigoplus \widetilde{S}^\bullet(j) \right)$, and
$\mu$'s are induced from the differential $d$ and the product structure $- \cdot
-$ of $R^!_{dg}$.

If our algebra $R$ is Koszul, equivalently, the relations are of quadratic,
then the cohomology algebra $H(R^!_{dg})$ is nothing but the Koszul dual
$R^!$ of $R$.
Hence, the dg Koszul dual $R^!_{dg}$ is a generalisation of the Koszul dual
to general path algebras over finite directed quivers.
The dg Koszul dual $R^!_{dg}$ can be reconstructed from $\mathcal{A}(R)^!_\infty$
by $\displaystyle R^!_{dg} = \bigoplus_{i, j \in \Delta_0} \hom_{\mathcal{A}(R)^!_{dg}}(\mathcal{S}(i),
\mathcal{S}(j)) = \hom_{mod(\mathcal{A}(R)^!_\infty)}\left( \bigoplus \mathcal{S}(i),
\bigoplus \mathcal{S}(j) \right)$, where the last two direct sums are taken
over $\Delta_0$.

\end{ex}

We finish this subsection by collecting some useful lemmas from \cite{Se08}.

\begin{lem}[(5n) in \cite{Se08}]\label{lem:FunctorialityOfDirSubcat}

Let $\mathcal{F} \colon \mathcal{A} \to \mathcal{B}$ be a cohomologically
full and faithful (c-full and faithful in short) $A_\infty$-functor and $\boldsymbol{Y}
= (Y_1, Y_2, \dots , Y_n)$ be a collection of objects in $\mathcal{A}$.
Then, there exists a canonical quasi-isomorphism between $\mathcal{A}^\to(\boldsymbol{Y})$
and $\mathcal{B}^\to(\mathcal{F}\boldsymbol{Y})$, where $\mathcal{F}\boldsymbol{Y}
= (\mathcal{F}Y_1, \mathcal{F}Y_2, \dots \mathcal{F}Y_n)$.

\end{lem}

\begin{lem}[Lemma 5.21 in \cite{Se08}]\label{lem:EquivOfDirSubcat}

Let $\boldsymbol{Y}$ and $\boldsymbol{Y}'$ be collections of objects in $\mathcal{A}$
and these objects are pairwise quasi-isomorphic, i.e. $Y_j
\cong Y'_j$ in $H^0(\mathcal{A})$ for every $j$.
Then, the associated directed subcategories $\mathcal{A}^\to(\boldsymbol{Y})$
and $\mathcal{A}^\to(\boldsymbol{Y}')$ are quasi-isomorphic.

\end{lem}

\subsection{$A_\infty$-Koszul duals and twists}

In this subsection, we develop the method to compute an $A_\infty$-Koszul dual
of a given directed $A_\infty$-category $\mathcal{A}$.
All the details and precise definitions can be found in Chapter I of \cite{Se08}.

First, we fix some notations. For an $A_\infty$-category $\mathcal{A}$, we
define the category of $\mathcal{A}$-modules $\mathcal{Q} \coloneqq mod(\mathcal{A})
= fun(\mathcal{A}^{op}, C(k))$.
For such categories, we can define the Yoneda embedding functor $\iota \colon
\mathcal{A} \to \mathcal{Q}$, by setting $(\iota X)(Y) = \hom_\mathcal{A}(Y,
X)$.
We set the triangulated $A_\infty$-category $Tw \mathcal{A}$ by the full
subcategory generated as triangulated $A_\infty$-category by the
objects which are lying in the image of the Yoneda embedding $\iota(Ob(\mathcal{A}))$.
Now, we have three embeddings of $A_\infty$-categories, $\mathcal{A} \hookrightarrow
Tw\mathcal{A} \hookrightarrow \mathcal{Q} = mod(\mathcal{A})$. These three
embeddings are known to be c-full and faithful.

For $X \in Ob(\mathcal{A})$ and $\mathcal{M} \in Ob(\mathcal{Q})$, we can
define the \textit{twist of} $\mathcal{M}$ \textit{along} $X$ , which is
wtritten by $\mathcal{T}_X\mathcal{M}$, by the mapping cone of the evaluation morphism $\iota X \otimes
\hom_\mathcal{Q}(\iota
X, \mathcal{M}) \to \mathcal{M}$. This is a generalisation of the case when
$\mathcal{A} = \mathcal{A}(R)$ as in Example \ref{ex:DirAlgAsAInftyCat}.
If there exists $Z \in Ob(\mathcal{A})$ such that $\iota Z$ and $\mathcal{T}_X
(\iota Y)$ are quasi-isomorphic, we write $Z = T_X Y$ and call it a \textit{twist
of} $Y$ \textit{along} $X$.
This is a fact that $Tw \mathcal{A}$ is closed under twist.
There are two remarks on the notion of twists.
The first one is that such a $Z$ may not be unique. Therefore whenever we write $T_XY$, we choose
one of such objects.
The second one is that when we write $T_X Y$, we always assume the existence of the representative of $\mathcal{T}_X \iota Y$.
Finally, the following holds:

\begin{lem}[Lemma 5.24. in \cite{Se08}]\label{lem:ResultOfTwistIsSimple}

Let $\mathcal{A}$ be a directed $A_\infty$-category with object set $\{ X_0
< X_1 < \cdots < X_n \}$, and set $\mathcal{S}'(j) \coloneqq T_{X_0}T_{X_1}
\cdots T_{X_{j-1}}X_j \in Ob(Tw\mathcal{A}) \hookrightarrow Ob(\mathcal{Q})$.
Then the resulting object $\mathcal{S}'(j)$ is quasi-isomorphic (in $\mathcal{Q}$)
to the simple module $\mathcal{S}(X_j)$.

\end{lem}

This lemma is a generalisation of the case that the category $\mathcal{A}$
is a directed $A_\infty$-category $\mathcal{A}(R)$ associated with a path
algebra with relations $R$ over a finite directed quiver.
By this lemma, we can compute an $A_\infty$-Koszul dual by iteration of twists.
We abbreviate $\mathcal{S}'(j)$ into $\mathcal{S}(j)$.
Together with the definition of $\mathcal{A}^!_\infty$ for directed $A_\infty$-category
$\mathcal{A}$, one has a natural isomorphism between $\mathcal{A}^!_\infty$
and $(Tw\mathcal{A})^\to(\boldsymbol{\mathcal{S}})$, where $\boldsymbol{\mathcal{S}}
= (\mathcal{S}(n), \mathcal{S}(n-1), \dots , \mathcal{S}(0))$.

We finish this section by recalling useful lemmas.

\begin{lem}[Lemma 5.6 in \cite{Se08}]\label{lem:FunctorialityOfTwist}

Suppose $\mathcal{F} \colon \mathcal{A} \to \mathcal{B}$ be a c-full and faithful $A_\infty$-functor and these  $Y_0$ and $Y_1$ be objects in $\mathcal{A}$.
Then, there exists a canonical isomorphism in $H^0(\mathcal{B})$ between
$\mathcal{F}(T_{Y_0}Y_1)$ and $T_{\mathcal{F}(Y_0)}\mathcal{F}(Y_1)$.

\end{lem}

\begin{lem}[Lemma 5.11 in \cite{Se08}]\label{lem:SphericalTwistIsEquiv}

Suppose that $Y_0$ is a spherical object in $\mathcal{A}$.
Then, $T_{Y_0}$ is a quasi-equivalence from $\mathcal{A}$ to itself.

\end{lem}

\begin{cor}\label{cor:TwistFormula}

Let $Y_0$ and $Y_1$ be objects in $\mathcal{A}$ and $Y_0$ is spherical.
Then, for any object $Z \in Ob(\mathcal{A})$, there exists a natural quasi-isomorphism between $T_{T_{Y_0}Y_1}Z$ and $T_{Y_0}T_{Y_1}T_{Y_0}^{-1}Z$.

\end{cor}

Here, the definition of a spherical objects can be found in (5h) in \cite{Se08}.

For a collection of objects $\boldsymbol{Y} = (Y_0, Y_1, \dots , Y_n)$ in
an $A_\infty$-category $\mathcal{A}$, we define a new collection $\mathcal{L}_j
\boldsymbol{Y} \coloneqq (Y_0, \dots, , Y_{j-1}, T_{Y_j}Y_{j+1}, Y_j, Y_{j+2},
\dots , Y_n)$ in $ \mathcal{A}$ and call it a \textit{mutation} of $\boldsymbol{Y}$.

\begin{lem}[Lemma 5.23 in \cite{Se08}]\label{lem:MutationAndDerEquiv}

Let $\boldsymbol{Y} = (Y_0, Y_1, \dots , Y_n)$ be a collection of spherical
objects in an $A_\infty$-category $\mathcal{A}$ and define $\boldsymbol{U} \coloneqq
\mathcal{L}_j \boldsymbol{Y}$.
Then there is a quasi-equivalence between $Tw \big( \mathcal{A}^\to(\boldsymbol{Y})
\big)$ and $Tw \big( \mathcal{A}^\to(\boldsymbol{U}) \big)$.
In particular, there is a equivalence of derived categories between $D(\mathcal{A}^\to(\boldsymbol{Y}))$
and $D(\mathcal{A}^\to(\boldsymbol{U}))$ as triangulated
categories.

\end{lem}

\begin{lem}\label{lem:SphericalTwistInAVSInATo}

Let $\boldsymbol{Y} = (Y_0, Y_1, \dots , Y_n)$ be a collection of spherical
objects in an $A_\infty$-category $\mathcal{A}$ and define $\boldsymbol{U} \coloneqq
\mathcal{L}_j \boldsymbol{Y}$.
Let us write $Y_i$ as an object in $\mathcal{A}^\to \coloneqq \mathcal{A}^\to(\boldsymbol{Y})$
as $\widetilde{Y}_i$ , the collection of them as $\widetilde{\boldsymbol{Y}}
= (\widetilde{Y}_1, \widetilde{Y}_2, \dots , \widetilde{Y}_n)$, and the mutation
as $\widetilde{U} \coloneqq \mathcal{L}_j \widetilde{\boldsymbol{Y}}$ (the
twist takes place in $\mathcal{A}^\to$, not in $\mathcal{A}$).
Then, there exists a quasi-isomorphism between $\mathcal{A}^\to(\boldsymbol{U})$
and $(\mathcal{A}^\to)^\to (\widetilde{\boldsymbol{U}})$.

\end{lem}

Conceptually, this lemma says that for spherical objects the twist in $\mathcal{A}$
and $\mathcal{A}^\to$ are equivalent in the above sense.
This lemma is proved in the proof of Lemma 5.23 in \cite{Se08}.

\section{Geometric preliminaries}\label{sec:GeomDef}

In this section, we prepare the notation of the Fukaya categories of exact
Riemann surfaces and discuss the twists in the Fukaya categories.

\subsection{Definition of the Fukaya categories}

The definition itself can be found in \cite{Se08} and its combinatorial description
which we mainly use can also be found in \cite{Su16}.
However, we repeat the relevant parts of those papers for the sake of completeness.

An \textit{exact symplectic manifold} $M = (M, \omega, \theta, J)$ consists
of a symplectic manifold with non-empty boundary $(M, \omega)$, a primitive
$\theta$ of $\omega$, i.e. $\theta$ is an 1-form satisfying $d\theta = \omega$,
and an $\omega$-compatible almost complex structure $J$. We impose that the negative Liouville vector field
$X_\theta$, defined by $\omega (-, X_\theta) = \theta(-)$, points strictly
inward on the boundary $\partial M$.

Now, we see the definition of the Fukaya category $\mathcal{F} = Fuk(M)$
of a given exact Riemann surface $M$.
In fact, we only use the Fukaya category $\mathcal{F}$ of the form $\mathcal{F}^\to(\boldsymbol{L})$
for some collection of objects $\boldsymbol{L}$ in this paper, hence what we really need to define is as follows: the set of objects $Ob(\mathcal{F})$,
the hom spaces $\hom_\mathcal{F}(L^\#_0, L^\#_1)$ for two distinct objects
$L^\#_0, L^\#_1 \in Ob(\mathcal{F})$, and the higher composition maps $\mu^d \colon \hom_\mathcal{F}(L^\#_{d-1},
L^\#_d) \otimes \hom_\mathcal{F}(L^\#_{d-2}, L^\#_{d-1}) \otimes \cdots \otimes
\hom_\mathcal{F}(L^\#_0, L^\#_1) \to \hom_\mathcal{F}(L^\#_0, L^\#_d)$ for
mutually distinct objects.

To define the objects of the Fukaya category $\mathcal{F} =  Fuk(M)$ of an
exact Riemann surface $M$, we fix a trivialization of $TM$ as a complex line
bundle (this is possible since $M$ possesses non-empty boundary).
Thanks to the complex structure $J$, we can identify the trivialization with
a non-vanishing vector field $X$.
Let $L \cong S^1 \hookrightarrow \mathring{M}$ be a Lagrangian submanifold.
We say that $L$ is \textit{exact} when $\int_L \theta = 0$. Let $\eta \colon
[0, 1] \to M$ be a compositon $[0, 1] \hookrightarrow \mathbb{R} \to \mathbb{R}/\mathbb{Z}
= S^1 \hookrightarrow M$ representing $L$.
We choose a function $\widetilde\alpha \colon [0, 1] \to \mathbb{R}$ such
that $\displaystyle \frac{d\eta}{dt}(t) \in \mathbb{R}_{>0}\cdot(e^{\pi i \widetilde\alpha
(t)}X_{\eta(t)}) \subset T_{\eta(t)}M$ holds.
Set $w(L) \coloneqq \widetilde{\alpha}(1) - \widetilde{\alpha}(0)$ and call
it the \textit{writhe} of $L$.
We say that $L$ is \textit{unobstructed} if $w(L) = 0$.
For an exact unobstructed Lagrangian submanifold $L$, we define its \textit{grading}
$\alpha \colon L \to \mathbb{R}$ by $\alpha (\eta(t)) = \widetilde\alpha (t)$.
We call a triple $L^\# = (L, \alpha , p)$ of an unobstructed Lagrangian
submanifold $L$, its grading $\alpha$, and arbitrary point $p \in L$ a \textit{Lagrangian
brane}.
Here, we call the third component of the Lagrangian brane $p$ a \textit{switching
point}.
Finally, we define the set of objects $Ob(\mathcal{F})$ of $\mathcal{F}$
by the set of all Lagrangian branes.
Note that a grading $\alpha$ of a Lagrangian brane $L^\#$ defines a new orientation
of $L$ by $(p \mapsto e^{\pi i \alpha (p)}X_p) \in \Gamma (TL)$. We call
it the \textit{brane orientation}.
We can see that a function $\alpha [n](p) \coloneqq \alpha (p) - n$ for $n
\in \mathbb{Z}$ is another grading of $L$. We call it the¥ $n$\textit{-fold
shift} of $\alpha$

\begin{rem}

There are few differences in the definition of objects of Fukaya categories
between Seidel \cite{Se08} and this paper.
In Seidel's definition, one uses a quadratic volume form $\eta_M^2$ which
is a section of $(\wedge^{top}T^*M)^{\otimes 2}$ (where the tensor product
is taken over $\mathbb{C)}$ while we use a non-vanishing vector field $X$.
The relation of these two is given by $\eta_M^2(X \otimes X) = 1$.
Then, our grading $\alpha$ is nothing but a grading $\alpha^\#$ of Seidel's
sense.
The relevant constructions are very simplified from Seidel's notation in this
paper since we only treat with exact Riemann surfaces (while Seidel considered
exact symplectic manifolds of any dimension).

A Lagrangian brane in Seidel's sense is a triple $L^\# = (L, \alpha^\# ,
P^\#)$, here $L$ and $\alpha^\#$ are the same with our definition but $P^\#$
is a Pin structure of $L$.
In the  definition of the Fukaya categories, the Pin structures are used for the determination
of the orientation of the moduli spaces. Hence they are used for the determination
of the sign of the higher composition maps $\mu$'s.
To determine the orientation of the moduli spaces, Seidel uses a real line bundle
$\beta$ associated to the Pin structure $P^\#$.
In the case of exact Riemann surfaces, the Pin structure $P^\#$ must be
non-trivial in order to achieve Theorem \ref{thm:AlgTwistIsDehnTwist} so we assume that. Therefore the real line bundle $\beta$
in our paper is always the non-trivial one.

Corresponding to that, a fixed point $p \in L$ is used as follows.
Since our real line bundle $\beta$ is not trivial, we can not trivialize
$\beta$ on whole $L$ but we can on $L \setminus \{ p \}$. With this trivialization,
we can consider that the orientation of $\beta$ changes when we go through
the point $p$. This is the meaning of $p$.
The choice of a point $p$ does not cause the difference of objects, i.e.
two Lagrangian branes $L^\#_0 = (L, \alpha , p_0)$ and $L^\#_1 = (L, \alpha
, p_1)$ are quasi-isomorphic in $\mathcal{F}$.
In the comparison with the definition by Seidel, we fix a trivialization
of a real line bundle $\beta$ instead of fixing of the Pin structure,
so we have as ``$S^1$ times many" objects as Seidel's definition. 

\end{rem}

Next, we define the hom set from $L^\#_0 = (L_0, \alpha _0, p_0)$ to $L^\#_1
= (L_1, \alpha _1, p_1)$.
From now on, we assume that any collection of Lagrangian branes
is in general position unless otherwise stated, i.e. any two submanifolds
intersect transversally, there is no triple point, and the switching point of one Lagrangian brane never contained in the other Lagrangian branes.
In this assumption, the hom set is defined by $\displaystyle \hom_\mathcal{F}(L^\#_0, L^\#_1)
\coloneqq \bigoplus _{p \in L_0\cap L_1} k \cdot [p]$ as a vector space.
Here, $[p]$ is a formal symbol corresponds to $p$. We sometimes abbreviate
$[p]$ into $p$.
For an intersection point $p \in L_0 \cap L_1$ as a morphism from $L^\#_0$ to $L^\#_1$, we define its index by $i(p)
= [\alpha _1(p) - \alpha _0(p)]+1$ and set $\displaystyle \hom_\mathcal{F}^d(L^\#_0,
L^\#_1) \coloneqq \bigoplus _{\substack{p \in L_0\cap L_1, \\ i(p) = d}}
k \cdot [p]$.

Finally, we define the $A_\infty$-structure $\mu$'s.
This is just a repetition but we only define the maps $\mu^d \colon \hom_\mathcal{F}(L^\#_{d-1},
L^\#_d) \otimes \hom_\mathcal{F}(L^\#_{d-2}, L^\#_{d-1}) \otimes \cdots \otimes
\hom_\mathcal{F}(L^\#_0, L^\#_1) \to \hom_\mathcal{F}(L^\#_0, L^\#_d)$ under
the condition that $L_i \pitchfork L_j$ and $p_i \not \in L_j$ for $i \neq
j$.

\begin{figure}[hbt]
\centering
\includegraphics[width=8cm]{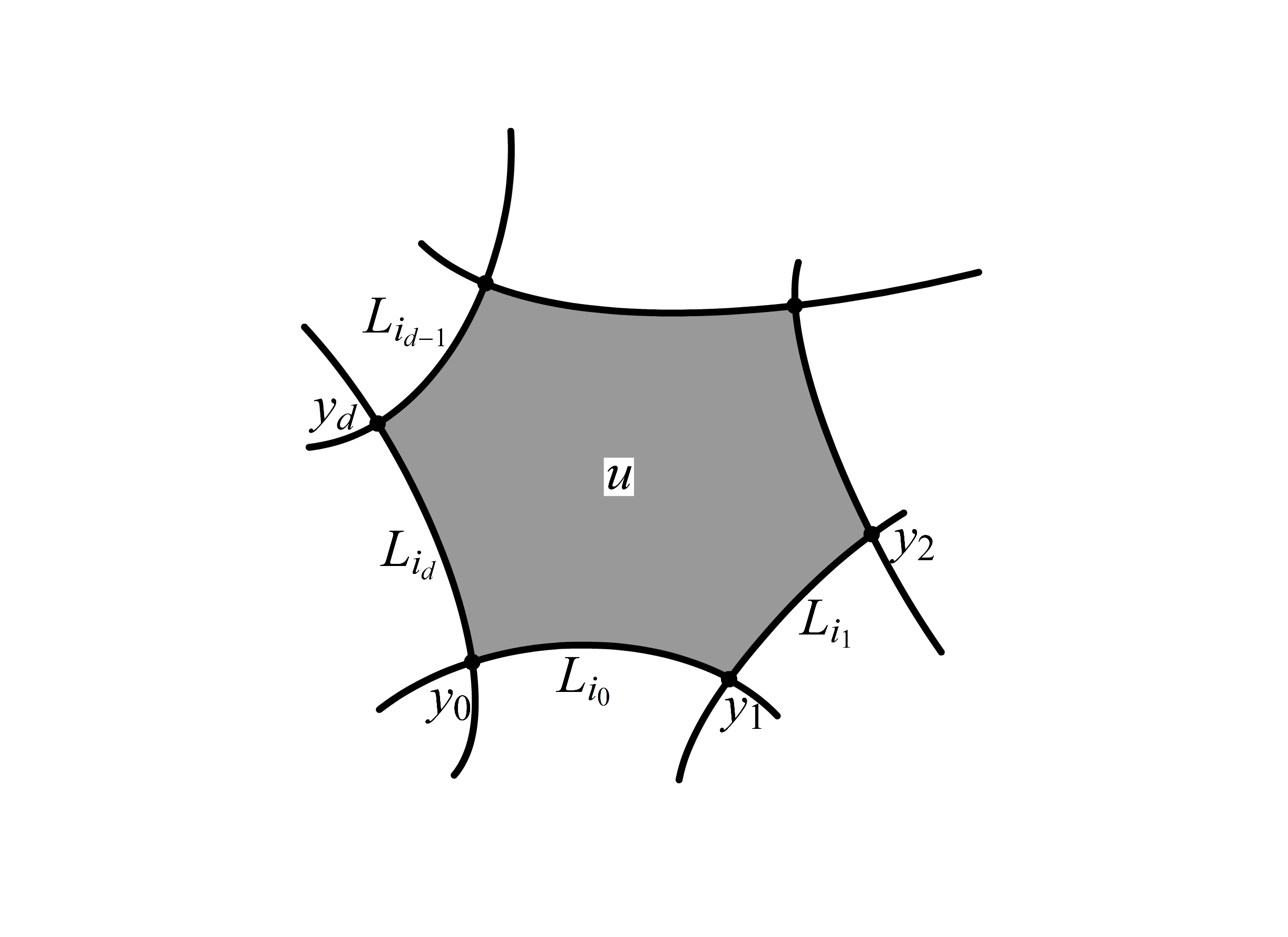}
\caption{(d+1)-gon in $M$}\label{fig:d+1gon}
\end{figure}

Let $\Delta ^{d+1}$ be a $(d+1)$-gon. We name its vertices $v_0, v_1, \dots
, v_d$ counterclockwise, the vertices connecting $v_{j-1}$ and $v_j$ by $[v_j,
v_{j+1}]$ $(0 \leq j < d)$, and the vertex connecting $v_0$ and $v_d$ by
$[v_d, v_0]$.
We define a set $\widetilde{\mathcal{M}}^{d+1}(y_1, y_2, \dots , y_d;
y_0)$ for $y_j \in L_{j-1} \cap L_j$ $(0 < j \leq d)$ and $y_0 \in L_0 \cap
L_d$ by the set of all orientation preserving immersions $u \colon \Delta
^{d+1} \to M$ satisfying that $u(v_j) = y_j$ for $0 \leq j \leq d$, $u([v_j,
v_{j+1}]) \subset L_j$ for $0 \leq j < d$, and $u([v_d, v_0]) \subset L_d$.
There exists a natural action of the group of diffeomorphisms of $\Delta^{d+1}$
fixing all vertices pointwise and the orientation.
Let us write the quotient space by $\mathcal{M}^{d+1}(y_1, y_2, \dots
, y_d; y_0)$ and call it a moduli space.
This moduli space is a set of $(d+1)$-gon as in Figure \ref{fig:d+1gon}.

It is known that this moduli space is empty unless $i(y_0) = i(y_1)
+ i(y_2) + \cdots + i(y_d) + (2-d)$. If it is not empty, then it is a finite
set.
For more detail, please refer Remark 3.22 in \cite{Su16} or Section 17 in
\cite{Se08}.

For an element $u \in \mathcal{M}^{d+1}(y_1, y_2, \dots , y_d; y_0)$,
we define its sign $(-1)^{s(u)} \in \{ \pm 1 \}$ as follows.
First, we assign $\pm 1$ to vertices. For a vertex $y_j$ of $u(\Delta^{d+1})$
with $0 < j \leq d$, we assign $-1$ to $y_j$ if the orientation of $L_j$ induced
from $\partial u$ and its brane orientation are opposite and $i(y_j)$ is
odd. Otherwise, we  assign $+1$ to it.
We assign $-1$ to $y_0$ if the orientation of $L_d$ induced from $\partial
u$ and its brane orientation are opposite and $i(y_0)$ is odd. Otherwise,
we assign $+1$ to it.
For each edge, $u([v_j, v_{j+1}])$ or $u([v_d, v_0])$, we assign $-1$ to it if it contains one of $\{ p_j\}$. Otherwise, we assign $+1$ to it.
Now $(-1)^{s(u)}$ is the product of all the $\pm 1$ above.
Finally, we define the map $\mu^d$ by 
\[ \mu^d(y_d, y_{d-1}, \dots , y_1) \coloneqq \sum_{\substack{y_0 \in L_0
\cap L_d\\ u \in \mathcal{M}^{d+1}(y_1, y_2, \dots , y_d; y_0)}} (-1)^{s(u)}
y_0. \]
It is known that this defines $A_\infty$-structure:

\begin{lem}

Let $L^\#_0, L^\#_1, \dots , L^\#_d$ be Lagrangian branes in general position.
For points $y_j \in L_{j-1} \cap L_j$ $(0 < j \leq d)$, we have the $A_\infty$-relation:
\begin{equation*}
\sum _{i, j, l} (-1)^{\bigstar _i} \mu ^l (y_d, \dots , y_{i+j+1}, \mu ^j
(y_{i+j}, \dots , y_{i+1}), y_i, \dots y_1) = 0.
\end{equation*}

\end{lem}

We don't prove this lemma. For the proof, please refer Theorem 3.25 in \cite{Su16}
or Section 12 and 13 in \cite{Se08}.

\begin{rem}

As we said, we don't define the hom set $\hom_\mathcal{F}(L^\#_0, L^\#_1)$
when two underlying spaces don't intersect transversally, especially the
case when $L^\#_0 = L^\#_1$. Also, we don't define  $\mu^d \colon \hom_\mathcal{F}(L^\#_{d-1},
L^\#_d) \otimes \hom_\mathcal{F}(L^\#_{d-2}, L^\#_{d-1}) \otimes \cdots \otimes
\hom_\mathcal{F}(L^\#_0, L^\#_1) \to \hom_\mathcal{F}(L^\#_0, L^\#_d)$ when
the objects are not in general position, especially in the case that
there exists $i \neq j$ such that $L^\#_i = L^\#_j$.
In these cases, we have to perturb the Lagrangian branes by Hamiltonian diffeomorphisms
and the definition itself becomes more complicated.
We can find the detail, in Chapter II of \cite{Se08} and many relevant papers and books for this topic.

\end{rem}

\begin{rem}\label{rem:AnyObjIsSpherical}

Any object in $\mathcal{F} = Fuk(M)$ for an exact Riemann surface $M$ is
spherical.
This is because the underlying space of any object is diffeomorphic
to the one dimensional sphere $S^1$ and some properties of Floer cohomology groups, for example, the
PSS isomorphism.

\end{rem}

\subsection{Algebraic twists versus geometric twists}

Let $V \cong S^1$ be an unobstructed exact Lagrangian submanifold of an exact
Riemann surface $M$.
Then the (right handed) Dehn twist $\tau_V$ can be lifted to a graded automorphism
of $M$. (The concept of graded automorphisms  appears in (12i) and the
existence of the lift is proved in the argument in (16f) of \cite{Se08}.)
Hence, for a Lagrangian brane $L^\#$, we can obtain a new Lagrangian brane
$\tau_V L^\# = \tau_V (L^\#)$.
If $V$ is an underlying space of a Lagrangian brane $V^\#$, then we write
the Dehn twist as $\tau_{V^\#}$.

\begin{thm}[(simplified version of) Theorem 17.16 in \cite{Se08}]\label{thm:AlgTwistIsDehnTwist}

Let $L^\#_0$ and $L^\#_1$ be two Lagrangian branes in an exact symplectic
manifold $M$.
Then, there exists an isomorphism between the Dehn twist $\tau_{L^\#_0}L^\#_1$
and the algebraic twist $T_{L^\#_0}L^\#_1$ in $H^0(TwFuk(M)) = DFuk(M)$.

\end{thm}

\begin{rem}

This theorem is essentially established in \cite{Se01}. This is  very fundamental
and crucial to define the Fukaya-Seidel categories of exact Lefschetz fibrations
which are studied with great attention in the context of homological mirror
symmetry.
See, for example, \cite{HV00}, \cite{Se00} for the mirror of $\mathbb{P}^2$,
and \cite{AKO08}.

\end{rem}

\section{Main results}\label{sec:MainResults}

In this section, we state the main theorem (Theorem \ref{thm:MAIN1TheFormulaOfAInftyKoszulDual}). The proofs will be presented in Section \ref{sec:GeomConstr}, \ref{sec:DirFukForRiemannDiag} and \ref{sec:CompOfDehnTwists}.
\subsection{Computation of $A_\infty$-Koszul duals}

\begin{thm}\label{thm:AInfKoszulDualByFukayaCat}

Let $\mathcal{A}$ be a directed $A_\infty$-category with the object set $Ob(\mathcal{A})
= \{ X_0 < X_1 < \cdots < X_n \}$. Suppose that there exist an exact Riemann
surface $M$ and a collection of Lagrangian branes $\boldsymbol{L}^\# = (L^\#_0,
L^\#_1, \dots , L^\#_n)$ such that $\mathcal{A}$ and $\mathcal{F}^\to(\boldsymbol{L}^\#)$
are quasi-isomorphic.
Then, $\mathcal{F}^\to(\boldsymbol{S}^\#)$ is an $A_\infty$-Koszul dual of
$\mathcal{A}$, 
where $\boldsymbol{S}^\# = (S^\#_n,
S^\#_{n-1}, \dots , S^\#_0)$ is a collection of objects defined by the iteration of Dehn twists $S^\#_j \coloneqq
\tau_{L^\#_0}\tau_{L^\#_1} \cdots \tau_{L^\#_{j-1}}L^\#_j$.

\end{thm}

\begin{proof}

What we have to show is that $\mathcal{A}^!_\infty$ and $\mathcal{F}^\to(\boldsymbol{S}^\#)$
are quasi-isomorphic.
First, $\mathcal{A}^!_\infty$ is naturally isomorphic to $(Tw\mathcal{A})^\to(\boldsymbol{\mathcal{S}})$,
where $\boldsymbol{\mathcal{S}} = (\mathcal{S}(n), \mathcal{S}(n-1), \dots
, \mathcal{S}(0))$ is a collection of objects $\mathcal{S}(j) = T_{X_0}T_{X_1}
\cdots T_{X_{j-i}}X_j$.
By assumption, there exists a quasi-isomorphism $\mathcal{G} \colon \mathcal{A}
\to \mathcal{F}^\to(\boldsymbol{L}^\#)$ and hence we have a quasi-isomorphism
$Tw\mathcal{G} \colon Tw\mathcal{A} \to Tw(\mathcal{F}^\to(\boldsymbol{L}^\#))$.
Both functors send $X_j$ to $L^\#_j$.
By Lemma \ref{lem:FunctorialityOfTwist}, $Tw\mathcal{G}(\mathcal{S}(j))$
is quasi-isomorphic to $\widetilde{S}^\#_j \coloneqq \widetilde{T}_{L^\#_0}\widetilde{T}_{L^\#_1}
\cdots \widetilde{T}_{L^\#_{j-1}} L^\#_j$, where $\widetilde{T}$ represents
the twist in $\mathcal{F}^\to(\boldsymbol{L})$.
Hence, $(Tw\mathcal{A})^\to(\boldsymbol{\mathcal{S}})$ is quasi-isomorphic
to $\big(Tw(\mathcal{F}^\to(\boldsymbol{L}^\#))\big)^\to(\widetilde{\boldsymbol{S}}^\#)$
by Lemma \ref{lem:FunctorialityOfDirSubcat}, \ref{lem:EquivOfDirSubcat},
and \ref{lem:FunctorialityOfTwist}, where $\widetilde{\boldsymbol{S}}^\#
= (\widetilde{S}^\#_n, \widetilde{S}^\#_{n-1}, \dots , \widetilde{S}^\#_0)$.

Now, by the iterated application of Lemma \ref{lem:SphericalTwistInAVSInATo}
(and Lemma \ref{lem:FunctorialityOfDirSubcat} and \ref{lem:FunctorialityOfTwist})
and Theorem \ref{thm:AlgTwistIsDehnTwist}, $\big(Tw(\mathcal{F}^\to(\boldsymbol{L}^\#))\big)^\to(\widetilde{\boldsymbol{S}}^\#)$
is quasi-isomorphic to $\mathcal{F}^\to(\boldsymbol{S}^\#)$.
This completes the proof.
\hfill $\Box$

\end{proof}

Thanks to the above theorem, we can compute an $A_\infty$-Koszul dual via the Fukaya categories and Dehn twists.
As in Example \ref{ex:DirAlgAsAInftyCat}, we can consider a path algebra
$R$ with relations over a finite directed quiver as a directed $A_\infty$-categorgy
$\mathcal{A} = \mathcal{A}(R)$.
The main point of this paper is to compute its $A_\infty$-Koszul dual by
this technique.

The followings are corollaries of Lemma \ref{lem:MutationAndDerEquiv}, Lemma \ref{lem:SphericalTwistInAVSInATo},
and Remark \ref{rem:AnyObjIsSpherical}.

\begin{cor}\label{cor:DerEquivOfKoszulDuals}

Let $\mathcal{A}$ and $\mathcal{F}^\to(\boldsymbol{S}^\#)$ be as in Theorem
\ref{thm:AInfKoszulDualByFukayaCat}.
Then, there exists a quasi-isomorphism between $Tw\mathcal{A}$ and $Tw(\mathcal{F}^\to(\boldsymbol{S}^\#))$,
hence there exists a equivalence of derived categories
between $D\mathcal{A}$ and $D(\mathcal{F}^\to(\boldsymbol{S}^\#))$ as triangulated categories.

\end{cor}

\begin{cor}\label{cor:KoszulDualIsDuality}

Let $\mathcal{A}$ and $\mathcal{F}^\to(\boldsymbol{S}^\#)$ be as in Theorem
\ref{thm:AInfKoszulDualByFukayaCat}.
Then, $\mathcal{A}$ is an $A_\infty$-Koszul dual of $\mathcal{F}^\to(\boldsymbol{S}^\#)$.

\end{cor}

\begin{proof}

Since $\tau_{S^\#_j}$ and $(\tau_{L^\#_1} \tau_{L^\#_2} \cdots \tau_{L^\#_{j-1}})
\tau_{L^\#_j} (\tau_{L^\#_1} \tau_{L^\#_2} \cdots \tau_{L^\#_{j-1}})^{-1}$
are isotopic, we choose a representative of $\tau_{S^\#_j}$ so that these
two coincide.
Then, we have $I_j \coloneqq \tau_{S^\#_n} \tau_{S^\#_{n-1}} \cdots \tau_{S^\#_{j+1}}
S^\#_j = \tau_{L^\#_0} \tau_{L^\#_1} \cdots \tau_{L^\#_n} L_j$.
Hence, we have a canonical isomorphism between $\mathcal{A} = \mathcal{F}^\to(\boldsymbol{L}^\#)$
and $\mathcal{F}^\to(\boldsymbol{I}^\#)$, where $\boldsymbol{I}^\# = (I^\#_1,
I^\#_2, \dots , I^\#_n)$.

\hfill $\Box$
\end{proof}

\subsection{Combinatorial setup}\label{subsec:CombinatorialSetup}

In this subsection, we prepare notations to describe $A_\infty$-Koszul duals of path algebras with relations with its Gabriel quiver is the directed
$A_n$-quiver $\ora{\Delta}_n$,

\[\xymatrix@R=0pt{
\bullet \ar[r]^{\alpha_1} & \bullet \ar[r]^{\alpha_2 } & \bullet \ar[r]^{\alpha_3}
& \cdots \ar[r]^{\alpha_n} & \bullet \, .\\
0 & 1 & 2 & \dots & n \, \, 
}\]

Let $R$ be a path algebra with relations $R = k\ora{\Delta} / (\rho_1, \rho_2
, \dots , \rho_m)$, where each relation $\rho_j = \alpha_{t_j} \alpha_{t_j-1}
\cdots \alpha_{s_j+1}$ is a path from $s_j$ to $t_j$ for $s_j, t_j \in [0,
n]_\mathbb{Z} \coloneqq \{ 0, 1, \dots , n\}$.
We call $s_j$ and $t_j$ a source point and a target point of $\rho_j$ respectively.
We assume that the length of any relation is greater or equal to two, i.e.
$t_j - s_j \geq 2$.
We call $\rho_j$ a relation corresponds to $[s_j, t_j]_\mathbb{Z}$.
Now, we can assume that $[s_i, t_i]_\mathbb{Z} \not \subset [s_j, t_j]_\mathbb{Z}$
for $i\neq j$ and $s_1 < s_2 < \cdots < s_m$ (hence $t_1 < t_2 < \cdots <
t_m$), so we assume them.
We write $S = \{ s_1, s_2, \dots s_m \}$, $T = \{ t_1, t_2, \dots , t_m \}$
and write $R = R_{S, T}$ to emphasise $S$ and $T$.
From now, we fix $n$, $S$, and $T$.

We define key items to describe an $A_\infty$-Koszul dual of $\mathcal{A}_{S,
T} \coloneqq \mathcal{A}(R_{S, T})$.
First, we define a map $d \colon [0, n]_\mathbb{Z} \to [0, n]_\mathbb{Z}
\sqcup \{ -\infty \}$ by $d(p) = \max \{ s_j \, | \, t_j \leq p \} = \max
\{ s \, | \, s < p, \hom_{\mathcal{A}_{S, T}}(s, p) = 0 \}$.
This is the nearest point $s$ smaller than $p$ such that $\hom_{\mathcal{A}_{S,
T}}(s, p)$ vanishes.
We define a finite decreasing sequence $\{ a^{(p)}_i \}_{0 \leq i \leq l_p}$
as follows.
First, we set $a^{(0)}_0 = 0$ and $l_0 = 0$.
For $p \geq 1$, we define $a^{(p)}_0 = p, a^{(p)}_1 = p-1$.
Suppose we have defined $a^{(p)}_q$ for $q < i$.
If $d(a^{(p)}_{i-2}) \neq d(a^{(p)}_{i-1})$, then we set $a^{(p)}_i \coloneqq
d(a^{(p)}_{i-2})$.
If $d(a^{(p)}_{i-2}) = d(a^{(p)}_{i-1})$, then we set $l_p = i-1$ and finish
the definition.

\begin{lem}

The sequences $\{ a^{(p)}_i\}_{0 \leq i \leq l_p}$ are strictly decreasing
and non-negative, i.e. $0 \leq a^{(p)}_i < a^{(p)}_{i-1}$.

\end{lem}

\begin{proof}

The inequality $a^{(p)}_i < a^{(p)}_{i-1}$ for $i = 1$ follows from the definition
of the sequence.
By definition of $d$, one can show $d(p) < p-1$ so we have the inequality for the case of $i =\ 2$.
Now, we prove the inequality $a^{(p)}_i < a^{(p)}_{i-1}$ for $i \geq
3$ (in the case when $i \leq l_p$).
Assume that $a^{(p)}_q < a^{(p)}_{q-1}$ holds for $q < i$.
By definition, we have $a^{(p)}_i = d(a^{(p)}_{i-2})$ and $a^{(p)}_{i-1}
= d(a^{(p)}_{i-3})$.
Since $d$ is non-decreasing, we have $a^{(p)}_i \leq a^{(p)}_{i-1}$. Moreover,
we have $a^{(p)}_i \neq a^{(p)}_{i-1}$ since $i \leq l_p$. Thus we have $a^{(p)}_i
< a^{(p)}_{i-1}$.

By definition, $a^{(p)}_i$ is an element in $d([0, n]_{\mathbb{Z}})$ such
that there exists $b \in d([0, n]_{\mathbb{Z}})$ satisfying $b < a^{(p)}_i$.
Therefore we have $a^{(p)}_i \neq -\infty$.
\hfill $\Box$

\end{proof}

\subsection{$A_\infty$-Koszul duals of path algebras}\label{subsec:FormulaOfAInfKoszulDual}

For $n, S, T$, we define a new directed $A_\infty$-category
$\mathcal{B}_{S,T}$ as follows.
Define $Ob(\mathcal{B}) = \{ B(n) < B(n-1) < \cdots < B(0) \}$, $\hom_{\mathcal{B}_{S,
T}}^i(B(p), B(a^{(p)}_i)) = k \cdot \eta_p^i$ (where $\eta_p^i$ is a formal
symbol), and other $\hom$'s are zero.
Let us write $\eta_p^l = \tilde{\eta}_{p, a^{(p)}_l}$. Then, we have that $\tilde{\eta}_{p, q}
\in \hom_{\mathcal{B}_{S, T}}(B(p), B(q))$.
Finally, we define $\mu$'s as follows:

\begin{align*}
\mu^d(\tilde{\eta}_{j_{d-1}, j_d} , \tilde{\eta}_{j_{d-2}, j_{d-1}}, \cdots
, \tilde{\eta}_{j_0, j_1}) =
\begin{cases}
(-1)^{\left( \big| \tilde{\eta}_{j_{d-1}, j_d} \big| + 1 \right)  \, \big|
\tilde{\eta}_{j_0, j_d} \big|} \tilde{\eta}_{j_0, j_d} & (\text{if it can
be non-zero}) \\
0 & (\text{otherwise}).
\end{cases}
\end{align*}

Here, ``it can be non-zero" means that $\hom_{\mathcal{B}_{S, T}}(B(j_0),
B(j_d))$ is non-zero and the relevant morphisms satisfy the degree condition
$\big| \tilde{\eta}_{j_0, j_d} \big| = \big| \tilde{\eta}_{j_0, j_1} \big| + \big| \tilde{\eta}_{j_1, j_2}
\big| + \cdots + \big| \tilde{\eta}_{j_{d-1}, j_d} \big| + (2-d)$, where $|x|$ stands
for the degree of $x$.
Then, this defines a directed $A_\infty$-category.

Now, the following theorem is the main theorem of this paper:

\begin{thm}\label{thm:MAIN1TheFormulaOfAInftyKoszulDual}

$\mathcal{B}_{S, T}$ is an $A_\infty$-Koszul dual of $\mathcal{A}_{S, T}$.

\end{thm}

\begin{cor}

An $A_\infty$-algebra $B_{S, T} \coloneqq \bigoplus _{i, j} \hom_{\mathcal{B}_{S,
T}}(B(i), B(j))$ is quasi-isomorphic to $(R_{S, T})^!_{dg}$.

\end{cor}

The proof is given in the following sections.
The outline is as follows: first, we find an exact Riemann surface $M$ and
a collection of Lagrangian branes $\boldsymbol{L}^\#$ such that $\mathcal{A}_{S,
T}$ and $\mathcal{F}^\to(\boldsymbol{L}^\#)$ are quasi-isomorphic; next, we compute the Dehn twists and obtain an $A_\infty$-Koszul
dual as $\mathcal{F}^\to(\boldsymbol{S}^\#)$.

Now, we study the structure of our $A_\infty$-Koszul dual $\mathcal{B}_{S,
T}$ with some concrete examples.
First, we study the case when $R = R_{S, T}$ is a quadratic algebra, i.e.
all the relations are of the form $[i, i+2]_\mathbb{Z}$.
By easy calculation, we can show that $\hom_\mathcal{B}(B(p), B(q))$ for $p \geq q$ is non-zero
only when $\{ q, q+1, \cdots , p-2 \} \in S$.
Moreover, when this is the case, the degree of the non-zero morphism is $p-q$.
By the condition of degree, we can show that $\mu^d = 0$ except for $d=2$.
Finally, we can conclude that $\mathcal{B}$ is isomorphic to $\mathcal{A}((R_{S^c,
T^c})^{op})$, where $S^c \coloneqq \{ 0, 1, \dots, n-2\} \setminus S$ and
$T^c \coloneqq \{ 2, 3, \dots, n\} \setminus T$. 

For example, if $\xi_j  (\neq 0) \in \hom_\mathcal{A}(j, j+1)$ and $\xi_{j+1}
(\neq 0) \in \hom_\mathcal{A}(j+1, j+2)$ satisfy that $\mu^2(\xi_{j+1}, \xi_j) = \mu_\mathcal{A}^2(\xi_{j+1}, \xi_j) \neq 0$, then we have $\mu^2(\tilde{\eta}_{j,
j+1}, \tilde{\eta}_{j+1, j+2})
= \mu_\mathcal{B}^2(\tilde{\eta}_{j, j+1}, \tilde{\eta}_{j+1, j+2})
= 0$.
Conversely,  if $\xi_j  (\neq 0) \in \hom_\mathcal{A}(j, j+1)$ and $\xi_{j+1}
(\neq 0) \in \hom_\mathcal{A}(j+1, j+2)$ satisfy that $\mu^2(\xi_{j+1}, \xi_j)
= \mu_\mathcal{A}^2(\xi_{j+1}, \xi_j) = 0$, then we have $\mu^2(\tilde{\eta}_{j, j+1}, \tilde{\eta}_{j+1, j+2})
= \mu_\mathcal{B}^2(\tilde{\eta}_{j, j+1}, \tilde{\eta}_{j+1, j+2})
\neq 0$.
Thus, we can observe that the products and relations of $\mathcal{A}$ and
$\mathcal{B}$ are ``reversed" as we have already saw.

Next, we see the case that $n = 3$ and we have only one relation corresponding to $[0, 3]_\mathbb{Z}$.
The algebra $R = R_{S, T}$ is no longer quadratic.
For this algebra, the duality emerges as the following form.
In this case, the formula defines $\mathcal{B}$ as follows:
hom spaces are all zero but $\hom_\mathcal{B}^0(B(j), B(j)) = k \cdot \tilde{\eta}_{j, j}$,
$\hom_\mathcal{B}^1(B(j+1), B(j)) = k \cdot \tilde{\eta}_{j, j+1}$, $\hom_\mathcal{B}^2(B(3),
B(0)) = \tilde{\eta}_{0, 3}$; $\mu$'s are all zero but $\mu^3(\tilde{\eta}_{0, 1},
\tilde{\eta}_{1, 2}, \tilde{\eta}_{2, 3}) = \tilde{\eta}_{0, 3}$.
This is nothing but the duality between product and relation.
This phenomenon cannot be captured in the dg settings because the dg-structure
lacks the structure of higher composition maps. 

Let us see the general cases of $R_{S, T}$. We can show  $l_p \geq 2 \Leftrightarrow p \in
T$ and when this is the case, there exists a relation corresponding to $[a^{(p)}_2, p]_\mathbb{Z}$  (we show this later
but this is not so hard). 
 We can see that the relation corresponds to $[s_j, t_j]_\mathbb{Z}$ in $R_{S, T}$ emerges in the structure
of $\mathcal{B}$ as the degree two morphism $\tilde{\eta}_{s_j, t_j}$ with
nontrivial higher composition $\mu^{t_j - s_j}$.

These are the typical examples:

\begin{cor}\label{cor:typicalEx}

Define $S_{n, k}$ and $T_{n, k}$ for $n >\ k$ by $S_{n, k} = \{ 0, 1, \cdots , n-k \}$ and $\{ k, k+1, \cdots , n \}$. Let us write $A_{n, k} \coloneqq
A_{S_{n, k}, T_{n, k}}$,
$\mathcal{A}_{n, k} \coloneqq \mathcal{A}_{S_{n, k}, T_{n, k}}$, and
$\mathcal{B}_{n, k} \coloneqq \mathcal{B}_{S_{n, k}, T_{n, k}}$.
Then, we have the following:
\begin{enumerate}
\item For $\mathcal{B}_{n, k}$, we have
\begin{align*}
\hom_{\mathcal{B}_{n, k}}^d(B(p), B(q)) = 
\begin{cases}
k \cdot \tilde{\eta}_{p, q} & (d \geq 0 \text{ is even and } p - q = kl \text{ for } d = 2l) \\
k \cdot \tilde{\eta}_{p, q} & (d \geq 0 \text{ is odd and } p - q = kl + 1 \text{ for } d = 2l + 1) \\
0 & (\text{otherwise})
\end{cases} 
\end{align*}
and $\mu^k(\eta^1_p, \eta^1_{p-1}, \dots , \eta^1_{p-k+1}) = \eta^2_p \colon S^\#_p \to S^\#_{p-k}[2]$.
(There are many other collections of morphisms with non-vanishing higher compositions, but we omit to write.)
\item Especially, for $n = k$, our category $\mathcal{B}_n \coloneqq \mathcal{B}_{n, n}$ is described as follows: $Ob(\mathcal{B}_n) = \{ B(n) > B(n-1) > \cdots > B(0) \}$, 
\begin{align*}
\hom_{\mathcal{B}_{n, k}}^d(B(p), B(q)) = 
\begin{cases}
k \cdot \eta^0_p & (d = 0, p = q) \\
k \cdot \eta^1_p & (d = 1, p - q = 1) \\
k \cdot \eta^2_n & (d = 2, p = n, q = 0) \\
0 & (\text{otherwise}),
\end{cases} 
\end{align*}
 and $\mu$'s are all zero but $\mu^2$ with identity morphisms and $\mu^n(\eta^1_n, \eta^1_{n-1}, \dots , \eta^1_1) = \eta^2_n \colon S^\#_n \to S^\#_0[2]$.
\end{enumerate}
\end{cor}

It is remarkable that the whole information of relations of $R_{S, T}$ can
be recovered (by hand) by the morphisms of $\mathcal{B}$ with degree less
than or equal to two and relevant higher compositions.
Thus, there emerges a natural question.

\begin{prob}

Find the properties of directed $A_\infty$-categories that determines $\mathcal{B}$ from its objects, morphisms with degree less than or equals to 2, and $\mu$'s between such a morphisms.

\end{prob}

\subsection{Combinatorial lemmas}

We prepare two lemmas which are used in the geometric computations in Section\ref{sec:CompOfDehnTwists}.

First, we define $d^\dagger \coloneqq [0, n]_\mathbb{Z} \to [0, n]_\mathbb{Z}
\sqcup \{ \infty \}$ to be the dual of $d$ by $d^\dagger(p) = \min \{ t_j \, | \, s_j \geq p \} = \min \{ s \, | \, s >
p \text{ and } \hom_\mathcal{A}(p, s) = 0 \}$, where $\mathcal{A} = \mathcal{A}_{S,
T} = \mathcal{A}(R_{S, T})$.
Next, we define the sequence $\{ a^{(p) \dagger}_j \}_{0 \leq j \leq l^\dagger_p}$
by replacing $d$ by $d^\dagger$ and setting $a^{(p)\dagger}_1 = p+1$ in the
definition of $\{ a^{(p)}_j \}_{0 \leq j \leq l_p}$.

\begin{lem}[Inversion formula]\label{lem:InversionFormula}

The sequences satisfy the following inversion formula: $a^{(a^{(p)}_j) \dagger}_j
= p$ and $a^{(a^{(p)\dagger}_j)}_j = p$ for $0 \leq j \leq l_p$ and $0 \leq j \leq l^\dagger_p$ respectively.

\end{lem}

\begin{proof}

First, we prove the former formula. We write $q = a^{(p)}_j$.
Since the statement for the case of $j = 0, 1$ is trivial, we consider the case $j \geq 2$ (so we are assuming that $l_p \geq j (\geq 2)$).
By the definition, there exists $t \in T$ such that $a^{(p)}_{i+1} < t \leq a^{(p)}_i$ for $0 \leq i \leq j-2$ since $d(a^{(p)}_{i+1}) \neq d(a^{(p)}_i)$.
For each $i$, we write the max of such $t$'s as $t_i$, i.e. $t_i \in T$,  $a^{(p)}_{i+1} < t_i \leq a^{(p)}_i$, and there is no element $t \in T$ such that $t_i < t \leq a^{(p)}_i$.

\begin{claim}\label{claim:inversionFormula}
$l^\dagger_q \geq j$ and $a^{(p)}_{j-l+1} < a^{(q)\dagger}_l \leq t_{j-l}$ for $2 \leq l \leq j.$
\end{claim}

We prove $l^\dagger_q \geq l$ and the above inequality by induction on $l$.
Let us consider the case of $l = 2$.
The inequality $l^\dagger_q \geq 2$ holds because of $q = d(a^{(p)}_{j-2}) $.
Also, the second inequality holds since $a^{(p)}_{j-1} < t_{j-2} = d^\dagger(q) = a^{(q)\dagger}_2$.

Next, we consider the case of $l = 3$.
Since $a^{(p)}_j < a^{(p)}_{j-1}$, we have $a^{(q)\dagger}_1 = q+1 \leq a^{(p)}_{j-1}$.
Since $a^{(p)}_{j-1} < a^{(p)}_{j-2}$ and $a^{(p)}_{j-1} = d(a^{(p)}_{j-3})$, there exists $s \in S$ such that $a^{(q)\dagger}_1 \leq s < a^{(q)\dagger}_2$ (one example of such an $s$ is $a^{(p)}_{j-1}$).
Thus, we have $l^\dagger_q \geq 3$.
Since $d^\dagger$ is non-decreasing, we have $a^{(q)\dagger}_3 = d^\dagger(a^{(q)\dagger}_1) \leq d^\dagger(a^{(p)}_{j-1}) = t_{j-3}$.
By the definition, there exists a relation corresponds to $[c, a^{(q)\dagger}_3]_\mathbb{Z}$ such that $c \in S$ is the smallest element in $S$ greater than or equals to $a^{(q)\dagger}_1 = q + 1$.
Together with $a^{(q)\dagger}_3 \leq t_{j-3}$, we have $q < c \leq d(t_{j-3}) = a^{(p)}_{j-1} \, (\in S)$.
Thus, we have $t_{j-2} < a^{(q)\dagger}_3 \leq t_{j-3}$.
Now, since $t_{j-2}$ is the lergest element in $T$ less than or equals to $a^{(p)}_{j-2}$, so we have $a^{(p)}_{j-2} < a^{(q)\dagger}_3 \leq t_{j-3}$.

Let us assume that $l^\dagger_q \geq s$ and the latter inequality in the Claim is true for $l \leq s$ for some $s$ with $3 \leq s < j$.
Then, we have the following inequality:

\[ a^{(p)}_{j-s+2} < a^{(q)\dagger}_{s-1} \leq t_{j-s+1} \leq a^{(p)}_{j-s+1} < a^{(q)\dagger}_s \leq t_{j-s} \leq a^{(p)}_{j-s} < t_{j-(s+1)}. \]

(Here, first and second inequality follows from the case of $l = s - 1$, third, sixth, and seventh inequality follows from the definition of $t$'s, fourth and fifth inequality follows from the case of $l = s$.)
Since $a^{(p)}_{j-s+1} \in S$ lies in $[a^{(q)\dagger}_{s-1}, a^{(q)\dagger}_s - 1]_\mathbb{Z}$, we have $l^\dagger_q \geq s + 1$.
Now, the inequality $a^{(p)}_{j-s} < a^{(q)\dagger}_{s+1} \leq t_{j-(s+1)}$ follows from applying $d^\dagger$ on $a^{(p)}_{j-s+2} < a^{(q)\dagger}_{s-1} \leq a^{(p)}_{j-s+1}$ and the maximality of $t_{j-s}$.
This completes the proof of Claim \ref{claim:inversionFormula}.

By substituting $l = j$ into the inequality of Claim \ref{claim:inversionFormula}, we have that $p-1 = a^{(p)}_1 < a^{(q)\dagger}_j \leq a^{(p)}_0 = p$. Thus we have $a^{(q)\dagger}_j = p$.
 
The latter formula of the Lemma \ref{lem:InversionFormula} can be proven by the argument obtained by interchanging symbols with $\dagger$ and without $\dagger$.
\hfill $\Box$
\end{proof}

\begin{lem}\label{lem:dddaggerlemma}

$d(d^\dagger(p)-1) \leq p-1$.

\end{lem}

\begin{proof}

We can write $d^\dagger(p) = t_j$ with $j = \min \{ i \, | \, p \leq s_i
\}$.
If $d(d^\dagger(p)-1) \geq p$, there exists $j'$ such that $p \leq s_{j'}
< t_{j'} \leq d^\dagger(p)-1$.
Since $t_{j'} \leq d^\dagger(p)-1 < t_j$, we have $s_{j'} < s_j$.
However, this contradicts with the minimality of $j$.
Thus we have the conclusion.

\hfill $\Box$
\end{proof}

\section{Construction of Riemann surfaces and Lagrangian branes}\label{sec:GeomConstr}

Our goal in this section is to construct an exact Riemann surface $M$ and
a collection of Lagrangian branes $\boldsymbol{L}^\# = (L^\#_1, L^\#_2, \dots
, L^\#_n)$ such that $\mathcal{A}_{S, T}$ and $\mathcal{F}^\to(\boldsymbol{L}^\#)$
are isomorphic. 

We use many small positive $\varepsilon$'s. We assume that they all are small
enough. We change such $\varepsilon$'s smaller if
necessary without any notification to avoid unnecessary complexity and confusion.

\subsection{Lemmas for construction}

First, we prepare some notations. 
Let $M$ be a two-dimensional manifold with non-empty boundary.
A \textit{two-tailed Lagrangian submanifold} $L^\pm = (L, \gamma^+, \gamma^-)$
is a triple of a one-dimensional submanifold $L$ diffeomorphic to $S^1$ and
tails $\gamma^\pm$.
Here, these tails are embeddings $\gamma^\pm \colon [0, 1] \to M$ such that
$\gamma^\pm(0) \in L$, $\gamma^\pm(1) \in \partial M$, $\gamma^\pm((0, 1))
\cap (L\ \cup \partial M) = \varnothing$, $\gamma^\pm \pitchfork L$ at $\gamma^\pm(0)$, the orientations of $\gamma^+$
and $L$ defines positive orientation of $M$ at $\gamma^+(0),$ and the other pair $\gamma^-$
and $L$ defines negative orientation of $M$ at $\gamma^-(0)$.
A collection of two-tailed Lagrangian submanifolds $\boldsymbol{L}^\pm
= (L^\pm_1, L^\pm_2, \dots , L^\pm_n)$ is called \textit{compatible} when
$\gamma^\pm_i([0, 1]) \cap \gamma^\pm_j([0, 1]) = \varnothing$ and $\gamma^\pm_i([0,
1]) \cap L_j = \varnothing$ for $i \neq j$.

\begin{lem}\label{lem:ConstructExactStr}

For a compatible collection of two-tailed Lagrangian submanifolds $\boldsymbol{L}^\pm$,
there exists an exact symplectic structure $(\omega, \theta, J)$ on $M$ such
that the underlying Lagrangian submanifolds of $\boldsymbol{L}^\pm$ become
exact Lagrangian submanifolds.

\end{lem}

\begin{proof}

First, we take an arbitrary exact symplectic structure $(\omega_0, \theta_0, J)$.
We write the boundary as $\partial M = \bigsqcup S_j$ and fix their collar
neighbourhoods $\iota_j \colon S^1 \times [0, \varepsilon_1) \hookrightarrow
M$.
We can assume that $\iota_j^* \omega_0 = d\varphi \wedge dx$ and $\iota_j^*
\theta_0 = (B_j-x)d\varphi$ where $\varphi$ and $x$ is the natural coordinate
of $S^1 =\ \mathbb{R} / \mathbb{Z}$ and $[0, \varepsilon_1)$ respectively
and $B_j = \int_{S_j} \theta_0 \,$ (relevant argument can be found in the
proof of Lemma 5.4 in \cite{Su16}).

With this setting, we change $\omega_0$ and $\theta_0$ as follows.
First, let us assume $E_1 \coloneqq \int_{L_1} \theta_0 < 0$.
Choose a function $h \colon (-\varepsilon_2, \varepsilon_2) \to \mathbb{R}_{\geq
0}$ such that $h$ is compactly supported and $\int _{-\varepsilon_2}^{\varepsilon_2}
h(t) dt = 1$.
Construct a new manifold with exact symplectic form $(\widetilde{M}, \widetilde{\omega_0},
\widetilde{\theta_0})$ as follows. The new manifold $\widetilde{M}$ is constructed as a gluing
of $V_1 \coloneqq \{ (\phi, x) \in (-\varepsilon_2, \varepsilon_2) \times
\mathbb{R} \, | \,  E_1h(\phi) \leq x <\ \varepsilon_1\ \, \}$ and $M$. Here,
the gluing identifies  $(0, 0) \in V_1$ with $\gamma_1^+(1)$ and identifies
other points naturally with respect to the tublar neighbourhood. See Figure
\ref{fig:widetildeM}.

\begin{figure}[hbt]
\centering
\includegraphics[width=8cm]{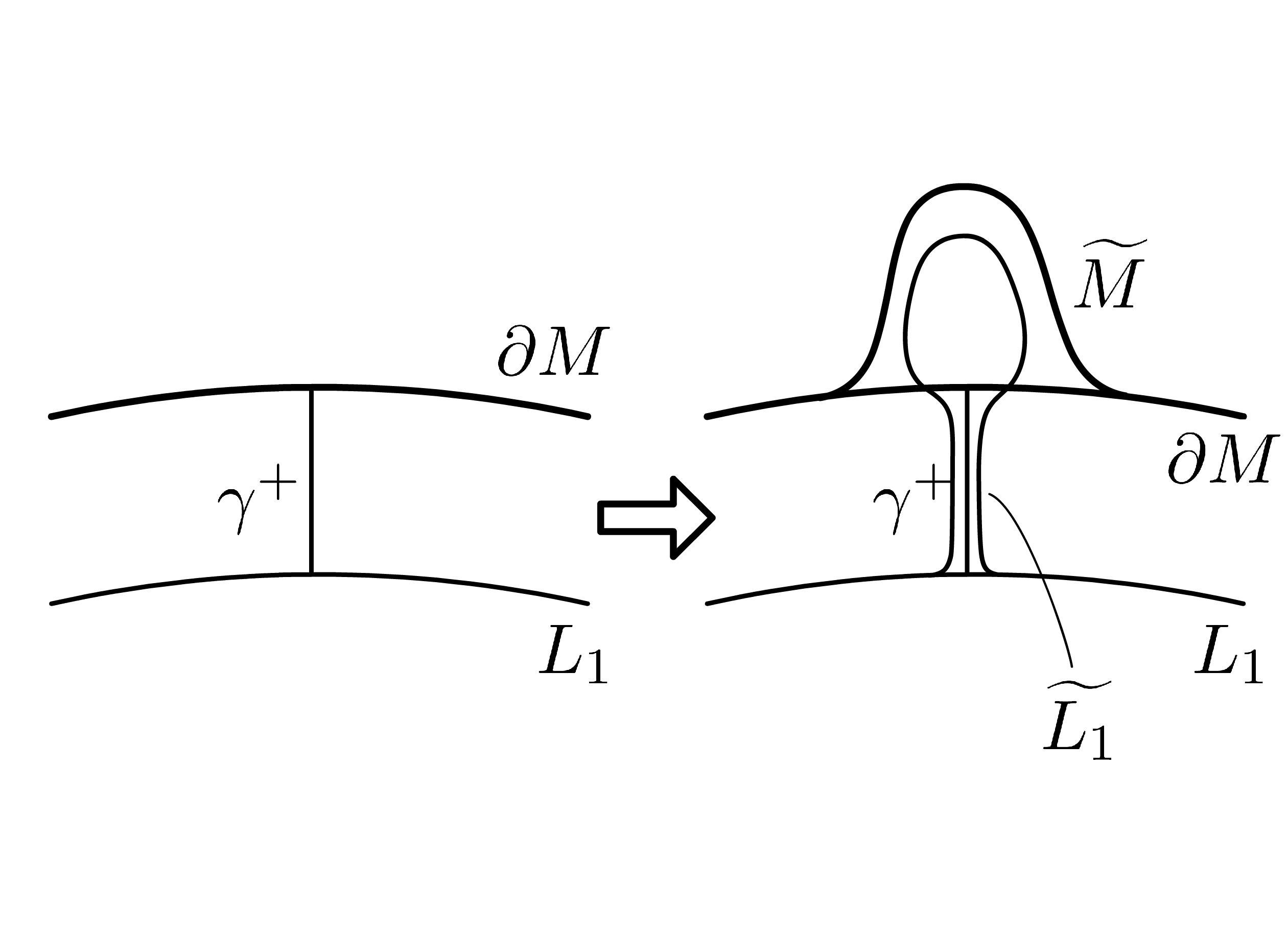}
\caption{$\widetilde{M}$ and $\widetilde{L}$'s}\label{fig:widetildeM}
\end{figure}

We set $\widetilde{\omega_0}|_{V_1} = d\phi \wedge dx$ and $\widetilde{\theta_0}
= (B_1 - x)d\phi$.
Then, the negative Liouville vector field can be written as $\displaystyle
X_\theta = (B_1 - x) \frac{\partial}{\partial x}$ hence it points strictly
inwards on $\partial \widetilde{M} \cap V_1$.

Next, we set submanifolds $\widetilde{L}_j$ in $\widetilde{M}$ to be the image
of the natural embedding $\iota \colon M \hookrightarrow \widetilde{M}$ for
$j \geq 2$.
For the case of $\widetilde{L}_1$, we define it as in Figure \ref{fig:widetildeM} such that
$\widetilde{L}_1$ is a deformation of $\iota(L_1)$ supported on very small
region aroung $\iota(\gamma^+(0))$ and $V_1$ so that the deformation does not cerate new intersection points of $\widetilde{L}$'s
and $\int_{\widetilde{L}_1} \widetilde{\theta_0} = 0$.

Now, there is a diffeomorphism $f \colon M \to \widetilde{M}$ such that $f$
coincides with the embedding away from small neighbourhood of the region
surrounded by $\iota(L_{1})$ and $\widetilde{L}_1$ and $f(L_j) = \widetilde{L}_j$.
We set $\omega_1 = f^*\widetilde{\omega_0}$ and $\theta_1 = \widetilde{\theta_0}$.
Then we can show that $(M, \omega_1, \theta_1, J)$ is an exact Riemann surface,
$\int_{L_1}\theta_1 = 0$, and $\int_{L_j} \theta_1 = \int_{L_j} \theta_0$
for $j \geq 2$.
Even when $E_1 = \int_{L_1} \theta_0 > 0$, we can find such $\omega_1$ and $\theta_1$  by almost the same construction which involves $\gamma^-_1$.
We iterate this construction and we can obtain the desired exact symplectic
structure. 
\hfill $\Box$

\end{proof}

\begin{lem}\label{lem:ConstructGrading}

Let $M$ be an exact Riemann surface and $\boldsymbol{L} = (L_1, L_2, \dots
, L_n)$ be a collection of exact Lagrangian submanifolds.
Suppose that $[L_1], [L_2], \dots , [L_n] \in H_1(M; \mathbb{Z})$ are linearly
independent and $H_1(M; \mathbb{Z}) / ([L_1], [L_2], \dots , [L_n])$ does
not have torsion other than two-torsion, then $TM$ admits a trivialization
$X \in \Gamma(TM)$ such that all the Lagrangian submanifolds $L_j$'s are unobstructed.

\end{lem}

\begin{proof}

Let us consider an exact Lefschetz fibtration $\pi \colon E \to D$ in the
sense of \cite{Se08} with vanishing cycles $L_1, L_2, \dots , L_n$ with a
suitable distinguish basis of vanishing paths.
Here, the taeget space $D$ is the unit disc in $\mathbb{C}$.
Then, the total space $E$ has the homotopy type of a two-dimensional CW-complex
which is obtained by attaching $n$ disks to $M$ along the vanishing cycles
\cite{Ka80}.
By the computation of the Mayer-Vietoris exact sequence, the assumptions on
homology classes of $L_j$'s induce that $H^2(E; \mathbb{Z}) \cong (\mathbb{Z}/2)^{\oplus
p}$ for some $p$.
Hence, we have $2c_1(E) = 0$.

Since the two-fold first Chern class $2c_1(E)$ vanishes, the total space
$E$ admits relative quadratic volume form $\eta_{E/D}^2$ as in the discussion
in (15c) of \cite{Se08}.
Hence, when we use the induced trivialization, all vanishing
cycles are unobstructed by the discussion in (16f) of \cite{Se08}.
\hfill $\Box$

\end{proof}

When a compatible collection of two-tailed Lagrangian submanifolds satisfies
the homological condition in Lemma \ref{lem:ConstructGrading}, we call such
a collection a \textit{perfect} collection of two-tailed Lagrangian submanifolds.
For a perfect collection of two-tailed Lagrangian submanifolds $\boldsymbol{L}^\pm
= (L_1^\pm , \dots , L^\pm_n)$, where $L^\pm_j = (L_j, \gamma^+_j, \gamma^-_j)$,
we can construct an exact symplectic structure and brane structures of each
submanifold of $\boldsymbol{L}^\pm$, namely $L^\#_j = (L_j, \alpha_j, p_j)$,
by the above lemmas.
We call the tuple $L^\maltese_j = (L_j, \gamma^+_j, \gamma^-_j, \alpha_j,
p_j)$ a two-tailed Lagrangian brane. The resulting two-tailed Lagrangian
branes form a collection of two-tailed Lagrangian branes $\boldsymbol{L}^\maltese$.
We define a directed $A_\infty$-category $\mathcal{F}^\to(\boldsymbol{L}^\maltese)$
by $\mathcal{F}^\to(\boldsymbol{L}^\maltese) \coloneqq \mathcal{F}^\to(\boldsymbol{L}^\#)$.

\subsection{Construction (1)}

In this and the next subsection, we construct an exact symplectic manifold
$M$ and a collection of two-tailed Lagrangian branes $\boldsymbol{L}^\maltese
= (L^\maltese_0, L^\maltese_1, \dots L^\maltese_n)$ so that $\mathcal{F}^\to(\boldsymbol{L}^\maltese)$
is isomorphic to $\mathcal{A}_{S, T}$.
In this subsection, we construct them for the case that $S = T = \varnothing$ as
a prototype of all the construction.

Thanks to the previous two lemmas, what we have to construct is reduced to
a two-dimensional manifold $M$ with non-empty boundary, a perfect collection
of two-tailed Lagrangian submanifolds $\boldsymbol{L}^\pm$, and brane
structures
on the underlying spaces of the two-tailed Lagrangian submanifolds.

We set $C_j \coloneqq S^1 = \mathbb{R} / \mathbb{Z}$ for $0 \leq j \leq n$
and define $M_0$ by the plumbing (and smoothing) of $D_j \coloneqq C_j \times
[-\varepsilon_3, \varepsilon_3]$.
Namely, let $\varphi_j$ and $x_j$ be the natural coordinate of $C_j$ and $[-\varepsilon_3,
\varepsilon_3]$.
Our plumbing is defined to identify two points $(\varphi_{j+1}, x_{j+1})
\in \{ (\varphi_{j+1}, x_{j+1}) \, | \, -\varepsilon_3 < \varphi_{j+1} <
\varepsilon_3\} \subset D_{j+1}$ and $\displaystyle \left( \frac12 + x_{j+1},
-\varphi_{j+1} \right) \in D_j$ for every $0 \leq j < n$.

\begin{figure}[hbt]
\centering
\includegraphics[width=8cm]{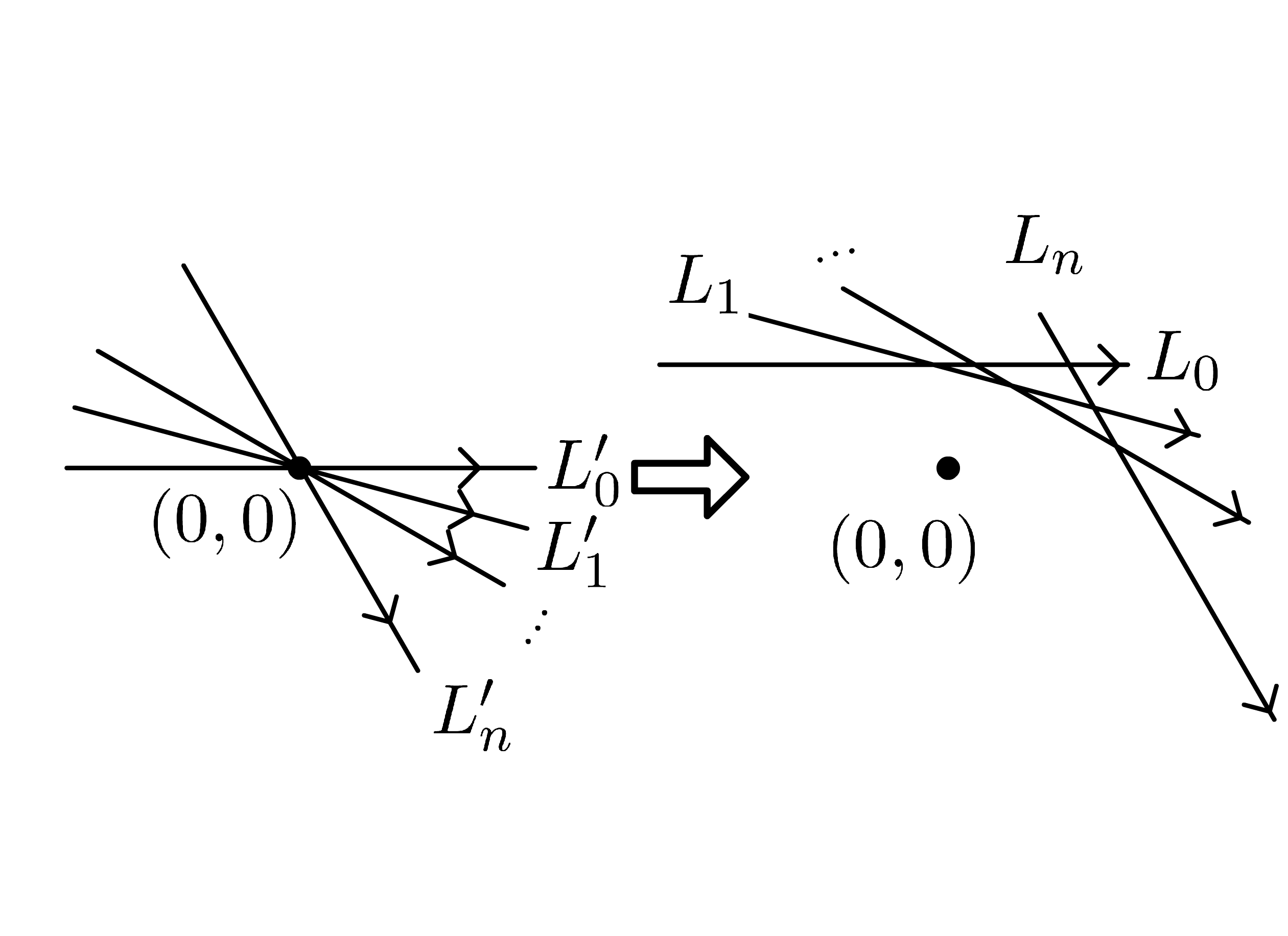}
\caption{$L'$'s and $L$'s}\label{fig:L''s}
\end{figure}

Fix smooth (right handed) Dehn twists $\tau_{C_j}$ along $C_j$ supported
in $D_j$.
We assume that $\tau_{C_j}|_{C_j}$ is the antipodal map.
We define submanifolds by $L'_j \coloneqq \tau_{C_1}^{-1} \tau_{C_2}^{-1}
\cdots \tau_{C_{j-1}}^{-1} C_j$.
Then, all these submanifolds pass through $(0, 0) \in D_0 \subset M_0$ and
there is no other intersection point.
We deform them to avoid $(0, 0)$ and pass through the left side of the point
with respect to their orientation as in Figure \ref{fig:L''s} and obtain
the resulting submanifolds $L_0, L_1, \dots , L_n$.

\begin{figure}[hbt]
\centering
\includegraphics[width=8cm]{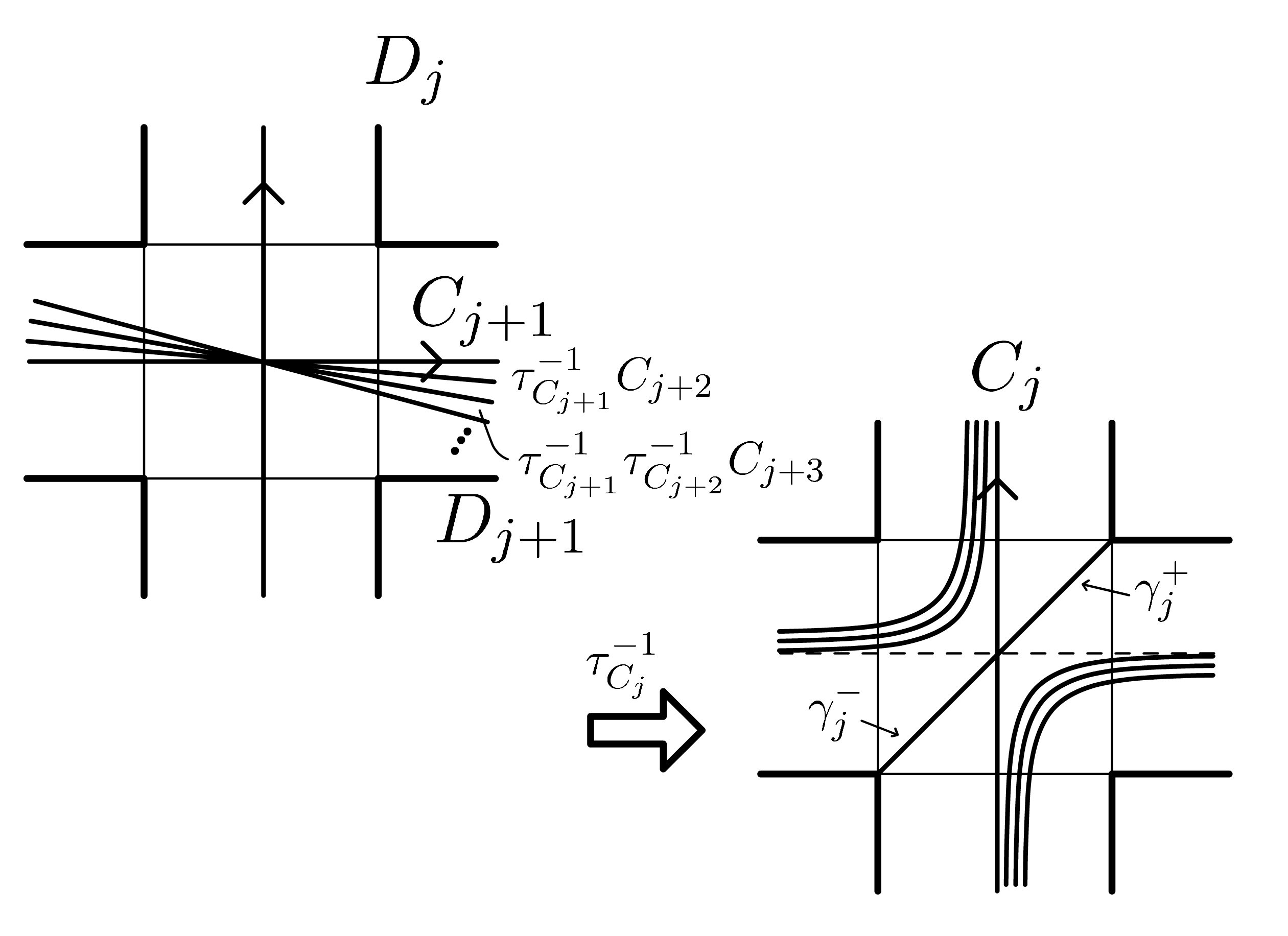}

\caption{$D_j \cap D_{j+1}$}\label{fig:DjAndDj+1}
\end{figure}

Next, we give tails to $L_j$'s such that they form a compatible collection
of two-tailed Lagrangian submanifolds. 
Let us consider the following sequence of collections of submanifolds $\boldsymbol{\Gamma}_j
\coloneqq (C_1, C_2, \dots, C_j, \,  \tau_{C_j^{-1}}C_{j+1}, \, \tau_{C_j}^{-1}\tau_{C_{j+1}}^{-1}C_{j+2},
\, \dots, \, \tau_{C_j}^{-1} \cdots \tau_{C_{n-1}}^{-1}C_n)$ for $1 \leq
j \leq n$.
Here, the two extreme cases are $\boldsymbol{\Gamma}_n = (C_1, C_2, \dots,
C_n)$ and $\boldsymbol{\Gamma}_1 = (L'_1, L'_2, \dots , L'_n)$.
The collection $\boldsymbol{\Gamma}_j$ is obtained by applying $\tau_{C_j}^{-1}$
on the latter $(n - j )$ submanifolds of $\boldsymbol{\Gamma}_{j+1}$.
Now, the submanifolds in $\boldsymbol{\Gamma}_{j+1}$ around $(0, 0) \in D_{j+1}$
is disposed as in the left half of Figure \ref{fig:DjAndDj+1}.

Next, we apply $\tau_{C_j}^{-1}$ to the appropriate submanifolds and
obtain $\boldsymbol{\Gamma}_j$, the submanifolds are transformed  as the right half of Figure \ref{fig:DjAndDj+1}.
After that, we perform Dehn twists $\tau_{C_{j-1}}^{-1}, \tau_{C_{j-2}}^{-1},
\dots , \tau_{C_1}^{-1}$ to obtain $\boldsymbol{\Gamma}_1$.
Now the subset $D_j \cap D_{j+1}$ is away from the supports of the Dehn twists $\tau_{C_{j-1}}^{-1}, \tau_{C_{j-2}}^{-1},
\dots , \tau_{C_1}^{-1}$ .
Since $L_j$ coincides with $L'_j$ on $\bigcup_{l} (D_l \cap D_{l+1})$ , we can define the tails $\gamma^\pm_j$ as in Figure \ref{fig:DjAndDj+1} for $j < n$.
For the case of $\gamma^\pm_n$, we define them by $\displaystyle \gamma^\pm_n(t)
\coloneqq \left( \frac12 , \mp \varepsilon_3 t \right)$.
Finally, we have a compatible collection of two-tailed Lagrangian submanifolds
$\boldsymbol{L}^\pm = (L^\pm_0, L^\pm_1, \dots , L^\pm_n)$.
By definition, we can check that the underlying submanifolds represent a
basis of $H_1(M; \mathbb{Z})$, so the collection is perfect.
Thus, we can obtain an exact symplectic manifold $M$ and a collection of two-tailed Lagrangian
branes $\boldsymbol{L}^\maltese = (L^\maltese_0, L^\maltese_1, \dots , L^\maltese_n)$.
Here we choose that the switching point $p_j \in L_j$ to be the
root $\gamma^\pm_j(0)$.

\begin{prop}\label{prop:BasicConstrForAnothingnothing}

$\mathcal{F}^\to(\boldsymbol{L}^\maltese)$ is isomorphic to $\mathcal{A}_{\varnothing,
\varnothing}$.

\end{prop}

Before we start the proof, we prepare a notation.
Our two-tailed Lagrangian submanifold $L^\pm_j$ has the feature that $\gamma^+_j(0)$
and $\gamma^-_j(0)$ coincide.
We call such a two-tailed Lagrangian subanifold just a \textit{tailed Lagrangian
submanifold} and call $\gamma^\pm_j(0)$ a \textit{root} of $L^\pm_j$.

\begin{lem}\label{lem:TailAndNoDisc}

Let $\boldsymbol{L}^\pm = (L^\pm_0, L^\pm_1, \dots , L^\pm_n)$ be a perfect
collection of tailed Lagrangian submanifolds in $M$.
Then, an immersion $u \colon \Delta^{d+1} \to M$ passing through at least one root
of $L^\pm_j$ does not appear in the moduli spaces.

\end{lem}

\begin{proof}

Since any class in the moduli space $\mathcal{M}(y_1, y_2, \dots\ , y_d;
y_0)$ can be represented by a holomorphic map $u \colon \Delta^{d+1} \to
M$ under a suitable complex structure on $\Delta^{d+1}$, we choose such a
holomorphic representative $u$.
Assume that $u$ passes through a root $\gamma^\pm_j(0)$ for some $j$.
Since $u$ is an immersion, the image of $u$ contains at least one of these
two points $\gamma^{\pm}_j(\varepsilon_4)$.
Let us assume that $\gamma^+_j(\varepsilon_4)$ is the point.
Then, $\gamma^+_j([0, 1]) \cap u(\Delta^{d+1})$ is open in $\gamma^+_j([0,
1])$ since $u(\partial \Delta^{d+1})$ is constrained in $\cup L_j$ and by
the maximum value principle of holomorphic functions.
Obviously, this set is closed and non-empty, thus we have $\gamma^+_j([0,
1]) \cap u(\Delta^{d+1}) = \gamma^+_j([0, 1])$ hence $\gamma^+_j(1)\in u(\Delta^{d+1})$.
However, $\gamma^+_j(1) \in \partial M$ so this contradicts with the maximum
value principle.
\hfill $\Box$

\end{proof}

For a perfect collection of tailed Lagrangian submanifolds $\boldsymbol{L}^\pm$,
we write the connected component of $L_j \setminus (\bigcup_{i \neq j} L_i)$
contains the root $\gamma^\pm_j(0)$ as $iL_j$ and call it the \textit{irrelevant
part}.
Define $cL_j \coloneqq L_j \setminus iL_j$ and call it the \textit{core}
of $L_j$.
By the above lemma, we can calculate the directed subcategory $\mathcal{F}^\to(\boldsymbol{L}^\#)$
by the information of the core $c\boldsymbol{L} \coloneqq
(cL_0, cL_1, \dots cL_n)$.

\begin{figure}[hbt]
\centering
\includegraphics[width=8cm]{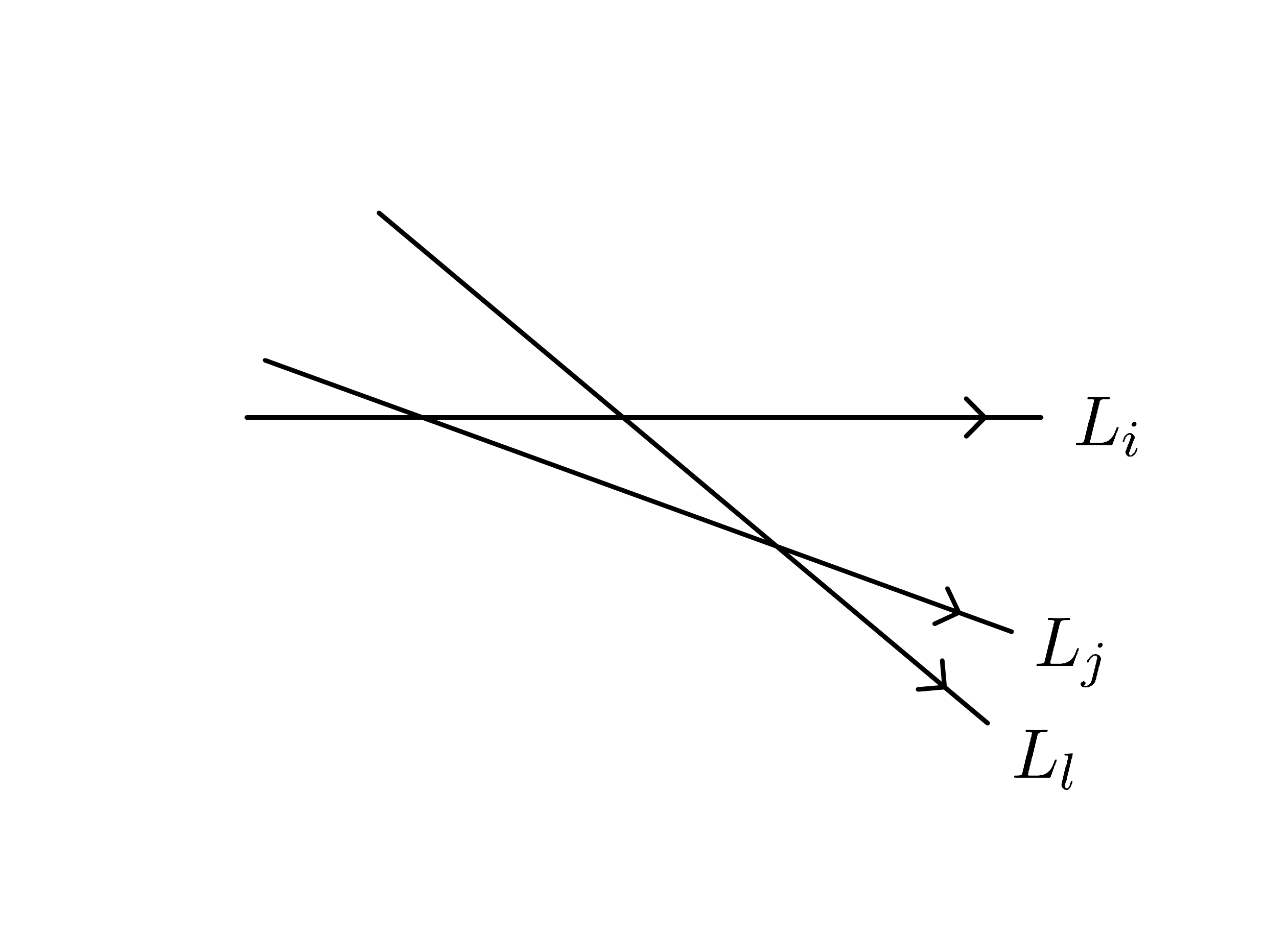}
\caption{Three Lagrangians}\label{fig:Lijk}
\end{figure}

Now, we are going to prove Proposition \ref{prop:BasicConstrForAnothingnothing}.
By the construction, the core $c\boldsymbol{L}$ is as in Figure \ref{fig:L''s}.
First, any two of submanifolds intersect at one point.
We write $p_{ij} \in L_i \cap L_j$ for $i < j$.
The differential $\mu^1$ is automatically zero because of the degree. We
can choose the gradings of $L^\#_j$'s so that the degree $|p_{j \, j+1}|$
of the morphism is zero.
Second, any three are in the position as in Figure \ref{fig:Lijk} for $i
< j < l$.
We have $\mu^2(p_{jl}, p_{ij})
= p_{il}$. Moreover, we can conclude that $|p_{ij}| = 0$ for any $i < j$.
Hence, $\mu^d$ for $d \geq 3$ are zero by the degree constraint.
Thus we have an isomorphism between $\mathcal{F}^\to(\boldsymbol{L}^\maltese)$
and $\mathcal{A}_{\varnothing , \varnothing}$.

\begin{rem}\label{rem:KoszulDualForNullNull}

By the above construction, we can show that $\tau_{L_{j-1}}L_j$ is homotopic
to $C_j$ so is $S_j = \tau_{L_1}\tau_{L_2} \cdots \tau_{L_{j-1}}L_j$.
By easy observation, we can conclude that $\hom_\mathcal{F}(S^\#_i, S^\#_j)$
is one-dimensional if $|i - j| = 1$ and zero if $|i - j| > 1$ when we choose
suitable representatives of the result of Dehn twists.
As we see later, the degree of the morphism $q^1_{j, j-1} \in \hom_\mathcal{F}(S^\#_j,
S^\#_{j-1})$ is one.
Therefore, we can conclude that the following directed $A_\infty$-category $\mathcal{B}$ is an $A_\infty$-Koszul
dual of $\mathcal{A}_{\varnothing, \varnothing}$.
Here, the directed $A_\infty$-category $\mathcal{B}$ is defined by $Ob(\mathcal{B})
= \{ B(n) > B(n-1) > \cdots > B(0) \}$, hom space is all zero but $\hom_\mathcal{B}(B(j),
B(j)) = k \cdot \eta^0_j$, $\hom_\mathcal{B}^1(B(j), B(j-1)) = k \cdot \eta^1_j$,
and $\mu$'s are all zero but $\mu^2$ with identity morphisms.
This coincides with the classical computation.

\end{rem}

\subsection{Construction (2)}\label{subsec:bypassing}

In this subsection, we construct a perfect collection of tailed Lagrangian submanifolds $\boldsymbol{L}^\pm_{S, T}$ in a two-dimensional manifold $M_{S, T}$ for $S = \{ s_1 < s_2 < \dots
< s_m \}$ and $T = \{ t_1 < t_2 < \dots < t_m \}$ from $\boldsymbol{L}^\pm$ and
$M$ in the previous section.

\begin{figure}[hbt]
\centering
\includegraphics[width=8cm]{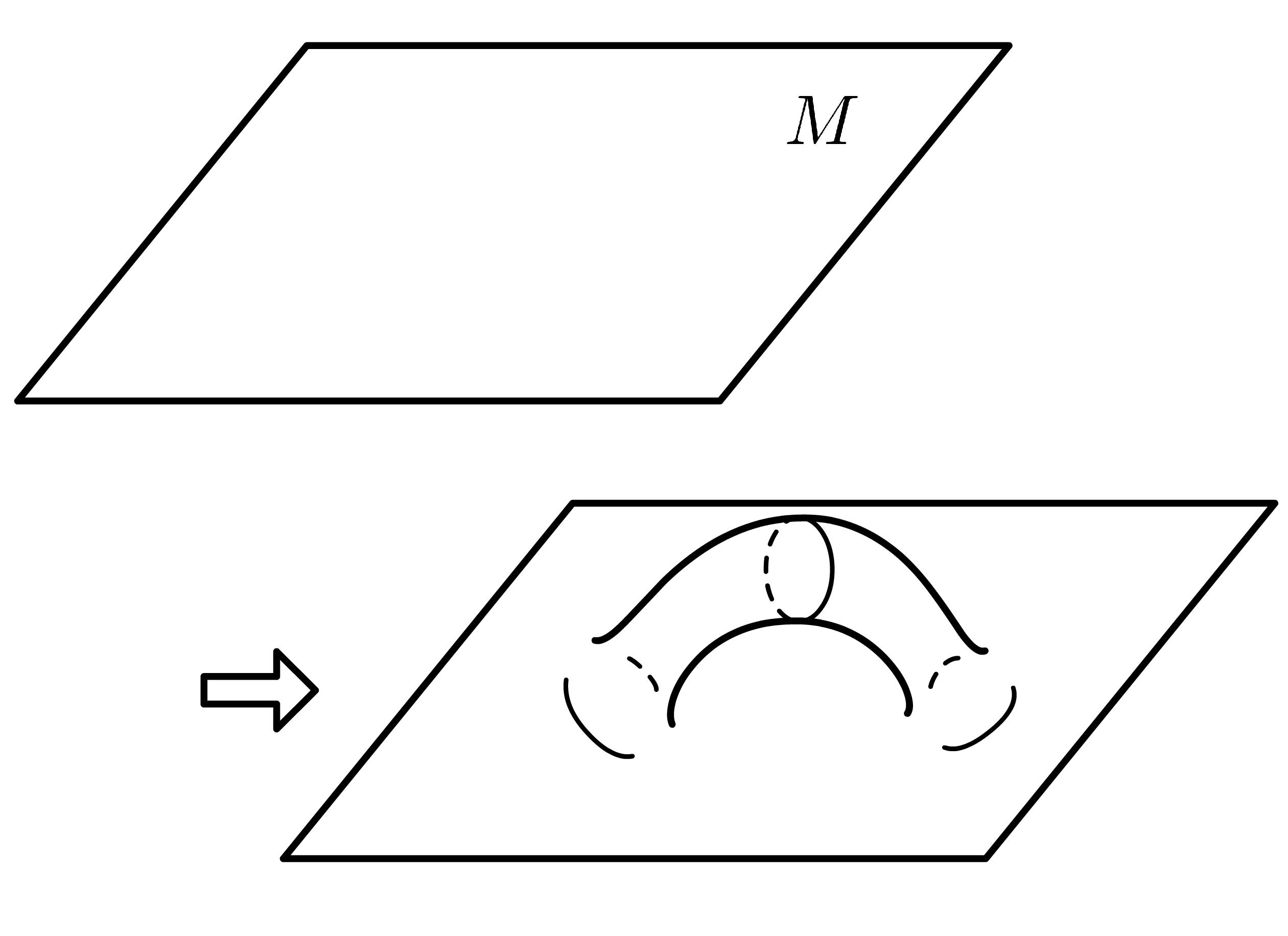}
\caption{Bypassing}\label{fig:bypassing}
\end{figure}

First, we call a surgery adding a genus as in Figure \ref{fig:bypassing}
a \textit{bypassing}.
We call the attached part a \textit{bypass}.
We identify the points in $M$ irrelevant to a bypassing with the corresponding
point in the result of the bypassing.

\begin{figure}[hbt]
\centering
\includegraphics[width=8cm]{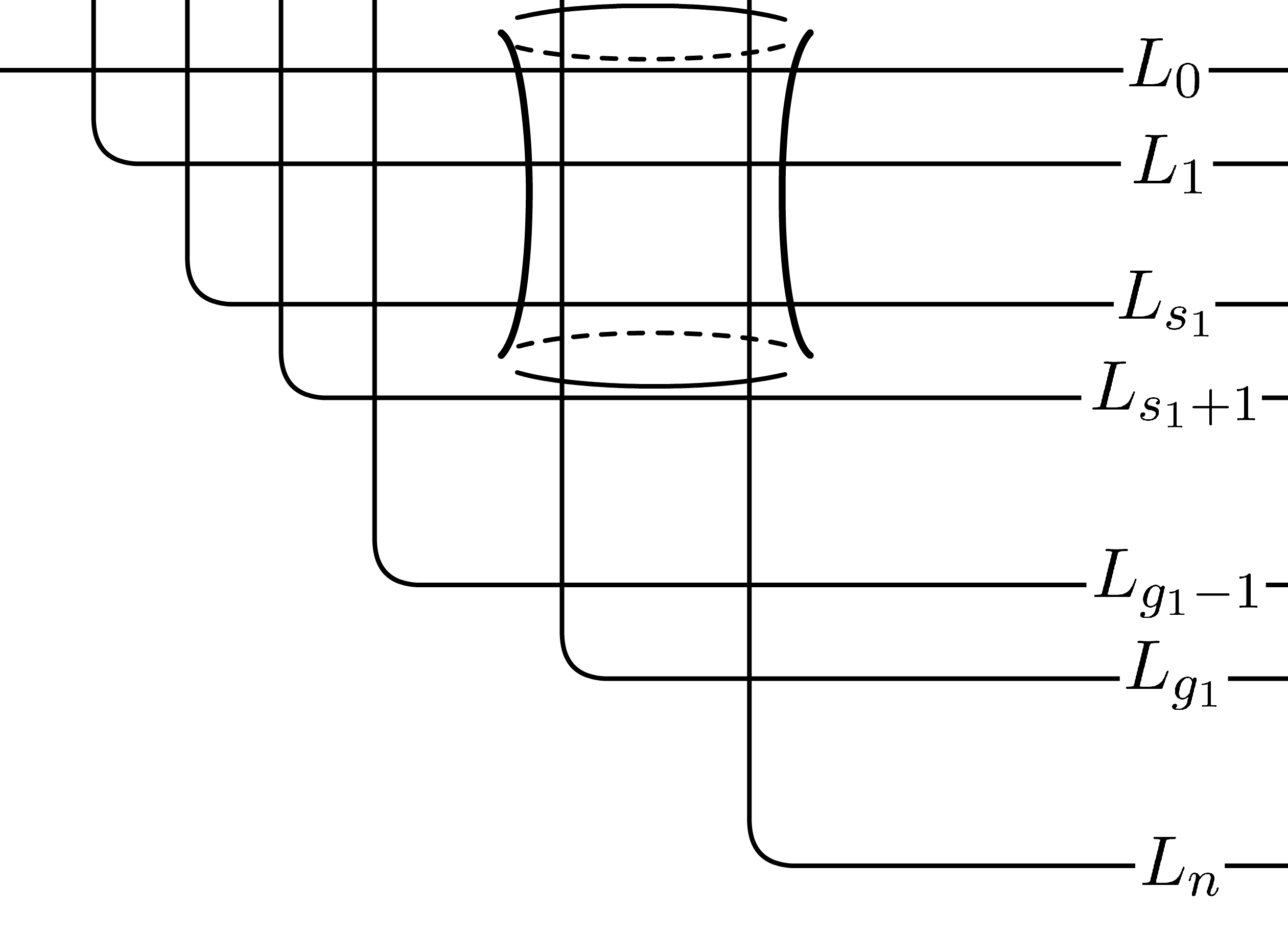}
\caption{Bypassing 2}\label{fig:bypassing2}
\end{figure}

\begin{figure}[hbt]
\centering
\includegraphics[width=8cm]{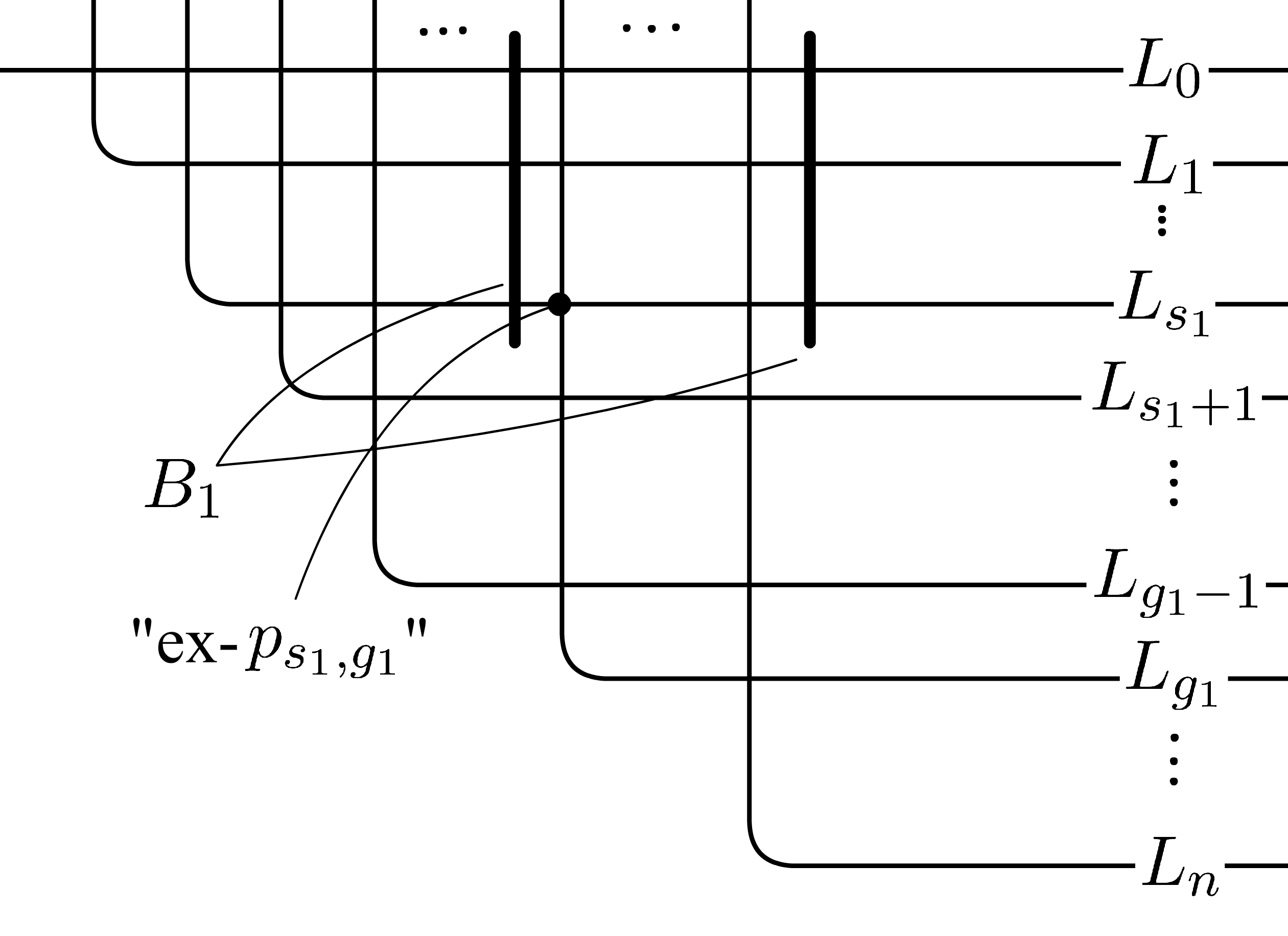}
\caption{Bypass diagram}\label{fig:bypassing3}
\end{figure}

The construction is as follows.
We write $M^{(0)} = M$ and $L^{(0)}_j \coloneqq L_j$.
First, we construct a bypass $B_1$ to remove the intersecton points $p_{ij}$
for $i \leq s_1$ and $j \geq t_1$.
The bypass $B_1$ is located around $c\boldsymbol{L}$, and the bypass across
the submanifolds $L^{(0)}_0, L^{(0)}_1, \dots , L^{(0)}_{s_1}$ as in Figure \ref{fig:bypassing2}.
We define $L^{(1)}_j$ for $0 \leq j \leq s_1$ by setting $L^{(1)}_j$ is almost
the same as $L^{(0)}_j$ but passes through under the bypass $B_1$.
For $s_1 < j < t_1$, we define $L^{(1)}_j$ to be the same as $L^{(0)}_j$.
For $j \geq t_1$, $L^{(1)}_j$ is defined to be a submanifold which is almost the same
as $L^{(0)}_j$ but go across the bypass $B_1$.
We simplify the diagram as in Figure \ref{fig:bypassing3}.

We write the resulting ambient manifold $M^{(1)}$.
By using the same tails, we have a perfect collection of tailed Lagrangian
submanifolds $\boldsymbol{L}^{(1)\pm} = (L^{(1)\pm}_0, L^{(1)\pm}_1, \dots
, L^{(1)\pm}_n)$.

\begin{prop}\label{prop:FirstBypassingAndCategory}

$(\mathcal{F}^{(1)})^\to(\boldsymbol{L}^{(1)\maltese})$ is isomorphic to
$\mathcal{A}_{\{s_1\} , \{ t_1 \} }$.

\end{prop}

\begin{proof}

It is enough to consider around the core $c\boldsymbol{L}^{(1)}$ in Figure
\ref{fig:bypassing3}.
Now, $L_i$ for $i \leq s_1$ and $L_j$ for $j \geq t_1$ no longer intersect
so we can acheive an isomorphism between $\hom_{\mathcal{F}^{(1)}}(L^{(1)\maltese}_i,
L^{(1)\maltese}_j)$ and $\hom_{\mathcal{A}_{\{s_1\} , \{ t_1 \} }}(i, j)$
by shift of the grading of Lagrangian branes if necessary.
The $A_\infty$-structures can be computed as in the same technique in the
proof of Proposition \ref{prop:BasicConstrForAnothingnothing}. Finally,
we have the desired isomorphism.
\hfill $\Box$

\end{proof}

\begin{figure}[hbt]
\centering
\includegraphics[width=8cm]{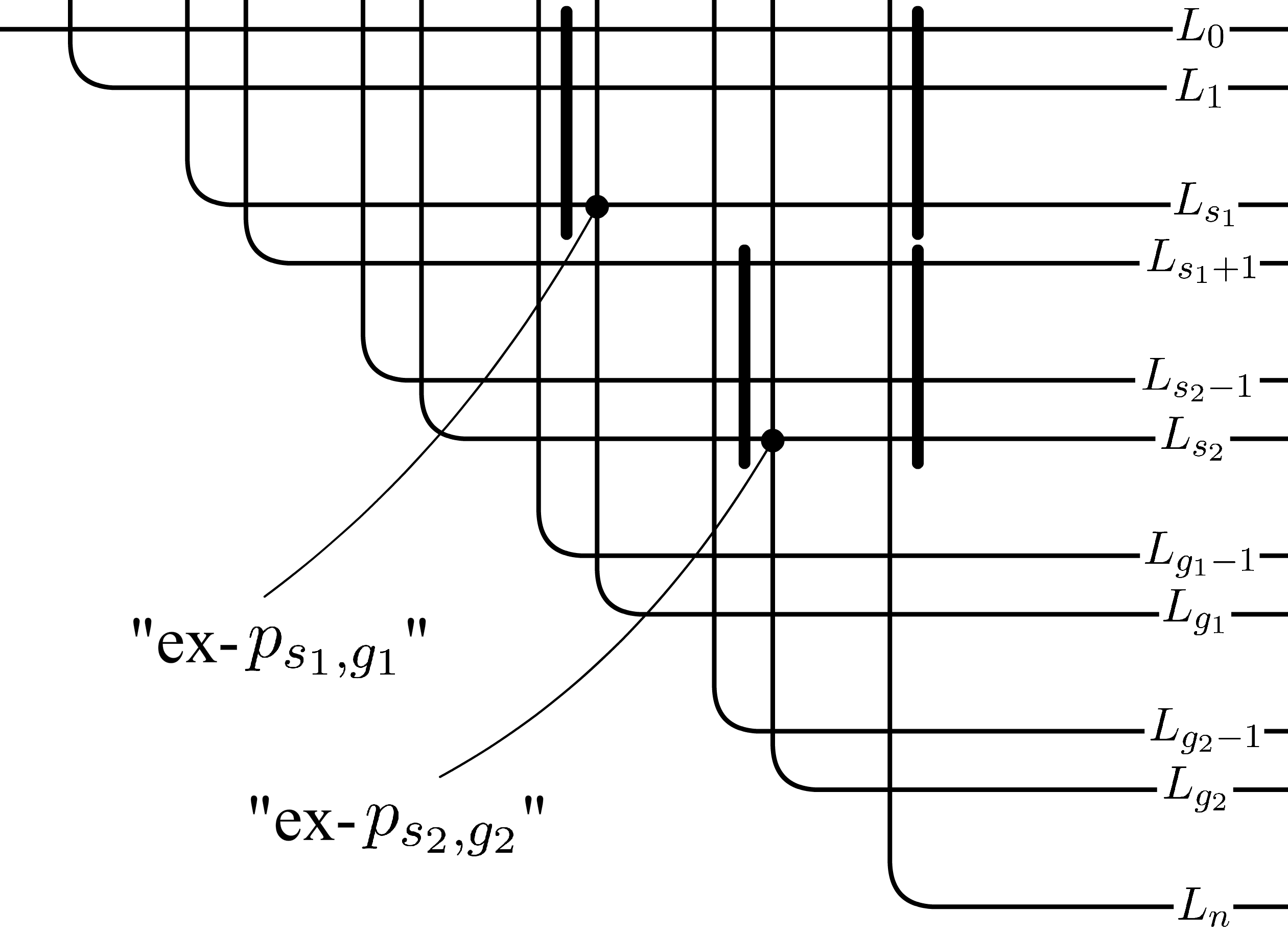}
\caption{Bypass diagram 2}\label{fig:bypassing4}
\end{figure}

Next, we construct the second bypass $B_2$ and related materials as follows.
We construct the bypass $B_2$ to be across the submanifolds $L^{(1)}_{s_1
+ 1}, L^{(1)}_{s_1 + 2}, \dots , L^{(1)}_{s_2}$.
We define new submanifolds $L^{(2)}_j$ as follows: \\
\begin{enumerate}[label=(\roman*)]
\item for $0 \leq j < s_1$, $L^{(2)}_j$ is the same as $L^{(1)}_j$;
\item for $s_1 \leq j \leq s_2$, $L^{(2)}_j$ is almost the same but passes across under the bypass $B_2$;
\item for $s_2 < j < t_2$, $L^{(2)}_j$ is the same as $L^{(1)}_j$; \\
\item for $t_2 \leq j \leq n$, $L^{(2)}_j$ is almost the same but passes across the bypass $B_2$
\end{enumerate}
as in Figure \ref{fig:bypassing4} (the figure illustrates the case $s_2 <
t_1$).
We name the resulting manifold $M^{(2)}$.

We iterate this process and obtain $M_{S, T} = M^{(m)}$ and $\boldsymbol{L}^\pm_{S,
T} = \boldsymbol{L}^{(m)\pm}$.
Moreover, by the same discussion in Proposition \ref{prop:FirstBypassingAndCategory},
we obtain the following proposition:

\begin{prop}

$(\mathcal{F}_{S, T})^\to(\boldsymbol{L}^\maltese_{S, T})$ is isomorphic
to $\mathcal{A}_{S, T}$.

\end{prop}

\clearpage

\section{Directed Fukaya categories for Riemann diagrams}\label{sec:DirFukForRiemannDiag}

In this section, we set up another construction of exact Riemann surfaces
and Lagrangian branes.
First of all, we consider a tuple $(D; l_1, l_2, \dots , l_n)$ where $D$
is an compact oriented surface with non-empty boundary and $l_j$ is an embedding
$l_j \colon [0, 1] \to D$ such that $l_j(0), l_j(1) \in \partial D$, $l_j \pitchfork \partial D$ at $l_j(0)$ and $l_j(1)$, and 2$N$
points $l_1(0), l_2(0), \dots , l_n(1) \in \partial D$ are pairwise distinct.
We call such a tuple $(D; l_1, l_2, \dots , l_n)$ a \textit{Riemann diagram}.

For a Riemann diagram $D = (D; l_1, l_2, \dots, l_n)$ we define a new compact
oriented surface $\widetilde{M}$ by attaching $N$ one-handles $H_j =[0,
1] \times [-\varepsilon_5, \varepsilon_5]$ and smoothing of the boundary.
Here, $j$-th handle is attached so that $(t, 0) \in H_j$ and $l_j(t) \in
\partial M$ are identified for $t = 0, 1$ and two distinct strips don't interfere
each other.
Next, we define a perfect collection of tailed Lagrangian submanifolds $\boldsymbol{L}^\pm
= (L^\pm_1, L^\pm_2, \dots , L^\pm_n)$ by smoothing of $l_j([0, 1]) \cup
[0, 1] \times \{ 0 \} (\subset H_j)$ and $\gamma^\pm_j(t) = (\frac12, \pm
\varepsilon_5 t) \in H_j \hookrightarrow \widetilde{M}$.
(The homological condition in Lemma \ref{lem:ConstructGrading} automatically
holds by the definition.)

Hence, we have an exact symplectic manifold and a collection of Lagrangian
branes. We write them $M_D$ and $\boldsymbol{L}^\maltese_D$.
Finally, we set $\mathcal{F}^\to_D \coloneqq (Fuk(M_D))^\to(\boldsymbol{L}^\maltese_D)$
and call it a \textit{directed Fukaya category associated with a Riemann diagram}
$D$.

\begin{figure}[hbt]
\centering
\includegraphics[width=8cm]{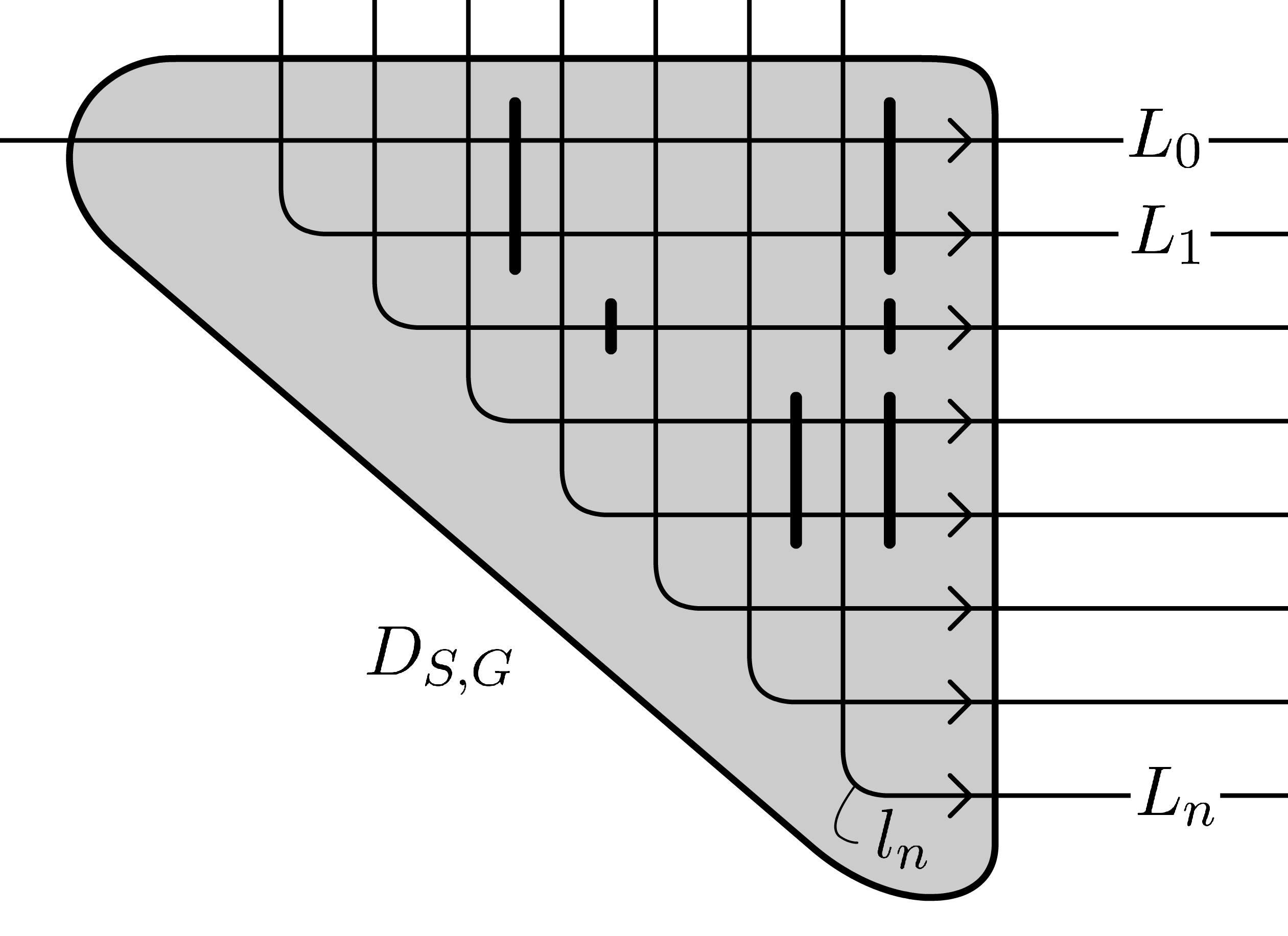}
\caption{$D_{S, T}$}\label{fig:DSG}
\end{figure}

\begin{rem}\label{rem:EquivalenceOfTwoConstruction}

Our previous construction can be reproduced when we choose a suitable
closed neighbourhood $D_{S, T}$ of the core $c\boldsymbol{L}^\pm$ of our
perfect collection of tailed Lagrangian submanifolds (and choose parametrizations
of $l_j = \colon [0, 1] \to L_j \cap D_{S, T}$ for $0 \leq j \leq n$) as
in Figure \ref{fig:DSG}.

\end{rem}

\section{Computation of Dehn twists}\label{sec:CompOfDehnTwists}

In this section, we compute an $A_\infty$-Koszul dual by using Theorem
\ref{thm:AInfKoszulDualByFukayaCat}.
We fix $n$, $S$ and $T$ and omit the subscripts, i.e.
we set $\mathcal{A} = \mathcal{A}_{S, T}$, $M = M_{S, T}$, $\mathcal{F} = \mathcal{F}_{S,
T} = Fuk(M_{S, T})$, and $\boldsymbol{L}^\maltese = \boldsymbol{L}^\maltese_{S,
T}$.
What we have to do is the following: (i) compute $\boldsymbol{S}^\# = (S^\#_n,
S^\#_{n-1}. \dots , S^\#_0)$ where $S^\#_j \coloneqq
\tau_{L^\#_0}\tau_{L^\#_1} \cdots \tau_{L^\#_{j-1}}L^\#_j$; (ii) compute
$\hom_\mathcal{F}(S^\#_i, S^\#_j)$, i.e. study $S_i \cap S_j$ and 
the degree; (iii) determine $\mu$'s i.e. count polygons.

First, let us study the intersections of $S^\#_j$'s.
However, we did not completely specify the underlying spaces of $S^\#_j$'s yet because
the twists and Dehn twists are defined only up to quasi-isomorphisms and
up to Hamiltonian isotopies respectively.
But, by Lemma \ref{lem:EquivOfDirSubcat}, we can change the representatives
of quasi-isomorphism classes of such $S^\#_j$'s in the Fukaya category.
Therefore, we can fix the convenient representatives in the following discussions.
Some of the statements in this section must start with the declaration ``with our choice of representatives" but we sometimes omit it for simplicity.

Before we begin the computation of Dehn twists, we assume one more condition: for a perfect collection of tailed Lagrangian submanifolds
$\boldsymbol{L}^\pm = (L^\pm_0, L^\pm_1, \dots , L^\pm_n)$ in $M$, there
exists a small closed neighbourhood $F^\pm_j$ of $\gamma^\pm([0, 1])$ such
that $F^\Box_i \cap F^\triangle_j \not =\varnothing$ only when $\Box = \triangle$ and $i = j$, $F^\pm_i \cap L_j = \varnothing$ for $i \neq j$, and
$\int_{F^\pm_j} \omega$ is large enough.
We can assume this by operating the surgeries in the proof of Lemma \ref{lem:ConstructExactStr}.

When this is the case, we can deform $L^\maltese_j$ freely away from $F^+_j
\cup F^-_j$ under keeping the condition that $L^\maltese_j$ is tailed Lagrangian brane, by adjustment in $F^\pm_j$.
Let us explain this.
Suppose we deform $L_j$ into $\widetilde{L_j}$ so that $L_j \cap F^\pm_j
= \widetilde{L_j} \cap F^\pm_j$.
In general, $\widetilde{L_j}$ may not satisfy $\int_{\widetilde{L_j}} \theta
= 0$.
By assumption, we can deform $\widetilde{L}$ into $\widetilde{\widetilde{L}}$
such that the deformation $\varphi$ is supported in $F^+_j \cup F^-_j$ and
$\widetilde{\widetilde{L}}$ satisfies $\int_{\widetilde{\widetilde{L}}} \theta
= 0$.
Let us consider replacing of the exact structure $(\omega , \theta , J)$ into
$(\varphi^* \omega, \varphi^* \theta, (\varphi_*)^{-1} J \varphi_*)$.
With  the new exact structure, $\widetilde{L}$ satisfy $\int_{\widetilde{L}}
\varphi^*\theta = 0$.
Moreover, by construction, $\mathcal{F}^\to(\boldsymbol{L}^\maltese)$,
$\mathcal{F}^\to(\widetilde{\widetilde{\boldsymbol{L}}}^\maltese)$,
and $\mathcal{F}'{}^\to(\widetilde{\boldsymbol{L}}^\maltese)$
are isomorphic, where $\mathcal{F}' = Fuk(M, \varphi^*\omega, \varphi^*\theta,
(\varphi_*)^{-1} J \varphi_*)$, $\widetilde{\boldsymbol{L}}^\maltese$ and
$\widetilde{\widetilde{\boldsymbol{L}}}^\maltese$ are collection of Lagrangian
branes which are obtained by replacing $L^\maltese_j$ into $\widetilde{L_j}^\maltese$
and $\widetilde{\widetilde{L_j}}^\maltese$ in $(M, \varphi^*\omega, \varphi^*\theta,
(\varphi_*)^{-1} J \varphi_*)$ and $(M, \omega, \theta, J)$ respectively.
Therefore what we have is the following: 

\begin{lem}

We use the same symbol as above.
When we deform $\boldsymbol{L}$ into $\widetilde{\boldsymbol{L}}$, there exist an exact symplectic structure $(\omega', \theta', J')$ on $M$ and a collection of tailed Lagrangian brane structure of $\widetilde{\boldsymbol{L}}$ in $M' = (M, \omega', \theta', J')$ such that $\mathcal{F}^\to(\boldsymbol{L}^\maltese) = \mathcal{F}'{}^\to(\widetilde{\boldsymbol{L}}^\maltese)$.

\end{lem}

We operate such deformations and replacements of the exact structure
without noticing in the following discussion.

\subsection{Choice of representataives}

\begin{figure}[hbt]
\centering
\includegraphics[width=8cm]{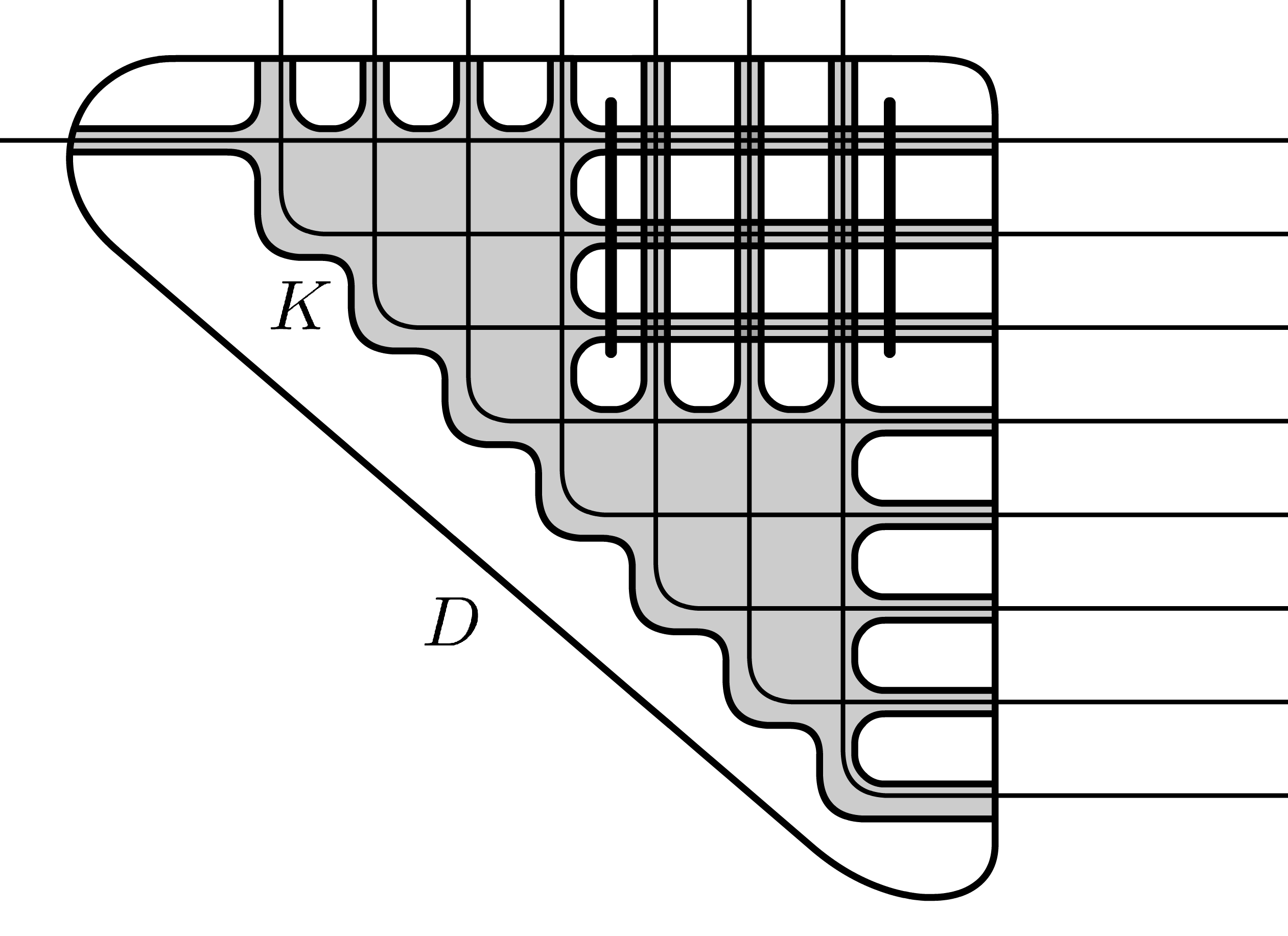}
\caption{$K \subset D$}\label{fig:K}
\end{figure}

To see the general case, we again consider $M = M_{S, T}$ as $M_{D_{S, T}}$ as
in Remark \ref{rem:EquivalenceOfTwoConstruction}.
First, we prepare some notations.
We define a closed subset $F = F_{S, T}$ of $D = D_{S, T}$ by the union of
$\bigcup_{j} L_j \cap D$ and triangles in $D$ encircled by $L_j$'s.
By definition, $F$ is contractible.
Then, we choose a small closed neighbourhood $K$ of $F$ and fix it such that
$K$ is diffeomorphic to the unit disc and is contained in $\varepsilon_6$-neighbourhood
of $F$.
Figure \ref{fig:K} illustrates the situation.
(The neighbourhood $K$ in Figure \ref{fig:K} contains many corners for the sake of  the simplification of the figure, but we consider that the actual $K$ does not have such corners.)

We fix the orientations of $l_j$'s as in Figure \ref{fig:DSG}.
The orientation of the forthcoming $S^1$'s are induced by these orientations.
We use the symbol $L^{(p)}_j$ again but another meaning as in subsection
\ref{subsec:bypassing}.
We define $L'{}^{(p)}_j$ as follows:
\begin{align*}
L'{}^{(p)}_j \coloneqq
\begin{cases}
\tau_{L_p}\tau_{L_{p+1}} \cdots \tau_{L_{j-1}}L_j & \text{if } j > p \\
L_j & \text{if } j \leq p
\end{cases}
\end{align*}
and we define $\boldsymbol{L}'{}^{(p)} \coloneqq (L'{}^{(p)}_0, L'{}^{(p)}_1,
\dots , L'{}^{(p)}_n)$.
Now we assume the following two conditions: the first is that the support
of each Dehn twist is enough thin, i.e. $\text{Supp}(\tau_{L_j}) \cap D \subset
K$.
Roughly speaking, this means that the ``width" of the support is enough smaller
than the ``length of edges of the grid in $D = D_{S, T}$"(See Figure \ref{fig:K}).
The second is that $\tau_{L_p} L_p = L_p$ and $\tau_{L_p}(L_p \cap D) \cap
D = \varnothing$.

A collection $\boldsymbol{N} = (N_0, N_1, \dots , N_n)$ of one-dimensional submanifolds in $M$ is said to have the \textit{property $P_p$} when the following
conditions are satisfied:
\begin{enumerate}[label=(\roman*)]
\item for distinct $i, j \geq p$, the intersection
$N_i \cap N_j$ does not contained in $D$;
\item for any $i \in \{ 0, 1, \dots , n \}$ and $j < p$, the intersection $N_i \cap N_j$ is contained in $K$.
\end{enumerate}
First of all, we prove the following lemma:

\begin{lem}

The collection $\boldsymbol{L}'{}^{(p)}$ has the property $P_p$.

\end{lem}

\begin{proof}

We prove this by induction.
First, for $p = n$, our collection $\boldsymbol{L}'{}^{(p)} = \boldsymbol{L}'
= (L'_0, L'_1, \dots , L'_n)$ has property $P_n$ by the definition itself.
Now, we prove the property $P_{p-1}$ for $\boldsymbol{L}'{}^{(p-1)}$ under
the assumption that $\boldsymbol{L}'{}^{(p)}$ has the property $P_p$.
Since the support $\text{Supp}(\tau_{L_{p-1}})$ is enough thin, we can assume
that $\text{Supp}(\tau_{L_{p-1}}) \cap H_j = \varnothing$ for $j \geq p$
and $\text{Supp}(\tau_{L_{p-1}}) \cap (L'{}^{(p)}_i \cap L'{}^{(p)}_j) =
\varnothing$ for distinct $i, j < p-1$.
Thus, we can deduce that the intersections $L'{}^{(p-1)}_i \cap L'{}^{(p-1)}_j$
for distinct $i, j \geq p$ remain in $M \setminus D$ and the intersections
$L^{(p-1)}_i \cap L^{(p-1)}_j$ for distinct $i, j < p-1$ remain in $K$.

Now, we study the intersection of $L'{}^{(p-1)}_{p-1} \cap L'{}^{(p-1)}_i$.
By the condition (ii) of property $P_p$, we have $L'{}^{(p)}_{p-1} \cap L'{}^{(p)}_i
\in K \subset D$ for any $i \neq p-1$. By the assumption that $\tau_{L_{p-1}}(L_{p-1}
\cap D)
\cap D = \varnothing$, we have $L'{}^{(p-1)}_{p-1} \cap L'{}^{(p-1)}_i \in
M \setminus D$.
Thus we have proved the property $P_{p-1}$ of $\boldsymbol{L}'{}^{(p-1)}$.
\hfill $\Box$
\end{proof}

Now we modify $\boldsymbol{L}'{}^{(p)}$ by isotopies.
First, we set $\boldsymbol{L}^{(n)} = \boldsymbol{L}'{}^{(n)}$.
Before we construct isotopy, we prove the following lemma:

\begin{lem}

$\boldsymbol{L}'^{(n-1)}$ can be a collection of underlying spaces of a perfect
collection of tailed Lagrangian submanifolds.

\end{lem}

\begin{figure}[hbt]
\centering
\includegraphics[width=8cm]{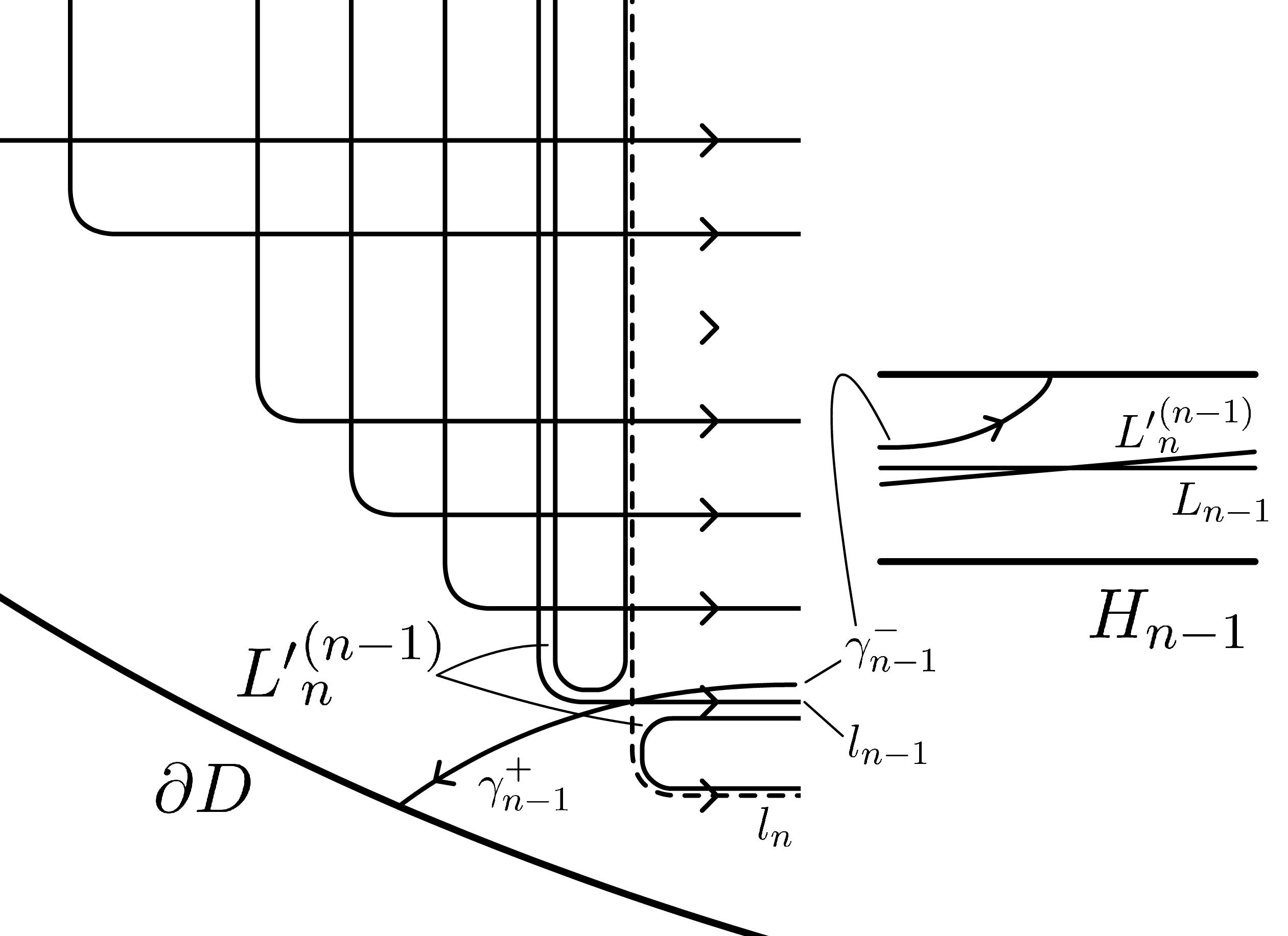}
\caption{$\boldsymbol{L}'{}^{(n-1)}$}\label{fig:L'(n-1)}
\end{figure}

\begin{proof}

For $j \neq n-1$, we can define tails $\gamma^\pm_j$ so that their image
are in $H_j$ since $L_i \cap H_j \neq \varnothing$ if and only if $i = j$.
For $j = n-1$, we can define $\gamma^\pm_{j-1}$ as in Figure \ref{fig:L'(n-1)}.
\hfill $\Box$
\end{proof}

\begin{figure}[hbt]
\centering
\includegraphics[width=8cm]{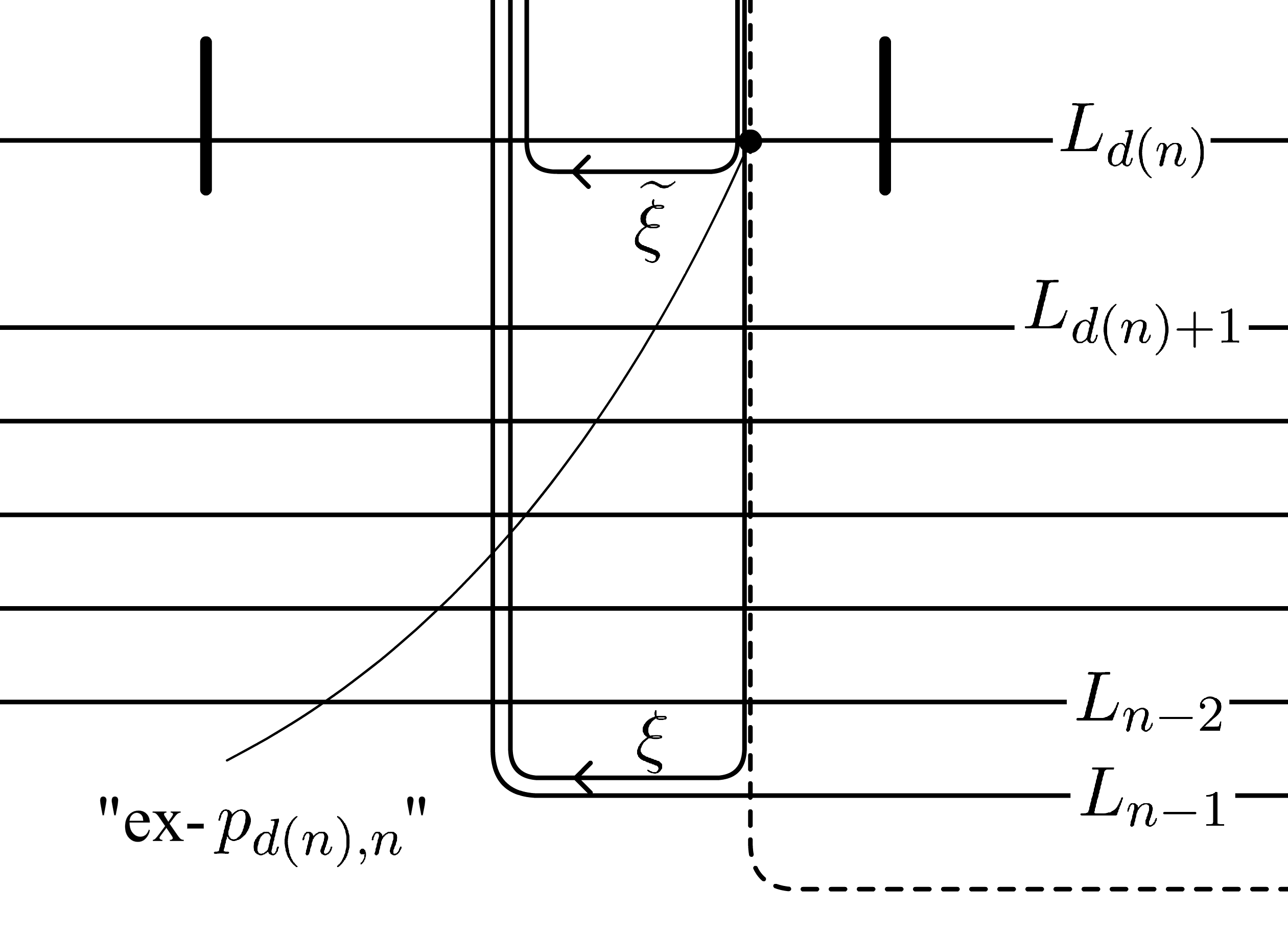}
\caption{$\xi$ and $\widetilde\xi$}\label{fig:xiandwidetildexi}
\end{figure}

Now we can deform $L'{}^{(n-1)}_n$ freely away from $\gamma^\pm_n$ as in the sense of the first discussion of this section. The first
isotopy for $\boldsymbol{L}'{}^{(n-1)}$ is constructed as follows.
By the assumptions, the connected component $\xi$ of $L'{}^{(n-1)}_n \cap
D = \tau_{L_{n-1}}L_n \cap D$ which lies in the left side of $l_{n-1}$ in
$K$ has the shape that $\xi$ comes from $l_n(0)$, go along $l_n$, turn right
just before $l_n$ reaches the $p_{n-1, n} \in L_{n-1} \cap L_n$, go along
the left side (with respect to the orientation of $l_{n-1}$) of $l_{n-1}$,
and finally go out from $D$ at the left side of $l_{n-1}(0)$ as in  Figure
\ref{fig:xiandwidetildexi}.
We call such a path $\xi$ whose endpoints are near $l_{n-1}(0)$ and $l_n(0)$
a \textit{path of type} $(n-1, n)$. (We always write the smaller one at the
left and bigger one at the right independent of the orientation of the
path.)

We make some observations about the intersections of $\boldsymbol{L}'{}^{(n-1)}$.
We can see that $\xi$ intersects $L_j$ for $d(n) < j < n-1$ twice.
We deform $\xi$ into $\widetilde{\xi}$ so that the pairs of intersections
disappear as in Figure \ref{fig:xiandwidetildexi} and deform $L'{}^{(n-1)}_n$ into
$L^{(n-1)}_n$ by $L^{(n-1)}_n = (L'{}^{(n-1)}_n \setminus \xi ) \cup \widetilde{\xi}$.
We call such a path $\widetilde{\xi}$ which come from a point near $l_n(0)$, go along
$l_{n-1}$, turn right and go ``between" $l_{d(n)}$ and $l_{d(n)+1}$, turn
right and go along $l_{n-1}(0)$, and go out at a point near $l_{n-1}(0)$, a path of
type $(n-1, n; d(n))$. 
Now, if $d(n-1) < d(n)$ then $L^{(n-1)}_n$ intersects with $L_{d(n)}$, and
if $d(n-1) = d(n)$ then $L^{(n-1)}_n$ intersects only with $L_{n-1}$.
Finally, we set $L^{(n-1)}_j = L'{}^{(n-1)}_j$ for $j \neq n$ and we have
$\boldsymbol{L}^{(n-1)}$.
The isotopy just reduces the intersections, so our collection $\boldsymbol{L}^{(n-1)}$
again has property $P_{n-1}$.

Moreover, $\boldsymbol{L}^{(n-1)}$ has the following properties:
\begin{enumerate}[label=(\roman*)]
\item for any $i \neq j$, the number of intersection points of $L^{(n-1)}_i$ and $L^{(n-1)}_j$ is at most one;
\item the collection can be a collection of underlying spaces
of a perfect collection of tailed Lagrangian submanifolds;
\item $L''{}^{(n-1)}_j \cap D \subset K$ for all $j$.
\end{enumerate}
When a collection of submanifolds has the property $P_p$ and the above property,
we say that the collection has the \textit{property} $Q_p$.

Next, we iterate the following procedure to construct $\boldsymbol{L}^{(p-1)}$
satisfying the property $Q_{p-1}$ from a collection $\boldsymbol{L}^{(p)}$
satisfying the property $Q_p$.
First, we set $L''{}^{(p-1)}_j \coloneqq L^{(p)}_j$ for $j < p$ and $L''{}^{(p-1)}_j
\coloneqq \tau_{L_{p-1}}L^{(p)}_j$ for $j \geq p$.

We write the left side of $l_{j-1}$ in $K$ by $K^\text{left}_{j-1}$.
Since $\boldsymbol{L}^{(p)}$ has the property $Q_p$, the connected component
$\xi$ of $L''{}^{(p)}_j \cap K^\text{left}_{p-1}$ for $j > p-1$ can be classified
into the following two cases:
\begin{enumerate}[label=(\alph*)]
\item one of the endpoints of $\xi$ is on the image of $l_{p-1}$ and the other is on $\partial K$ near $l_i(0)$ for some $i > p-1$, we call such a path a path of type $<i>$;
\item the two endpoints of $\xi$ are on $\partial K$ near $l_{i_1}(0)$ and $l_{i_2}(0)$ respectively for some $p-1 < i_1 < i_2$, i.e. $\xi$ is a path of type $(i_1, i_2; i_3)$ for some $i_3 < p-1$.
\end{enumerate}
We can assume that all paths of type $(i_1, i_2; i_3)$ do not intersect with the
support of $\tau_{L_{p-1}}$.

After applying of $\tau_{L_{p-1}}$ to $L^{(p)}_j$ for $j \geq p$, any connected
component $\xi$ of $L''{}^{(p-1)}_j \cap K^\text{left}_{j-1}$ is of type
$(i_1, i_2; i_3)$ with some $i_3 < p-1 \leq i_1 < i_2$.
Now, such a path $\xi$ of type $(i_1, i_2; i_3)$ intersects with $L_j$ more
than once only when $d(i_2) < j \leq i_3$.
At that time, the number of intersections is two and the intersection points
can be removed by isotopy as in the case of $p = n$. After the isotopy,
we obtain a path $\widetilde{\xi}$ of type $(i_1, i_2; d(i_2))$.
Observe that if $d(i_1) = d(i_2)$, then $\widetilde{\xi}$ does not intersects
with $L_j$ for $j < p-1$, and if $d(i_1) < d(i_2)$, then $\widetilde{\xi}$
intersects with $L_{d(i_2)}$.
We change all such $\xi$'s into $\widetilde{\xi}$'s by the isotopies.

If we apply the isotopies for suitable order, all the isotopies
just reduces intersections and not create new intersection points.
(Such an order can be constructed as follows. Any path $\xi$ divides $K$
into two regions and one is contained in $K^\text{left}_{p-1}$. We name the contained region
$K_\xi$. We define partial order of paths by $\xi < \xi' \Leftrightarrow
K_\xi \subset K_{\xi'}$. This is well-defined by the condition of property $P_{p-1}$. We add more relation and make it a total order.
This is what we want.)
Finally, we apply the isotopies for corresponding $L''{}^{(p-1)}_j$ and obtain
$L^{(p-1)}_j$.

\begin{figure}[hbt]
\centering
\includegraphics[width=8cm]{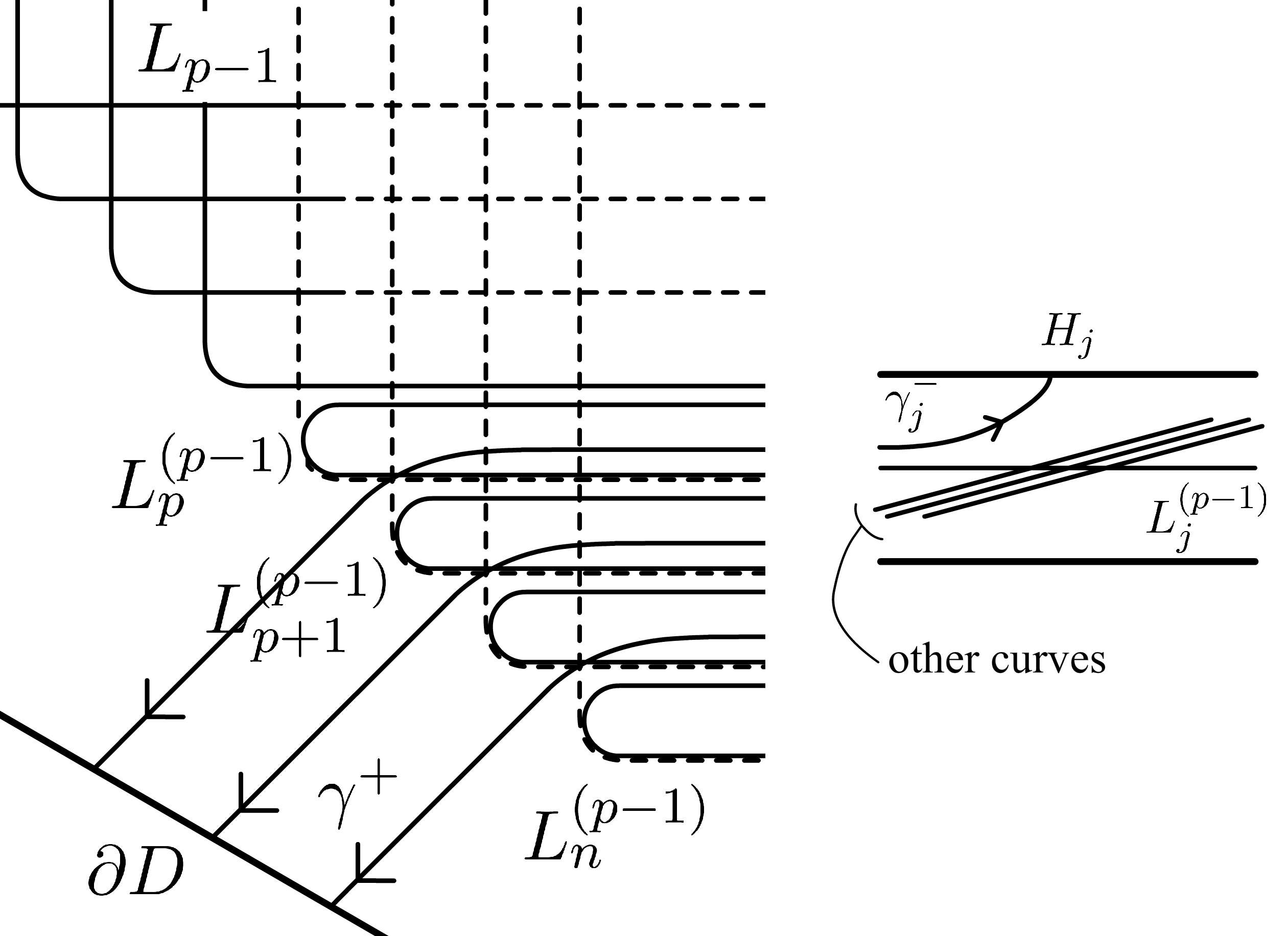}
\caption{$\boldsymbol{L}^{(p-1)}$ and their tails}\label{fig:L(k-1)}
\end{figure}

By the construction, we have shown that the collection $\boldsymbol{L}^{(p-1)}$
of such $L{}^{(p-1)}_j$'s has the property $Q_{p-1}$ except for the property
that $\boldsymbol{L}^{(p-1)}$ can be a collection of underlying spaces of
a perfect collection of tailed Lagrangian submanifolds.
Therefore, it's time to check the condition (ii) for $\boldsymbol{L}''^{(p-1)}$
and $\boldsymbol{L}^{(p-1)}$.
Since each isotopy takes place in $K^\text{left}_{j}$ for a suitable $p-1
\leq j < n$,
we have $L''{}^{(p-1)}_j \cap K^\text{right}_{p-1} = L^{(p-1)}_j \cap K^\text{right}_{p-1}
= \tau_{L_{j-1}}L_j \cap K^\text{right}_{p-1}$ for $j > p-1$.
Here, $K^\text{right}_{p-1}$ is the right part of $l_{p-1}$ in $K$.
Hence, $\boldsymbol{L}''^{(p-1)}$ and $\boldsymbol{L}^{(p-1)}$ can be drown
as in Figure \ref{fig:L(k-1)}.
Therefore, we can define $\gamma^\pm_j$ as in Figure \ref{fig:L(k-1)}.
Thus, we have constructed a collection $\boldsymbol{L}^{(p-1)}$ with property
$Q_{p-1}$.

Finally, we have a perfect collection of tailed Lagrangian submanifolds $\boldsymbol{L}^{(0)}{}^\pm
= (L^{(0)}_0{}^\pm, L^{(0)}_1{}^\pm, \dots , L^{(0)}_n{}^\pm)$.
The gradings are induced by those of $L^\maltese$'s.
In fact, each $L^{(0)}_j$ shares an interval with $L_j$ in $H_j$, so they have grading
$\alpha^{(0)}_j$ such that $\alpha^{(0)}_j = \alpha_j$ on $L^{(0)}_j \cap
L_j \cap H_j$.
We specify that the switching point $q_j$ of $L^{(0)}_j$ as the root of the
tail.

Thus, we have a perfect collection of tailed Lagrangian branes $\boldsymbol{L}^{(0)\maltese}$
and define another collection $\boldsymbol{S}^\maltese = (S^\maltese_n, S^\maltese_{n-1},
\dots , S^\maltese_0)$, where $S^\maltese_j = L^{(0)}_j{}^\maltese$.
Finally, by Theorem \ref{thm:AInfKoszulDualByFukayaCat},  $\mathcal{F}^\to(\boldsymbol{S}^\maltese)$
is an $A_\infty$-Koszul dual of $\mathcal{A} = \mathcal{A}_{S, T} = \mathcal{A}(R_{S,
T})$.

\begin{rem}

All the isotopies used to construct $\boldsymbol{L}^{(p)}$ are taken place
in $K$ so they don't affect the intersection of $L^{(p)}_j$'s for $j >
p$.
In fact, there exists a diffeomorphism $f \colon M \to M$ such that $f(L''{}^{(p)}_j)
= f(L^{(p)}_j)$ for $j > p$.
Thus, we have an isomorphism of $A_\infty$-categories between $\mathcal{F}^\to((L''{}^{(p)}_p{}^\maltese,
L''{}^{(p)}_{p+1}{}^\maltese, \dots , L''{}^{(p)}_n{}^\maltese))$ and
$\mathcal{F}^\to((L^{(p)}_p{}^\maltese,L^{(p)}_{p+1}{}^\maltese, \dots
, L^{(p)}_n{}^\maltese))$.

Moreover, the isotopies act on  $L''{}^{(p)}_j$ for $j > p$, we have
$L^{(p)}_j = L_j$ for $j \leq p$.

\end{rem}

\subsection{Intersections}

In this subsection, we prove the following propositions:

\begin{prop}\label{prop:WhenIntersect}

For $i < p$, $S_p \cap S_i \neq \varnothing$ if and only if there exists
$1 \leq j \leq l_p$ such that $a^{(p)}_j = i$.

\end{prop}

Together with the inversion formula (Lemma \ref{lem:InversionFormula}), we
have the following corollary:

\begin{cor}

For $i > p$, $S_i \cap S_p \neq \varnothing$ if and only if there exists
$1 \leq j \leq l^\dagger_p$ such that $a^{(p)\dagger}_j = i$.

\end{cor}

We name the unique intersection point of $S_j$ and $S_i$ for $i < j$ by $q_{j,
i}$

\begin{prop}\label{prop:OrderOfIntersection}

Along the orientation of $S_p$, the following points in $S_p$ appear in the
following order:
\begin{align*}
q_p,
q_{a^{(p)\dagger}_1, p}, q_{a^{(p)\dagger}_3, p},
\dots , q_{a^{(p)\dagger}_{l^\dagger_p}, p},\dots ,
q_{a^{(p)\dagger}_4,p}, q_{a^{(p)\dagger}_2, p},
q_{p, a^{(p)}_2}, q_{p, a^{(p)}_4},
\dots , q_{p,a^{(p)}_{l_p}}, es\dots ,
q_{p, a^{(p)}_3}, q_{p, a^{(p)}_1}.
\end{align*}

\end{prop}

We prove these propositions in the following discussion.

\subsubsection{Proof of Proposition \ref{prop:WhenIntersect}}

\begin{figure}[hbt]
\centering
\includegraphics[width=8cm]{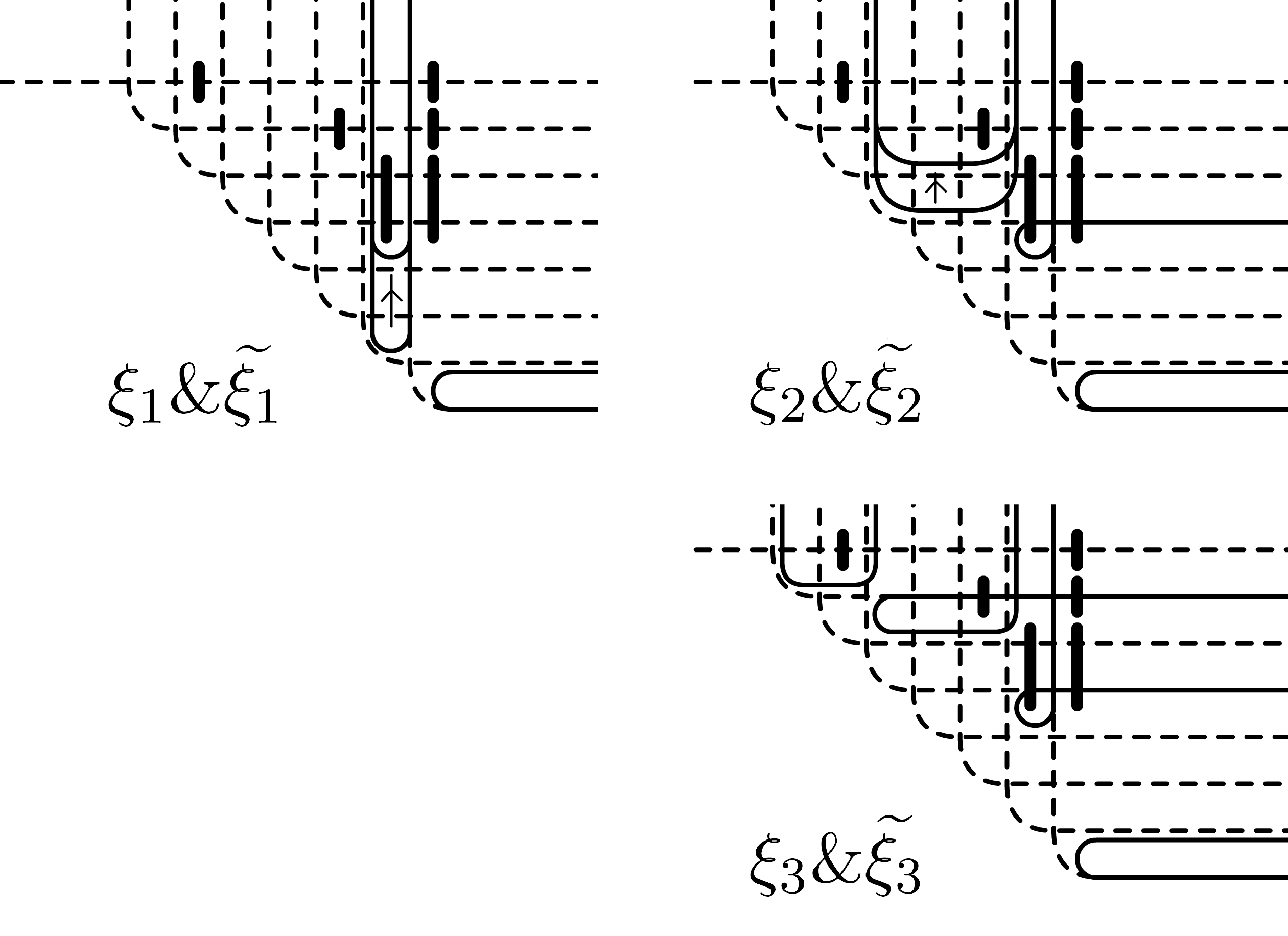}
\caption{$\xi$'s and $\widetilde\xi$'s}\label{fig:xi'sandwidetildexi's}
\end{figure}

Figure \ref{fig:xi'sandwidetildexi's} illustrates the key point of our proof so please refer
the figure when it is  needed.

When we consider the construction of $S_j$'s of $j \leq p$, we can ignore the submanifolds $L^{(*)}_{j'}$ for $j' > p$ since the
Dehn twists $\tau_{L_{j'}}$ and the isotopies used to construct $L^{(*)}_{j'}$
for $* > p$ is irrelevant for the construction.

Recall our construction of $S_j$'s. A path $\xi$ of type $(p, q)$ intersects with $L_j$ of $d(q) < j < p$ twice and intersects with $L_j$ with $d(p) < j \leq d(q)$ once.
If the path $\xi$ is a part of relevant submanifold appeared in the construction, we deform it into a path $\widetilde\xi$ of type $(p, q; d(q))$ to eliminate the pairs of intersection points with $L_j$ of $d(q) < j < p$.
As a result, $\widetilde\xi$ intersects with $L_j$ of  $d(p) < j \leq d(q)$ once and no longer intersects with $L_j$ of $d(q) < j < p$.

Let us consider the submanifold $L''{}^{(p-1)}_p = \tau_{L_{p-1}}L_p$.
By the construction, $L''{}^{(p-1)}_p \cap K^\text{left}_{p-1}$ has only one
connected component $\xi_1$ and this is a path of type $(p-1, p)$.
We deform this path into a path $\widetilde{\xi_1}$ of type $(p-1, p; d(p))$.
If $d(p) = d(p-1)$ i.e. $d(a^{(p)}_0)
= d(a^{(p)}_1)$ ($\Leftrightarrow l_p = 1$), then $\widetilde{\xi_1}$ no longer intersects with $L_j$
with $j < p-1$ neither does $L^{(p-1)}_p$.
Hence, the remaining Dehn twists and isotopies do not interact with $L^{(p-1)}_p$.
Therefore, we can deduce that $S_p = L^{(p-1)}_p$ and $S_p$ does not intersect with
$S_j$ with $j < p-1$.
Thus we have proved  Proposition \ref{prop:WhenIntersect} for the case $l_j
= 1$.

Now, suppose that $d(p) \neq d(p-1)$.
Recall that the path  $\widetilde{\xi_1}$ is a path of type $(p-1, p; d(p))$.
The path $\widetilde{\xi_1}$ intersects
with $L_{d(p)}$ and does not intersect with $L_j$ for $d(p) < j < p-1$.
Therefore, we have $L^{(d(p) + 1)}_p = L^{(d(p) + 2)}_p = \cdots = L^{(p-1)}_p$.
Let us consider $L''{}^{(d(p))}_p = \tau_{L_{d(p)}}L^{(p-1)}_p$.
By construction, $L''{}^{(d(p))}_p \cap K^\text{left}_{d(p)}$ has only one
connected component $\xi_2$, a path of type $(d(p), p-1)$.
The path $\xi_2$ is deformed into a path $\widetilde{\xi_2}$ of type $(d(p),
p-1; d(p-1))$ by our isotopy.
If $d(p-1) = d(d(p))$, equivalently $d(a^{(p)}_1) = d(a^{(p)}_2)$ ($\Leftrightarrow l_p = 2$), then $\widetilde{\xi_2}$
no longer intersects with $L_j$ of $j < d(p)$.
Hence, by the same argument, we can finish the proof of Proposition \ref{prop:WhenIntersect}
for the case $l_p = 2$.
If $d(a^{(p)}_1) \neq d(a^{(p)}_2)$, the path $\widetilde{\xi_2}$ intersects
with $L_{d(p-1)}$ and does not intersect with $L_j$ with $d(p-1) < j < d(p)$.
By the same argument, we can deduce the following: the submanifolds coincide
$L^{(d(p-1) + 1)}_p = \cdots = L^{(d(p))}_p$; the subset $L''{}^{(d(p-1))}_p \cap K^\text{left}_{d(p-1)} = \tau_{L_{d(p-1)}}L^{(d(p))}_p
\cap K^\text{left}_{d(p-1)}$ has only one connected component $\xi_3$ which
is a path of type $(d(p-1), d(p))$; the path $\xi_3$ is deformed into a path
$\widetilde{\xi_3}$ of type $(d(p-1), d(p); d(d(p))) = (a^{(p)}_3, a^{(p)}_2;
d(a^{(p)}_2))$.
If $d(a^{(p)}_2) = d(a^{(p)}_3)$, this the end of the proof of Proposition
\ref{prop:WhenIntersect} for the case $l_p = 3$.
If $d(a^{(p)}_2) \neq d(a^{(p)}_3)$, then we should iterate the above procedure.

The procedure is as follows.
Let us assume that $L^{(a^{(p)}_{j+1}+1)}_p \cap K^\text{left}_{a^{(p)}_{j+1}}$ is a path of type $<a^{(p)}_j>$ and $\bigg\{ l \, \bigg| \, L^{(a^{(p)}_{j+1}+1)}_p \cap L^{(a^{(p)}_{j+1}+1)}_l \neq \varnothing , a^{(p)}_{j+1} < l < p \bigg\}$ = $\left\{ a^{(p)}_i \, \big| \, 1 \leq i \leq j \right\}$.
(Whenever we write the symbol $a^{(p)}_{j+1}$, we assume that $l_p \geq j+1$.)
Here, recall that a path of type $<s>$ in $K^\text{left}_t$ for $s > t$ is a path with two end points, one of the end points is located in $\partial K$ near $l_s(0)$ and the other point is located in the image of $l_t$.
Note that this hypothesis with $j = 0$ is always true.
Define $L''{}^{(a^{(p)}_{j+1})}_p \coloneqq \tau_{L_{a^{(p)}_{j+1}}} L^{(a^{(p)}_{j+1}+1)}_p$.
By the first hypothesis, $L''{}^{(a^{(p)}_{j+1})}_p \cap K^\text{left}_{a^{(p)}_{j+1}}$ is a path $\xi_{j+1}$ of $(a^{(p)}_{j+1}, a^{(p)}_j)$.
Thus, by definition of $L^{(a^{(p)}_{j+1})}_p$, $\xi_{j+1}$ is deformed into a path $\widetilde{\xi_{j+1}}$ of type $(a^{(p)}_{j+1}, a^{(p)}_j; d(a^{(p)}_j))$ and we have $I_{j+1} = \left\{ l \, \big| \, \widetilde{\xi_{j+1}} \cap L_l \neq \varnothing , l < a^{(p)}_{j+1} \right\} = \left\{ l \, \big| \, d(a^{(p)}_{j+1}) < l \leq d(a^{(p)}_j) \right\}$.
Moreover, we have $\bigg\{ l \, \bigg| \, L^{(a^{(p)}_{j+1})}_p
\cap L^{(a^{(p)}_{j+1})}_l \neq \varnothing , a^{(p)}_{j+1} \leq l < p \bigg\}$ = $\left\{ a^{(p)}_i \, \big| \, 1 \leq i \leq j + 1 \right\}$

Let us consider two cases.
The first case is the case of $d(a^{(p)}_{j+1}) = d(a^{(p)}_j)$ ($\Leftrightarrow l_p = j+1$).
In this case, $I_{j+1} = \varnothing$.
Thus, the support of the remaining Dehn twists and isotopies to construct $\boldsymbol{L}^{(s)}$ for $s < a^{(p)}_{j+1}$ do not intersect with $L^{(a^{(p)}_{j+1})}_p$.
Thus we have $S_p = L^{(a^{(p)}_{j+1})}_p$, and $\{ l \, | \, S_p \cap S_l \neq \varnothing , l < p \} = \left\{ a^{(p)}_i \, \big| \, 1 \leq i \leq l_p \right\}$.

The second case is the case of $d(a^{(p)}_{j+1}) < d(a^{(p)}_j)$ ($\Leftrightarrow l_p > j+1$).
Since $I_{j+1} = \left\{ l \, \big| \, \widetilde{\xi_{j+1}}
\cap L_l \neq \varnothing , l < a^{(p)}_{j+1} \right\} = \left\{ l \, \big|
\, d(a^{(p)}_{j+1}) < l \leq d(a^{(p)}_j) = a^{(p)}_{j+2} \right\}$, we have $L^{(a^{(p)}_{j+2}+1)}_p = \cdots = L^{(a^{(p)}_{j+1})}_p$.
Moreover, by the construction, $L^{(a^{(p)}_{j+2}+1)}_p \cap K^\text{left}_{a^{(p)}_{j+2}}$ is a path of $<a^{(p)}_{j+1}>$ and $\bigg\{ l \, \bigg| \, L^{(a^{(p)}_{j+2}+1)}_p
\cap L^{(a^{(p)}_{j+2}+1)}_l \neq \varnothing , a^{(p)}_{j+2} < l < p \bigg\}$ = $\left\{ a^{(p)}_i
\, \big| \, 1 \leq i \leq j+1 \right\}$.
These conditions coincide with the formula which are obtained by replacing $j$ into $j+1$ in the first two conditions we assumed.

Finally, to prove Proposition \ref{prop:WhenIntersect}, we iterate the above procedure $l_p$-times.

\subsubsection{Proof of Proposition \ref{prop:OrderOfIntersection}}

Next, we study the order of intersections.
In this subsection, we study the Dehn twists and isotopies as above with
orientation of submanifolds.

First, we show that the subsets $\{ q_p \}$, $\{ q_{a^{(p)\dagger}_j, p }\}$,
$\{ q_{p, a^{(p)}_j} \}$ of $S_p$ appear in this order with respect to the orientation of $S_p$.
By the construction, we have that $q_{a^{(p)\dagger}_j, p} \in H_j$.
When one goes along $S_p$ from $q_p$, the first strip one goes through is $H_p$.
This shows that the subsets appear in the above order.
(Recall that the brane orientation of $L^{(p)}_p = L_p$ is the same with that induced from the orientations in Figure \ref{fig:L'(n-1)}.)
 
Next, we study the order of $q_{p, a^{(p)}_j}$'s.
Together with the orientation, the path $\widetilde{\xi_1}$ is a path from $l_p(0) = l_{a^{(p)}_0}(0)$ to a point near $l_{p-1}(0) =
l_{a^{(p)}_1}(0)$.
Thus, the points in  $L^{(p-1)}_p = L^{(a^{(p)}_1)}_p$ appear by the following
order: $q_p, d(p), d(p)-1, \dots , d(p-1)+1, p-1 \, ( = a^{(p)}_1)$.
Here, each number represents the intersection of $L^{(p-1)}_p$ and corresponding
submanifold.
Figure \ref{fig:xi'sandwidetildexi's} illustrates the situation.

Suppose that $l_p = 1 \, (\Leftrightarrow d(p) = d(p-1))$.
Then the order of the points is $q_p, a^{(p)}_1$.
Moreover, $L^{(p-1)}_p$ is away from the support of remaining Dehn twists
and isotopies.
Hence, we have the proof for the case $l_p = 1$.
(However, in fact, this case is trivial.)

Now, we consider the case $l_p \geq 2$.
The Dehn twist $\tau_{L_{d(p)}} = \tau_{L_{a^{(p)}_2}}$ acts on the path
$\widetilde{\xi_1}$ and obtain $\xi_2$.
We can see that the path $\widetilde{\xi_2}$ comes from a point near $l_{d(p)}(0)
= l_{a^{(p)}_2}(0)$, go to a point near $l_{p-1}(0) = l_{a^{(p)}_1}(0)$,
and intersects with $L^{(d(p))}_j = L^{(a^{(p)}_2)}_j$ by the following order of subscripts,
$d(d(p))+1, d(d(p))+2, \dots, d(p-1)$.
Hence, the points in  $L^{(a^{(p)}_2)}_p$ appear by the following
order: $q_p, d(p) \, (= a^{(p)}_2), d(d(p))+1 \, (=d(a^{(p)}_2) + 1), d(d(p))+2,
\dots, d(p-1) \, (= d(a^{(p)}_1)), p-1 \, (= a^{(p)}_1)$.
Figure \ref{fig:xi'sandwidetildexi's} illustrates the situation.

Suppose that $l_p = 2 \, (\Leftrightarrow d(a^{(p)}_1) = d(a^{(p)}_2))$.
Then, the order of the points is $q_p, a^{(p)}_2, a^{(p)}_1$.
By almost the same discussion as in the case of $l_p = 1$, we have proved
the case of $l_p = 2$.

Now, consider the case for $l_p \geq 3$.
By the same discussion, the intersection points in $L^{(a^{(p)}_3)}_p$ appear by the following
order: $q_p, a^{(p)}_2, d(a^{(p)}_2), d(a^{(p)}_2)-1, \dots, d(a^{(p)}_3)+1,a^{(p)}_3,
a^{(p)}_1$.
Again Figure \ref{fig:xi'sandwidetildexi's} illustrates the situation.

If $l_p = 3$, then we can finish the proof by the same argument, and if $l_p
> 3$, then we can finish the proof by the iteration of  the above discussion.

Next, we study the order of $q_{a^{(p)\dagger}_j, p}$'s.
We prove the statement about the order of them by induction on $p$.
For the case of $p = n$, the statement is trivial.
Now we assume that the statement is true for $p > s$ and prove the case of
$p = s$.

\begin{figure}[hbt]
\centering
\includegraphics[width=8cm]{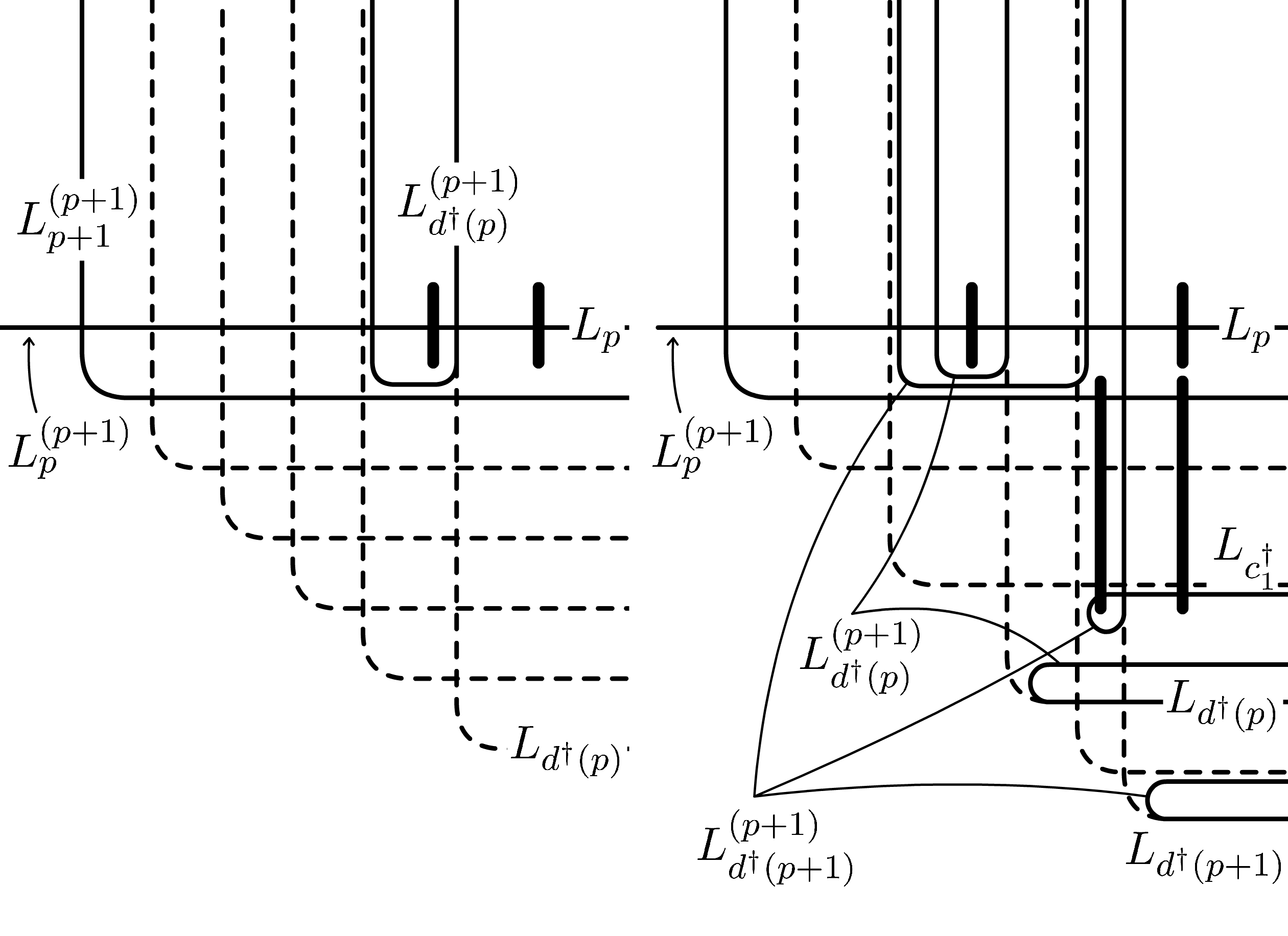}
\caption{$\boldsymbol{L}^{(s+1)}$}\label{fig:Lk+1}
\end{figure}

Let us see the case with small $l^\dagger_s$.
In the case of $l^\dagger_s = 1$, the statement is trivial.
In the case of $l^\dagger_s = 2$, there exists a relation $[s, d^\dagger(s)]_\mathbb{Z}$
and we have $a^{(s)^\dagger}_0 = s$, $a^{(s)^\dagger}_1 = s+1$, and $a^{(s)^\dagger}_2
= d^\dagger(s)$.
In this case, $\boldsymbol{L}^{(s+1)}$ is as in the left part of Figure
\ref{fig:Lk+1}.
The remaining Dehn twists and isotopies do not change the order of the intersection
points $q_{a^{(s)^\dagger}_j, s}$'s, so the statement for these points holds.

In the case of  $l^\dagger_s = 3$, then there exist two relations $[s, d^\dagger(s)]_\mathbb{Z}$
and $[c^\dagger_1, d^\dagger(s+1)]_\mathbb{Z}$ satisfying $s < c^\dagger_1
< d^\dagger(s)$.
The submanifold $L^{(s+1)}_{a^{(s)\dagger}_3}$ is isotopic to $\tau_{c^\dagger_1}\tau_{a^{(s)\dagger}_3-1}L_{a^{(s)\dagger}_3}$,
where $\tau_{j} = \tau_{L_j}$.
Hence, the intersection of $L^{(s+1)}_s$ and $L^{(s+1)}_{a^{(s)\dagger}_3}$
is just left (with respect to the orientation of $L_{c^\dagger_1}$) to the
intersection of $L^{(s+1)}_s$ and $L^{(s)}_{c^\dagger_1}$ as in the right
part of Figure \ref{fig:Lk+1}. 
By the construction of $\boldsymbol{L}^{(s)}$, the remaining Dehn twists
and isotopies do not change the order of the intersection
points $q_{a^{(s)^\dagger}_j, s}$'s, so the statement for this case holds.

\begin{figure}[hbt]
\centering
\includegraphics[width=8cm]{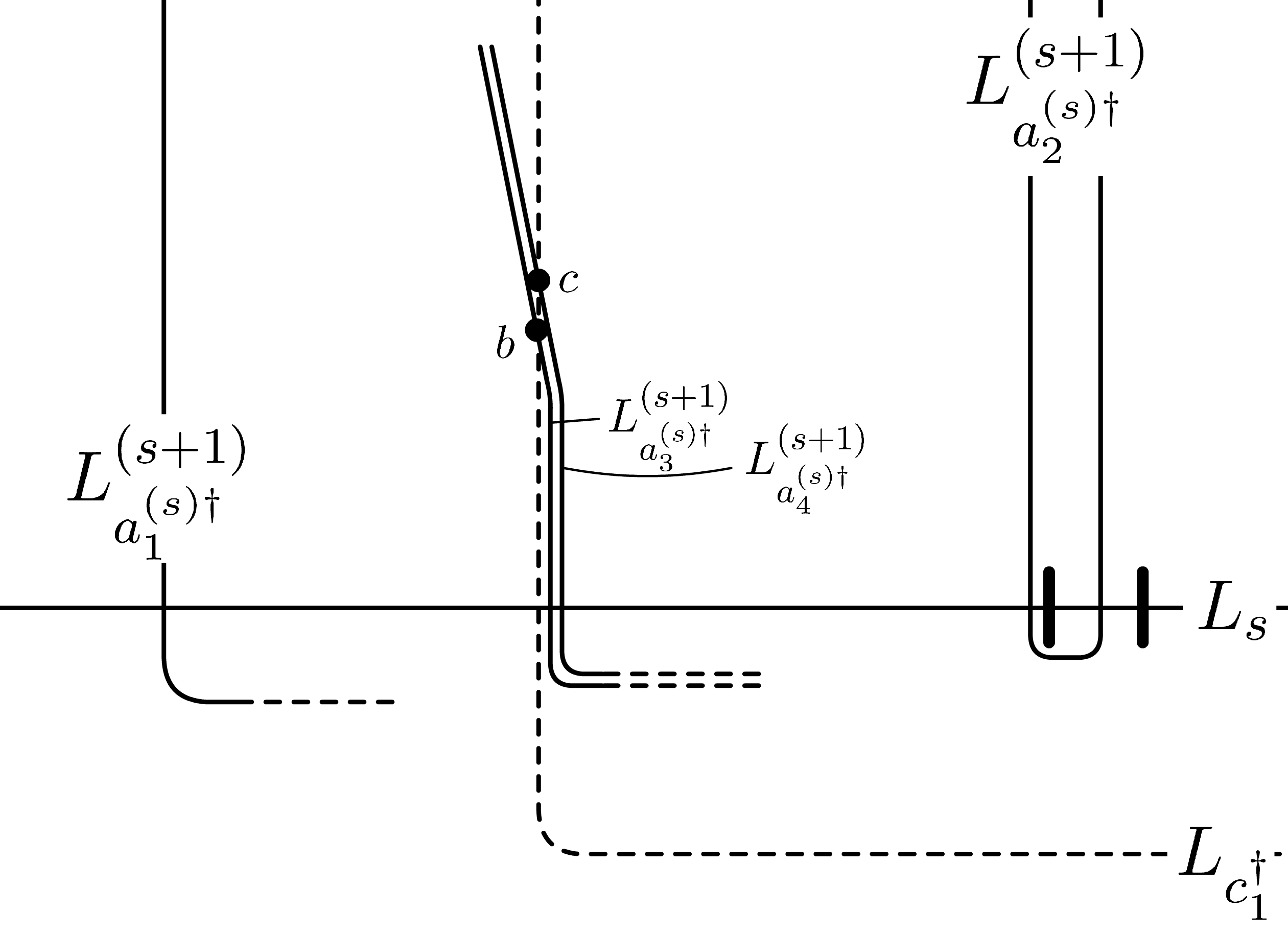}
\caption{The case when $c^\dagger_1 = d(d^{\dagger 2}(s)-1)$}\label{fig:CASE1}
\end{figure}

In the case of $l^\dagger_s = 4$, there exist three relations $[s, d^\dagger(s)]_\mathbb{Z}$,
$[c^\dagger_1, d^\dagger(s+1)]_\mathbb{Z}$, and $[b^\dagger_2, d^{\dagger
2}(s)]$ such that $s < c^\dagger_1 < d^\dagger(s)$ and $d^\dagger(s) \leq
b^\dagger_2 < d^\dagger(s+1)$.
The submanifold $L^{(s+1)}_{a^{(s)\dagger}_3}$ is same as in the case of
$l^\dagger_s = 3$ and $L^{(s+1)}_{a^{(s)\dagger}_4}$ is isotopic to $\tau_{d(a^{(s)\dagger}_4
- 1)}\tau_{b^\dagger_2}\tau_{a^{(s)\dagger}_4 - 1}L_{a^{(s)\dagger}_4}$.
Now, $d(a^{(s)\dagger}_4-1)$ satisfies the following inequality $c^\dagger_1
\leq d(a^{(s)\dagger}_4-1) < d^\dagger(s)$.
Here, the second inequality follows from Lemma \ref{lem:dddaggerlemma}.
The first inequality follows from the fact that $c^\dagger_1$ is the smallest
element in $S$ which is greater than $s$ and $s < d(d^{\dagger 2}(s)-1) =
d(a^{(s)\dagger}_4 - 1)$.
If $c^\dagger_1 < d(d^{\dagger 2}(s)-1)$, then the statement holds.
Now assume that $c^\dagger_1 = d(d^{\dagger 2}(s)-1)$.
By the assumption of the induction and the definition of $\boldsymbol{L}^{(s)}$,
three points $q_{c^\dagger_1}$, $c$, and $b$ in $L^{(c^\dagger_1)}_{c^\dagger_1}$
are located in this order, where $b$ is the unique intersection point with
$L^{(c^\dagger_1)}_{a^{(s)\dagger}_3}$ and $c$ is that with $L^{(c^\dagger_1)}_{a^{(s)\dagger}_4}$,
since $a^{(s)\dagger}_3 = a^{(c^\dagger_1)}_2$ and $a^{(s)\dagger}_4 = a^{(c^\dagger_1)}_3$.
Hence, as in Figure \ref{fig:CASE1}, the statement in this case also holds.

\begin{figure}[hbt]
\centering
\includegraphics[width=8cm]{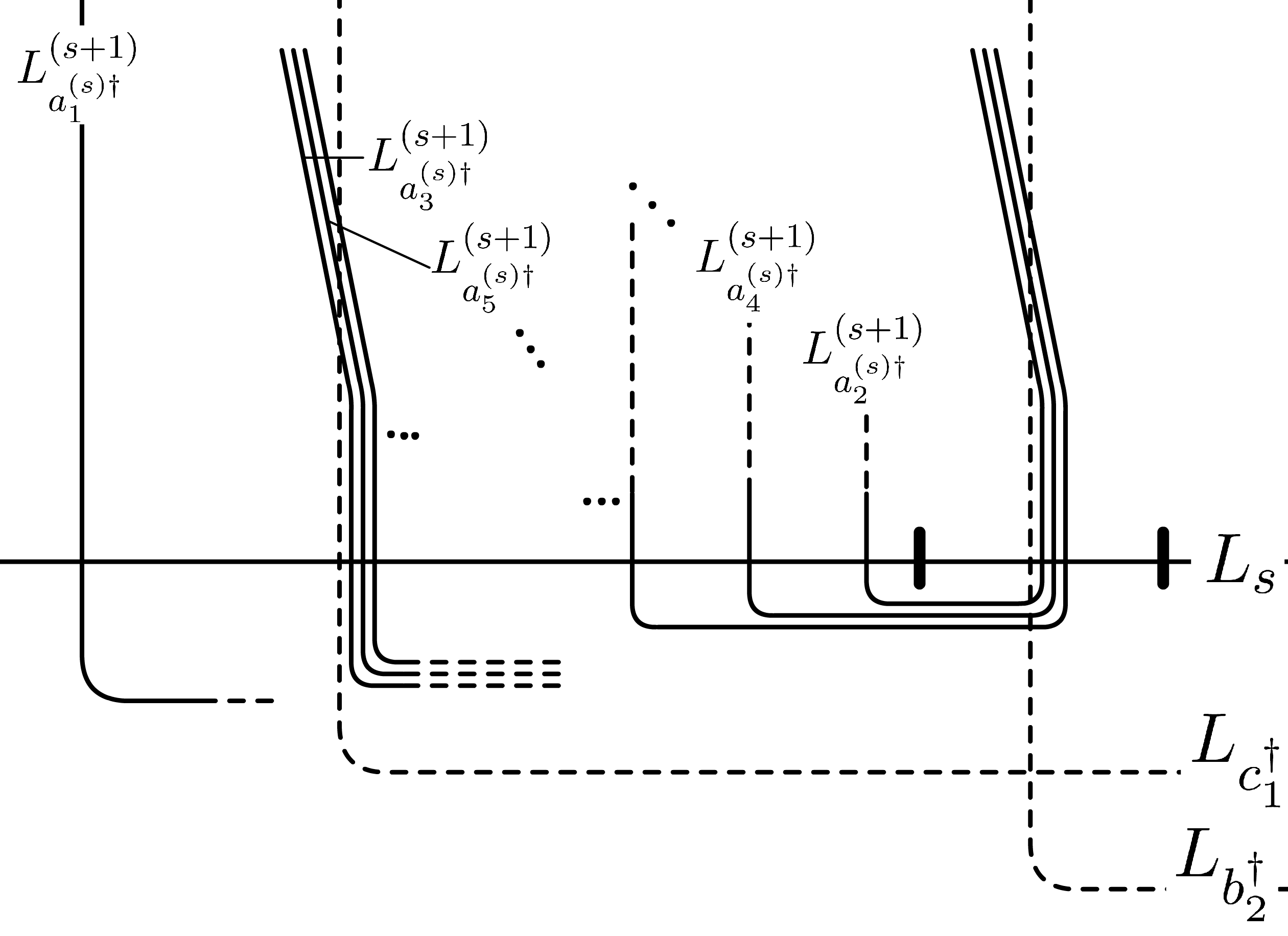}
\caption{The case when $l^\dagger_s \geq 5$}\label{fig:ldaggerdGEQ5}
\end{figure}

Finally, we consider the case when $l^\dagger_s \geq 5$.
As in the discussion above, there exist  relations $[s, d^\dagger(s)]_\mathbb{Z}$,
$[c^\dagger_1, d^\dagger(s+1)]_\mathbb{Z}$, and $[b^\dagger_2, d^{\dagger
2}(s)]$ with the same condition.
For odd number $2j + 1$ with $j \geq 1$, $L^{(c^\dagger_1)}_{a^{(s)\dagger}_{2j+1}}$
intersects with $L^{(c^\dagger_1)}_{c^\dagger_1}$ since $a^{\left(a^{(s)\dagger}_{2j+1}\right)}_{2j}
= d^jd^{\dagger j}(s+1) = c^\dagger_1$.
By the hypothesis of the induction and the construction of $\boldsymbol{L}^{(s)}$,
we can deduce the points $q_s, q_{a^{(s)\dagger}_1, s}, q_{a^{(s)\dagger}_3,
s}, q_{a^{(s)\dagger}_5, s}, \dots$ in $S_s$ are in this order along the
orientation of $S_s$ as in the left side of Figure \ref{fig:ldaggerdGEQ5}.

For even number $2j$ with $j \geq 1$, $L^{(b^\dagger_2)}_{a^{(s)\dagger}_{2j}}$
intersects with $L^{(b^\dagger_2)}_{b^\dagger_2}$ since $a^{(a^{(s)\dagger}_{2j})}_{2j-2}
= d^{j-1}d^{\dagger j}(s) = dd^{\dagger 2}(s) = b^\dagger_2$.
Additionally, $L^{(b^\dagger_2)}_{a^{(s)\dagger}_{2j}}$ intersects with $L_{r_j}$
with $r_j = d^{j-1}(a^{(s)\dagger}_{2j}-1) = d^{j-1}(d^{\dagger j}(s) - 1)$.
By Lemma \ref{lem:dddaggerlemma} and discussion of the case with $l^\dagger_s
= 4$, we have $c^\dagger_1 \leq \cdots \leq r_2 \leq r_1 < d^\dagger(s)$.
Again by the same argument as in the case $l^\dagger_s = 4$, we can show
that the following points are in this order: 
$q_s,
q_{a^{(s)\dagger}_1, s}, q_{a^{(s)\dagger}_3, s},
\dots , q_{a^{(s)\dagger}_{l^\dagger_s}, s},\dots ,
q_{a^{(s)\dagger}_4, s}, q_{a^{(s)\dagger}_2, s}$.
as in the right side of Figure \ref{fig:ldaggerdGEQ5}.
This completes the proof of Proposition \ref{prop:OrderOfIntersection}.

\clearpage

\subsection{Core of $\boldsymbol{S}^\maltese$}

Recall that the sets of start and goal points of relations are $S = \{ s_1,
s_2, \dots , s_m \}$ and $T = \{ t_1, t_2, \dots , t_m \}$.
We write the interval in $S_j$ between $q_{j, i}$ and $q_{j, l}$ contained
in the core of $\boldsymbol{S}$ by $[i, l]_j$.
(We don't care that which number is bigger.)

\begin{lem}\label{lem:basicPolygon}

For any $1 \leq j \leq m$, the intervals $[t_j, s_j+1]_{s_j}$, $[s_j, s_j
+ 2]_{s_j + 1}, \dots , [t_j - 2, t_j]_{t_j-1}$, $[t_j - 1, s_j]_{t_j}$ bounds
a $(t_j - s_j + 1)$-gon.

\end{lem}

\begin{proof}

\begin{figure}[hbt]
\centering
\includegraphics[width=8cm]{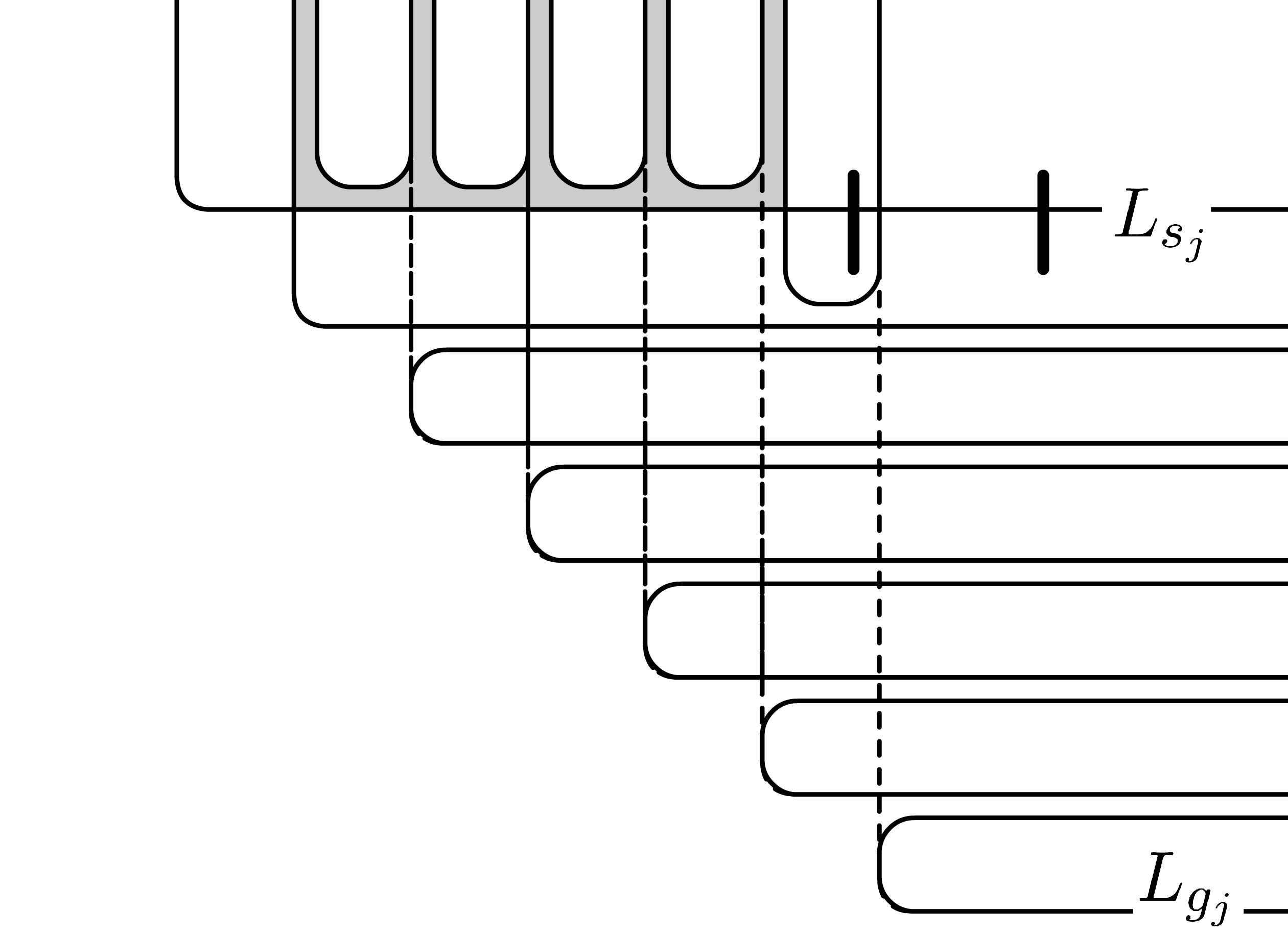}
\caption{the $(t_j - s_j + 1)$-gon}\label{fig:d+1-gon}
\end{figure}

By the constuction, the submanifolds $L^{(s_j + 1)}_i$ for $s_j \leq i \leq
t_j$ intersect as in Figure \ref{fig:d+1-gon}. Hence, the above submanifolds bound a $(t_j - s_j +1)$-gon.
Again by the construction, the remaining Dehn twists and isotopies do not
destroy the $(t_j - s_j + 1)$-gon.
This completes the proof.
\hfill $\Box$
\end{proof}

We write the above $(t_j - s_j + 1)$-gon by $Y_j$ and set $\displaystyle
X_j \coloneqq \bigcup_{i \leq j} cS_i \, \cup \, \bigcup_{t_i \leq j} Y_i$,
where $cS_i$ is the core of $S_i$.

\begin{lem}

$X_j$ is contractible.

\end{lem}

\begin{proof}

We prove this by induction.
For $j = 0$, This is true since $X_0 = cS_0$ is just an interval or a point.
Now we prove that $X_p$ is contractible under the condition that $X_{p-1}$
is contractible.

If $p \neq t_i$ for any $1 \leq i \leq m$, then $X_p = X_{p-1} \cup cS_p$.
In this case, since $d(p) = d(p-1)$ so we have $l_p = 1$. Thus,
$cS_p$ intersects with $X_{p-1}$ at just one point $q_{p, p-1} \in cS_{p-1} \subset X_{p-1}$. Hence
we have that $X_p$ is contractible.

Next, suppose that $p = t_i$ for some $i$.
Then, $X_p$ is the union $X_{p-1} \cup cS_p \cup Y_i$.
In this case, it is sufficient to show that $\partial Y_i \cap X_{p-1}$ is
connected because it implies that $X_p$ is the union of contractible set $X_{p-1}$ and contractible set $Y_i$ with contractible intersection.

By Proposition \ref{prop:OrderOfIntersection}, we have $cS_j = [j-1, j+1]_j$.
Now, we have $\displaystyle \partial Y_i \cap X_{p-1} = [s_i + 1, t_i]_{s_i}
\cup \bigcup_{s_i < j < t_i = p}[j-1, j+1]_j \cup ([p - 1, p + 1]_p \cap
X_{p-1})$. 
Therefore it suffices to show that $([p - 1, p + 1]_p \cap X_{p-1}) \setminus \{
q_{p, p-1} \}$ is connected.

If $l_p = 2$, then we have $[p - 1, p + 1]_p \cap X_{p-1} = \{ q_{p,
p-1}, q_{p, a^{(p)}_2} \}$.
This follows from the following discussion.
The intersection $cS_p \cap (\cup_{j < p} \, cS_j)$ is these two points hence
the interval $[p-1, a^{(p)}_2]_p$ is either contained in $X_{p-1}$ or just
intersects with $X_{p-1}$ only at the two end points.
Since $q_{p, p-1}$ and $q_{p+1, p}$ are away from $\cup_{t_i < p} Y_i$ by Proposition \ref{prop:OrderOfIntersection}
and Lemma \ref{lem:basicPolygon}, we have $[p- 1, p + 1]_p \cap
X_{p-1} = \{ q_{p,
p-1}, q_{p, a^{(p)}_2} \}$.

Now, we consider the case $l_p \geq 3$. In this case, there exist relations
$[d(p), p]_\mathbb{Z}$ and $[d(p-1), c_1]_\mathbb{Z}$ with $d(p) < c_1 <
p$.
Let us write $c_1 = t_{i'}$. Since $t_{i'} = c_1 < p$, we have $Y_{i'} \subset X_{p-1}$.
Let $\gamma'$ be a path from $q_{p, p-1}$ to $q_{p, a^{(p)}_2}$ along $S_p$
with its image $[p-1, a^{(p)}_2]_p$.
Consider a path $\gamma$ which is an extension of $\gamma'$ so that $\gamma$
go beyond $q_{p, a^{(p)}_2}$ a little along $S_p$.
Then $\gamma$ starts from $q_{p, p-1}$ first intersect with $\cup_{j < p} \, cS_j$ at $q_{p, a^{(p)}_3}$,
and the final intersection point with that set is $q_{p, a^{(p)}_2}$.
Because of the following three facts: $q_{p, p-1}$ is away from $Y_{i'}$, the intersection of $S_p$ and $\partial Y_{i'}$
is transversal, and $\{ q_{p, a^{(p)}_2}, q_{p, a^{(p)}_3} \} = S_p \cap \partial Y_{i'}$ by Proposition \ref{prop:OrderOfIntersection},
we can conclude that $\gamma$ go into the interior of $Y_{i'}$ at $q_{p, a^{(p)}_3}$.
If $l_p = 3$, $\gamma$ must stay in $Y_{i'}$ unless it reaches the point $q_{p,
a^{(p)}_2}$ so we have $([p - 1, p + 1]_p \cap X_{p-1}) \setminus \{
q_{p, p-1} \} = [a^{(p)}_2, a^{(p)}_3]_p$ and this is connected.

Finall,  the case of $l_p \geq 4$, we consider
the intersection points $q_{p, a^{(p)}_r}$ of $\gamma$ and
$\cup_{j < p} \, cS_j$ with $r \geq 4$.
Since $a^{(p)}_r < a^{(p)}_3 = s_{i'}$, these points are contained in $\mathring{Y}_{i'}$. Therefore again $\gamma$ must stay in $Y_{i'}$ unless it reaches the point $q_{p,
a^{(p)}_3}$ so we have $([p - 1, p + 1]_p \cap X_{p-1}) \setminus \{
q_{p, p-1} \} = [a^{(p)}_2, a^{(p)}_3]$ and this is connected.

\hfill $\Box$
\end{proof}

\begin{rem}

By the above lemmas and propositions, the core $cS$ of $S$ has the following
properties.

\begin{itemize}
\item There exists a $(t_j - s_j + 1)$-gon $Y_j$ corresponds to a relation $[s_j, t_j]_\mathbb{Z}$.
\item The root and intersection points in $S_j$ is distributed by the order
displayed in Proposition \ref{prop:OrderOfIntersection}.
\item Especially, the core satisfies $cS_j = [j-1, j+1]_j$ for $0 < j < n$.
\item If $l^\dagger_0 \geq 2$, then $cS_0 = [1, d^\dagger(0)]_0$, else $cS_0
= \{ q_{1, 0} \}$.
\item If $l_n \geq 2$, then $cS_n = [d(n), n-1]_n$, else $cS_n = \{ q_{n.
n-1} \}$
\item For any $0 \leq p \leq n$, $\displaystyle X_p = \bigcup_{0 \leq j \leq
p}cS_j \cup \bigcup_{1 \leq t_i \leq p} Y_i$ is contractible.
\end{itemize}

In many concrete examples, such a ``diagram" is unique up to ``isotopy", and the author couldn't find any counter examples. Therefore, the author believes that the
uniqueness holds under some justification.
However, we don't go into this direction.

\end{rem}

\subsection{Determination of degree}

\begin{lem}\label{lem:Deg1Mor}

The cohomology group $H(\hom_\mathcal{F}(S^\#_j, S^\#_{j-1}))$ is  of one-dimension
and  concentrated in degree one part.
 
\end{lem}

\begin{proof}

Since our Lagrangian branes are Hamiltonian isotopic to  $S'{}^\#_j = \tau_{L^\#_1}
\cdots \tau_{L^\#_{j-1}}L^\#_j$
and $S'{}^\#_{j-1} = \tau_{L^\#_1} \cdots \tau_{L^\#_{j-2}}L^\#_{j-1}$, it
is
enough to study the cohomology of $\hom_\mathcal{F}(L^\#_j, L^\#_{j-1})$
by applying $(\tau_{L^\#_1}\tau_{L^\#_2} \cdots \tau_{L^\#_{j-1}})^{-1}$.
 
Recall that $L_{j-1}$ and $L_j$ intersects only at $p_{j-1, j}$ and its degree
is zero so we have $-1 < \alpha_j(p_{j-1, j}) - \alpha_{j-1}(p_{j-1, j})
< 0$.
Thus we have $0 < \alpha_{j-1}(p_{j-1, j}) - \alpha_j(p_{j-1, j}) < 1$ and
hence $p_{j-1, j} \in \hom_\mathcal{F}^1(L^\#_j, L^\#_{j-1})$.
\hfill $\Box$

\end{proof}

\begin{lem}\label{lem:DegJMor}

$q_{p, a^{(p)}_j} \in \hom^j_\mathcal{F}(S^\#_p, S^\#_{a^{(p)}_j})$.

\end{lem}

\begin{proof}

We prove this lemma by induction on $j$.
The first case $j = 1$ is proved in Lemma \ref{lem:Deg1Mor}.

Suppose that the statement is true for the case $j = i-1$.
Consider a sequence of intervals $[a^{(p)}_i, a^{(p)}_{i-1}]_p$, $[p, a^{(p)}_{i-1}
- 1]_{a^{(p)}_{i-1}}$, $[a^{(p)}_{i-1}, a^{(p)}_{i-1} - 2]_{a^{(p)}_{i-1}
- 1}$, \dots , $[a^{(p)}_{i} + 1, p]_{a^{(p)}_{i}}$.
These intervals form a loop $\gamma$ and in fact this loop does not have self intersections since there is no relation corresponds to an interval contained in $[a^{(p)}_i, a^{(p)}_{i-1}]_\mathbb{Z}$ by the definition of $\{ a^{(p)}_j \}$.
We can show that this loop $\gamma$ go left at every corner because of the orientations of the inter¥vals and the degrees of intersections.
Moreover, since $X_n$ is contractible, $\gamma$ bounds a $(a^{(p)}_{i-1} - a^{(p)}_i
+ 2)$-gon.
This shows that
$\mathcal{M}(q_{a^{(p)}_i + 1, a^{(p)}_i}, q_{a^{(p)}_i+ 2 , a^{(p)}_i + 1}, \dots , q_{p, a^{(p)}_{i-1}}; q_{p, a^{(p)}_i}) \neq \varnothing$
and thus we have $\deg (q_{a^{(p)}_i}) = i.$

\hfill $\Box$
\end{proof}

\subsection{Counting discs}

In this subsection, we prove the following proposition:

\begin{prop}\label{prop:ExistDisc}

Assume that the integers $i_0 < i_1 < \cdots < i_l$ satisfy the following
condition: $\hom_\mathcal{F}(S^\#_{i_j}, S^\#_{i_{j-1}}) \neq 0$ for any $j$, $\hom_\mathcal{F}(S^\#_{i_l},
S^\#_{i_0}) \neq 0$, and $|q_{i_1, i_0}| + |q_{i_2, i_1}| + \cdots + |q_{i_l, i_{l-1}}|
+ (2-l) = |q_{i_l, i_0}|$.
Then, there exists just one $(l + 1)$-gon contributes to the higher
composition and $\mu^l(q_{i_1, i_0}, q_{i_2, i_1}, \dots , q_{i_l, i_{l-1}})
= (-1)^{|q_{i_l, i_{l-1}}| \, ( \, |q_{i_l, i_0}| + 1)} q_{i_l, i_0}$ holds.

\end{prop}

Before we begin to prove this proposition, we see the following lemma:

\begin{lem}\label{lem:NoBackwardDisc}

For any $i_0 < i_1 < \cdots < i_l$ $(l \geq 1)$, there is no clockwise null-homotopic
embedded loop $\gamma$ such that it goes along an interval contained in $cS_{i_l}$,
turns left or right and goes along an interval in $cS_{i_{l-1}}$, turns left
or right and goes along an interval in $cS_{i_{l-2}}$, iterate this procedure for $i_{l-3}$, $i_{l-4}, \cdots , i_1$, and finally goes along $cS_{i_0}$ and comes back to the start point.

\end{lem}

\begin{figure}[hbt]
\centering
\includegraphics[width=8cm]{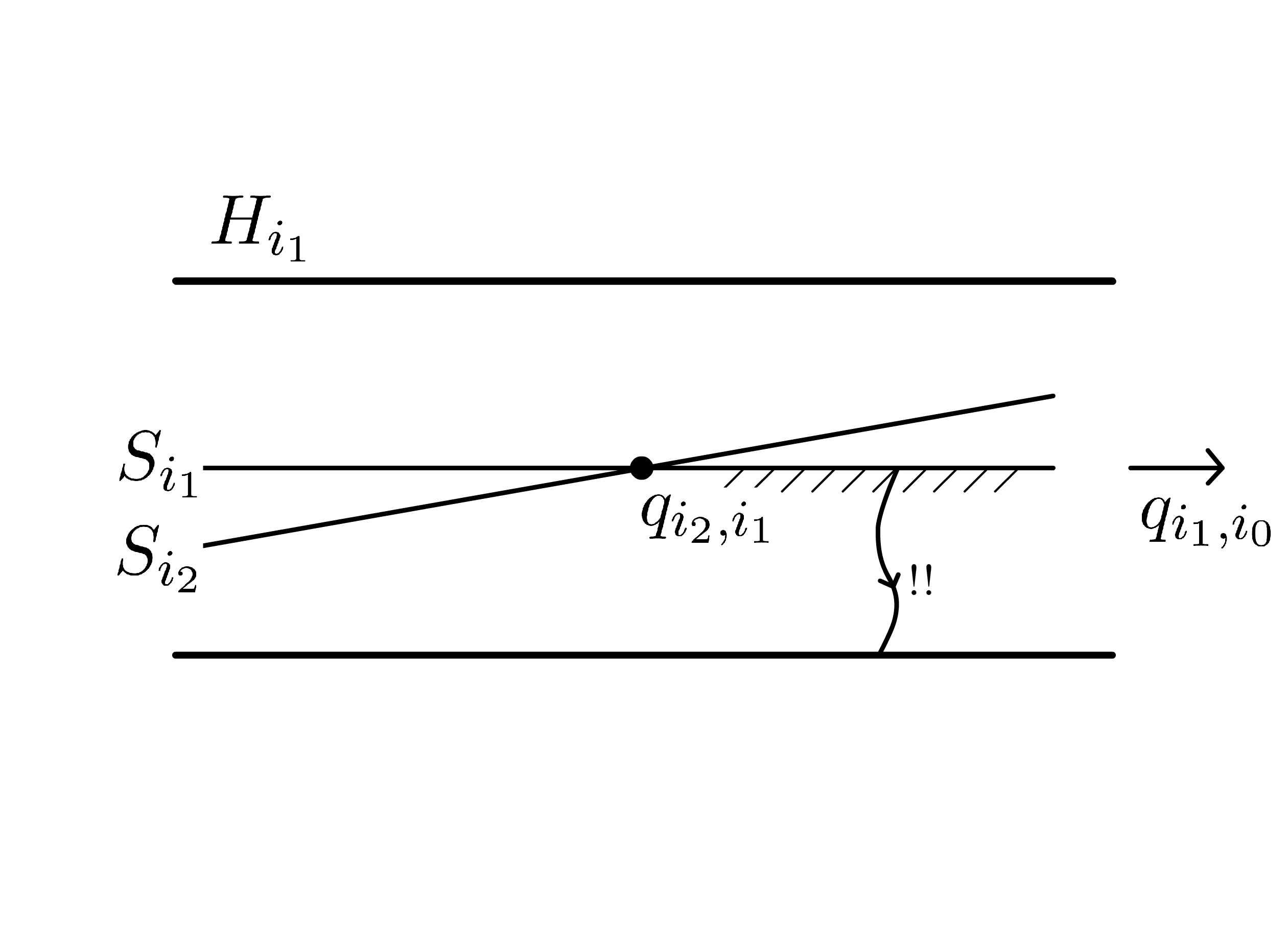}
\caption{Figure of $H_i$}\label{fig:Si's}
\end{figure}

\begin{proof}

In the case of $l = 1$, there is no  bigon since $S_{i_0} \cap S_{i_1}$ is
either $\varnothing$ or singleton.

Now we consider the case $l \geq 2$.
Assume contrarily that there exists a loop $\gamma$ as in the statement.
By definition, we have $q_{i_2, i_1} \in H_{i_1}$. By Proposition \ref{prop:OrderOfIntersection},
the orientation of $\gamma$ and that of $S^\#_{i_1}$ coincide.
Since $\gamma$ is a clockwise null-homotopic loop, $\gamma$ bounds a disc
in its right side.
However, it is impossible by the same argument of Lemma \ref{lem:TailAndNoDisc}
since $S_{i_0} \cap H_{i_1} = \varnothing$ (this is because $i_0 < i_1$) as depicted in 
Figure \ref{fig:Si's}.
\hfill $\Box$
\end{proof}

Now, we start the proof of the Proposition \ref{prop:ExistDisc}.
For integers $i_0 < i_1 < \cdots < i_l$ as in the proposition, we consider
a loop $\gamma$ which first starts from $q_{i_1, i_0}$, goes to $q_{i_2,
i_1}$ along $cS_{i_1}$, and at every corner, $\gamma$ turns left.
Obtain $\gamma'$ by perturbing $\gamma$ so that $\gamma'$ is smooth, and
we write its image by $L_\gamma$.
By the degree condition, we can compute the writhe of $L_\gamma$ as
$w(L_\gamma) = 2$.

\begin{figure}[hbt]
\centering
\includegraphics[width=8cm]{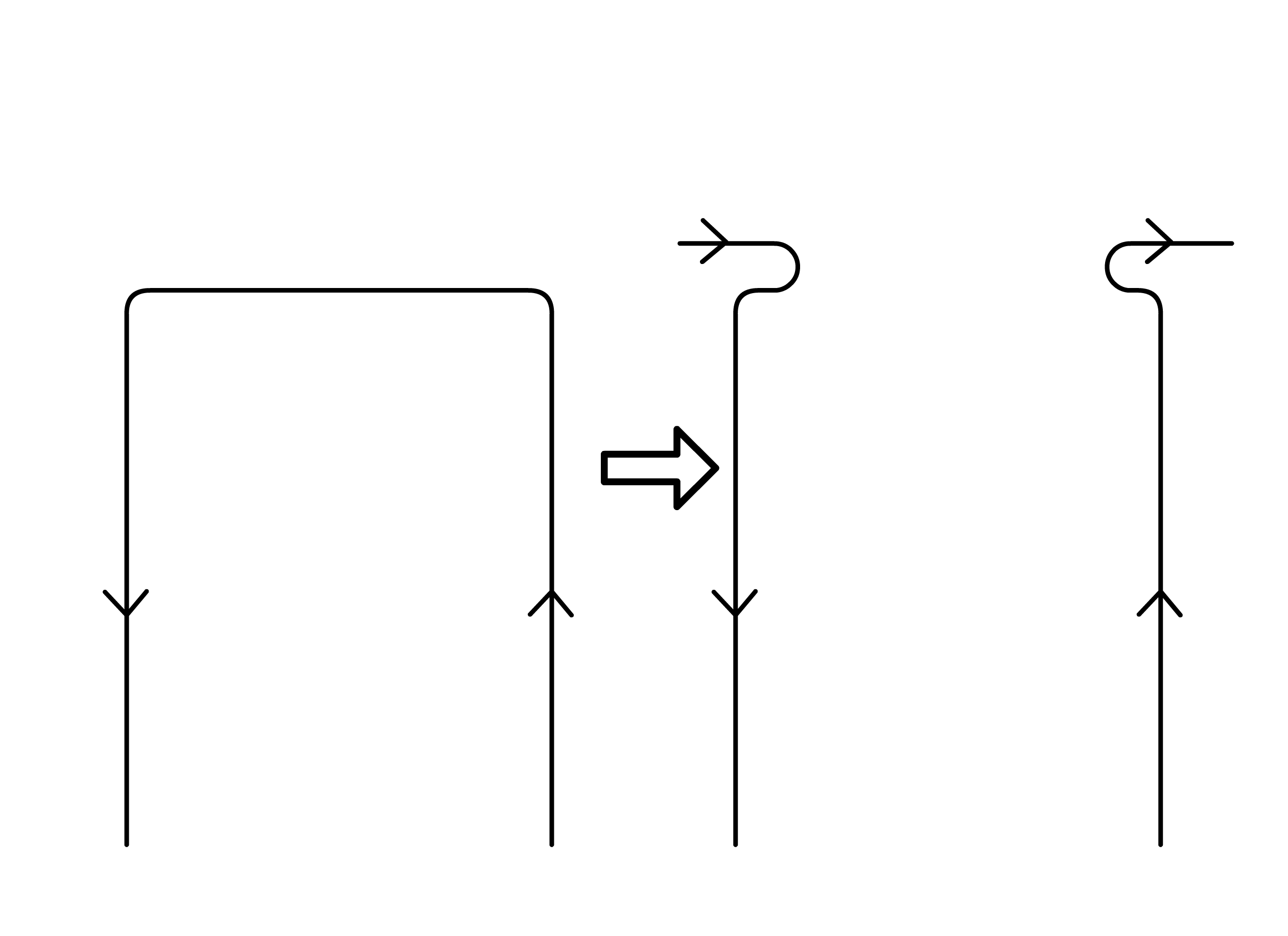}
\caption{making right turn into left turn}\label{fig:flipturns}
\end{figure}

First, we prove that $\gamma \subset c\boldsymbol{S}$.
Assume contrarily i.e. $I = \{ i_j \, | \, \gamma \cap S_{i_j} \not \subset
cS_{i_j}\}$ is non-empty.
We change the way of turning of $\gamma$ around $S_i$ for $i \in I$, so that
the resulting curve $\widetilde{\gamma}$ is contained in $c\boldsymbol{S}$.
At this constructoin, we change left turns into right turns as in Figure \ref{fig:flipturns}
in even times, so we have $w(L_{\widetilde{\gamma}}) \leq 0$, where $L_{\widetilde{\gamma}}$
is the image of smoothing $\widetilde{\gamma}'$ of $\widetilde{\gamma}$.

\begin{figure}[hbt]
\centering
\includegraphics[width=8cm]{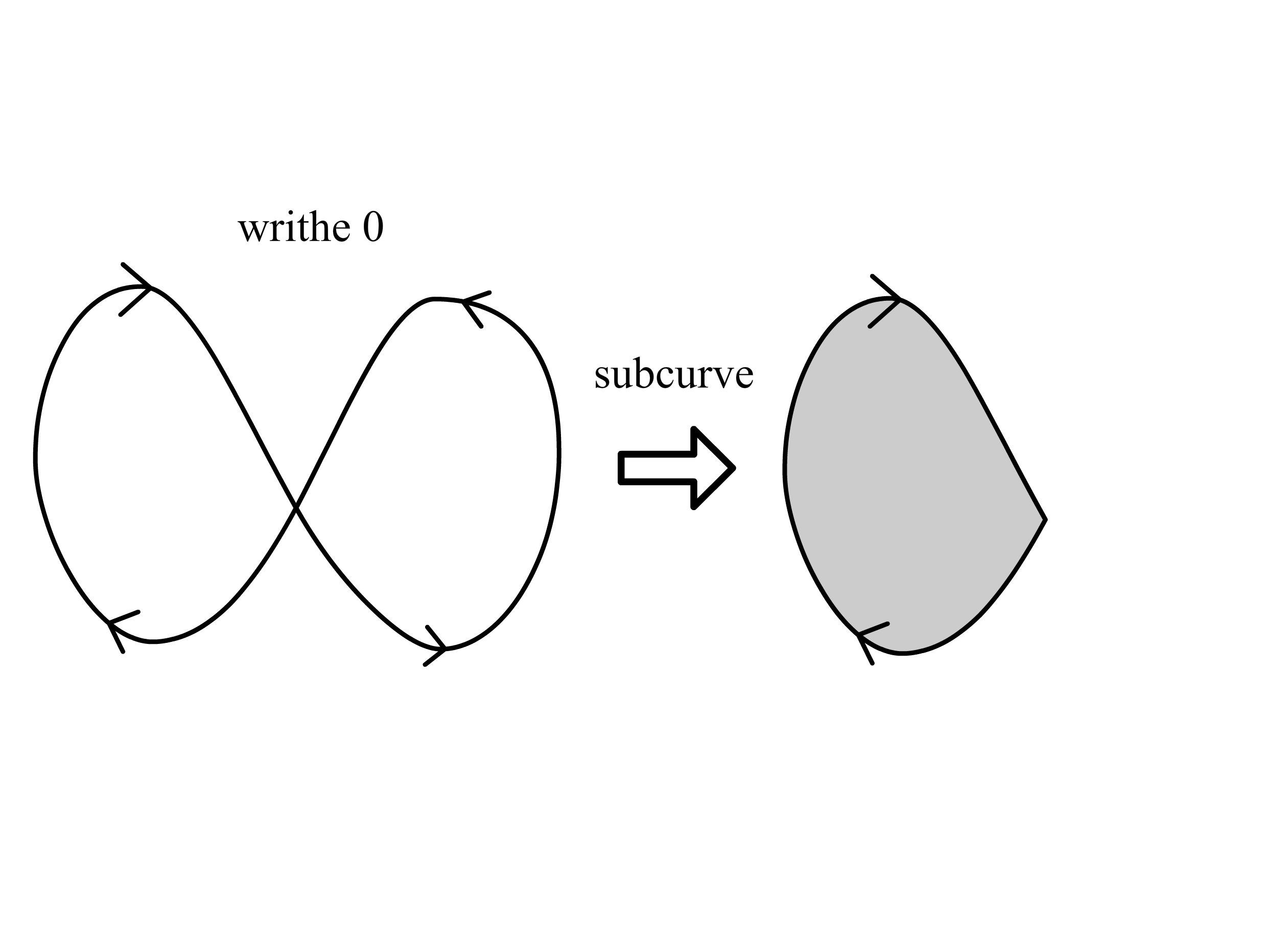}
\caption{Example of finding clockwise loop}\label{fig:clockwiseLoop}
\end{figure}

Now, $\widetilde{\gamma}$ is a piecewise smooth immersed curve in a contractible
region $X_n$ in $M$.
Moreover, all the self-intersection is transitive.
Hence, we can show that there exists sub curve which bounds a disc in its
right side. (The easiest case is presented in Figure \ref{fig:clockwiseLoop}.)
Therefore this contradicts with Lemma \ref{lem:NoBackwardDisc}.
Thus, we have shown that $\gamma \subset c\boldsymbol{S} \subset X_n$.

Now, $\gamma$ is a writhe two, piecewise smooth, embedded curve in contractible region $X_n$ (the lack of the self-intersection is deduced by the following discussion: if it has
self-intersections, there must exist a clockwise subloop).
Hence, $\gamma$ bounds a disc $u$ in its left side.

Since  $\boldsymbol{S}^\pm$ is a perfect collection of tailed Lagrangian submanifolds
and each pair of $S_j$'s intersect at most once, the number of such discs
is at most one. Hence the disc above is the unique element in $\mathcal{M}(q_{i_1,
i_0}, q_{i_2, i_1}, \dots , q_{i_l, i_{l-1}})$.

Finally, we study the sign $(-1)^{s(u)}$ of $u$.
By the definition of $\boldsymbol{S}^\maltese$, all switching points are
irrelevant to the sign.
As in the proof of Lemma \ref{lem:NoBackwardDisc}, the brane orientation
of $S^\maltese_{i_j}$ and $\partial u$ coincide for $0 < j < l$.
By Proposition \ref{prop:OrderOfIntersection},  the orientations of $S^\maltese_{i_l}$ and $\partial
u$ coincide if and only if $|q_{i_l, i_{l-1}}|$ is odd ($\Leftrightarrow$
$j$ in $i_{l-1} = a^{(i_l)}_j$ is odd).
Thus the sign is $(-1)^{s(u)} = (-1)^{(|q_{i_l, i_l-1}| + 1)(|q_{i_l, i_{l-1}}|
+ |q_{i_l, i_0}|)} = (-1)^{(|q_{i_l, i_l-1}| + 1)|q_{i_l, i_0}|}$.
This completes the proof of Proposition \ref{prop:ExistDisc}.

Now, what we have proved is that $\mathcal{F}_{S, T}^\to(\boldsymbol{S}_{S,
T}^\maltese)$ and $\mathcal{B}_{S, T}$ are isomorphic.
Therefore, this completes the proof of Theorem \ref{thm:MAIN1TheFormulaOfAInftyKoszulDual}.

\clearpage

\subsection{Some examples}\label{subsec:SomeExamples}

We see some examples of the core of $\boldsymbol{S}$.
First, we prepare an algebra which we compute its $A_\infty$-Koszul
dual.
We compute $A_\infty$-Koszul duals of $A_1 \coloneqq k(\ora\Delta_3,
(\alpha_3\alpha_2\alpha_1))$,
$A_2 \coloneqq k(\ora\Delta_6, (\alpha_3\alpha_2\alpha_1, \alpha_6\alpha_5\alpha_4))$,
and $A_3 \coloneqq k(\ora\Delta_6, (\alpha_3\alpha_2\alpha_1, \alpha_4\alpha_3,
\alpha_6\alpha_5\alpha_4))$.
Here, $\ora{\Delta}_n$ is the directed $A_n$-quiver.
We distinguish the relevant items like exact Riemann surface $M$, a
collection of Lagrangian submanifolds $\boldsymbol{L}$, and so on by giving
them subscripts like $M_i$ and $\boldsymbol{L}_i$.

First, we study the core of $\boldsymbol{S}_1$.
This collection consists of $S_3 \simeq \tau_0 \tau_1 \tau_2 L_3$, $S_2 \simeq \tau_0 \tau_1 L_2$, $S_1 \simeq \tau_0 L_1$, and $S_0 \simeq L_0$.
(Here, the symbol ``$\simeq$" represents that both sides are Hamiltonian isotopic.)
Now, we just want to investigate their    intersections and polygons, we apply them $\tau_0^{-1}$
and consider $S'_3 \simeq \tau_1 \tau_2 L_3$, $S'_2 \simeq \tau_1 L_2$, $S'_1 \simeq L_1$,
and $S'_0 \simeq L_0$.

\begin{figure}[hbt]
\centering
\includegraphics[width=8cm]{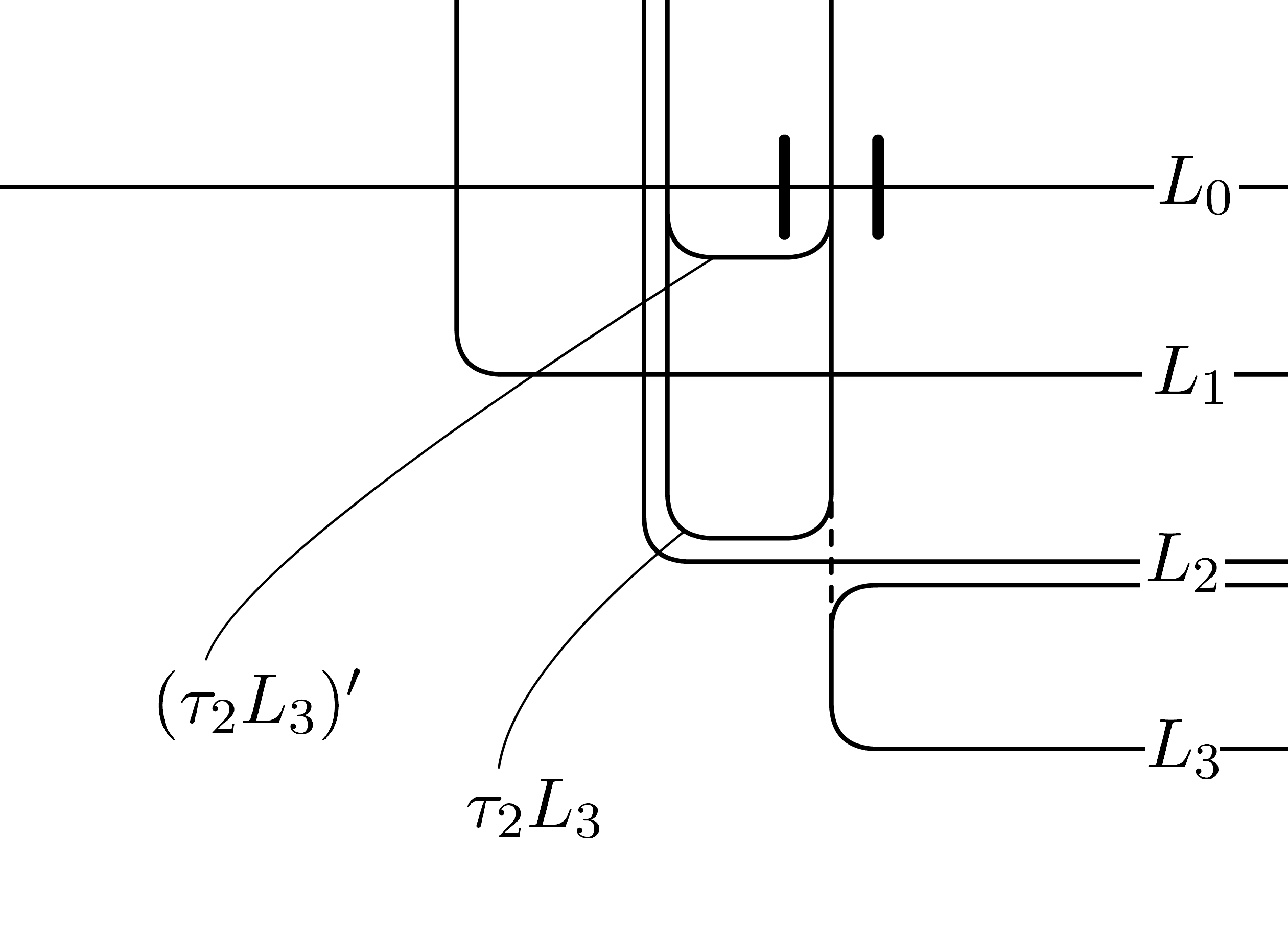}
\caption{$\tau_2 L_3$ and $(\tau_2 L_3)'$}\label{fig:S1_1}
\end{figure}

Now, we consider $\tau_2 L_3$. This is drawn in Figure \ref{fig:S1_1}. We
can see that $\tau_2L_3$ intersects with $L_1$ twice, but we can eliminate
the intersections by an isotopy. We write the resulting curve by $(\tau_2
L_3)'$.
Since $(\tau_2 L_3)' \cap L_1 = \varnothing$, we can define $S'_3 = (\tau_2
L_3)'$.
Next, we consider $\tau_1 L_2$.
This curve again intersects with $L_0$ twice so we eliminate the intersections
by the same way, obtain $(\tau_1 L_2)'$ and replace $S'_2$ by this curve.

\begin{figure}[hbt]
\centering
\includegraphics[width=8cm]{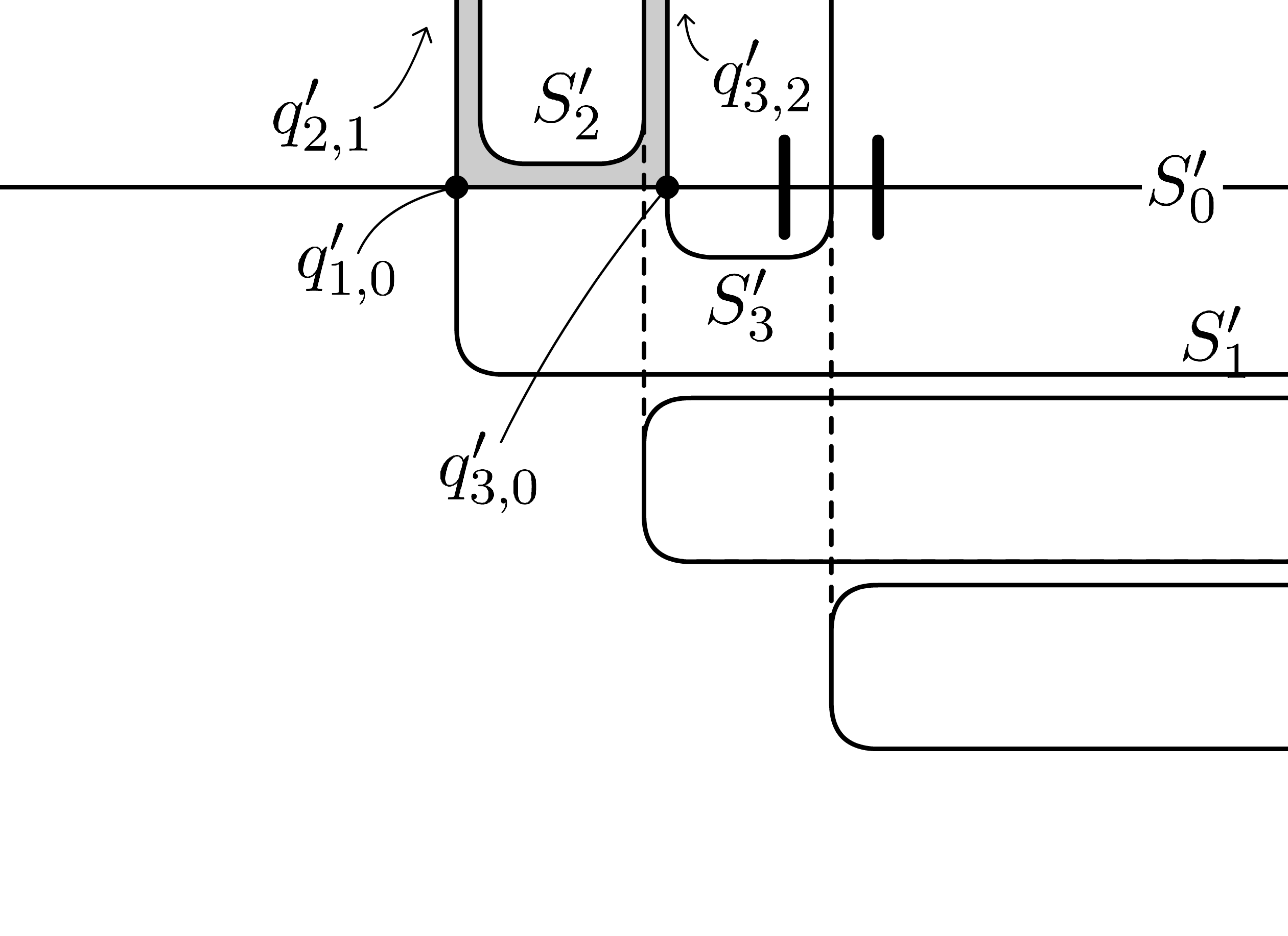}
\caption{$\boldsymbol{S}_1$}\label{fig:S1_2}
\end{figure}

Then, the resulting $\boldsymbol{S}'_1$ is as in Figure \ref{fig:S1_2}.
We can see that there are four points $q'_{i, j} \in S'_i \cap S'_j$ for
$(i, j) = (1, 0), (2, 1), (3, 2), (3, 0)$ and quadrangle with the four vertices.
In fact, there is no polygon other than this quadrangle, so the (boundary
of) quadrangle
is nothing but the core of $\boldsymbol{S}_1'$.
By the same computation of Dehn twists, we can see that there emerges one
$(l
+ 1)$-gon for every relation of length $l$.
Together with the degree, a directed $A_\infty$-category
$\mathcal{B}_1 \coloneqq \mathcal{F}_1^\to(\boldsymbol{S}^\#_1) \cong \mathcal{F}_1^\to(\boldsymbol{S}'_1{}^\#)$
can be represented as follows: $\hom_{\mathcal{B}_1}^d(S^\#_i, S^\#_j) =
0$ except for $\hom_{\mathcal{B}_1}^0(S^\#_j, S^\#_j) =k \cdot 1_{S^\#_j}$,
$\hom_{\mathcal{B}_1}^1(S^\#_j, S^\#_{j-1}) =\ k \cdot q_{j, j-1}$, and $\hom_{\mathcal{B}_1}^2(S^\#_3,
S^\#_0) = k \cdot q_{3, 0}$; $\mu$'s are all zero but $\mu^2$ with identity
morphisms and $\mu^3(q_{1, 0}, q_{2, 1}, q_{3, 2}) = q_{3, 0}$.
This formula coincides with that in Theorem \ref{thm:MAIN1TheFormulaOfAInftyKoszulDual}.

\begin{figure}[hbt]
\centering
\includegraphics[width=8cm]{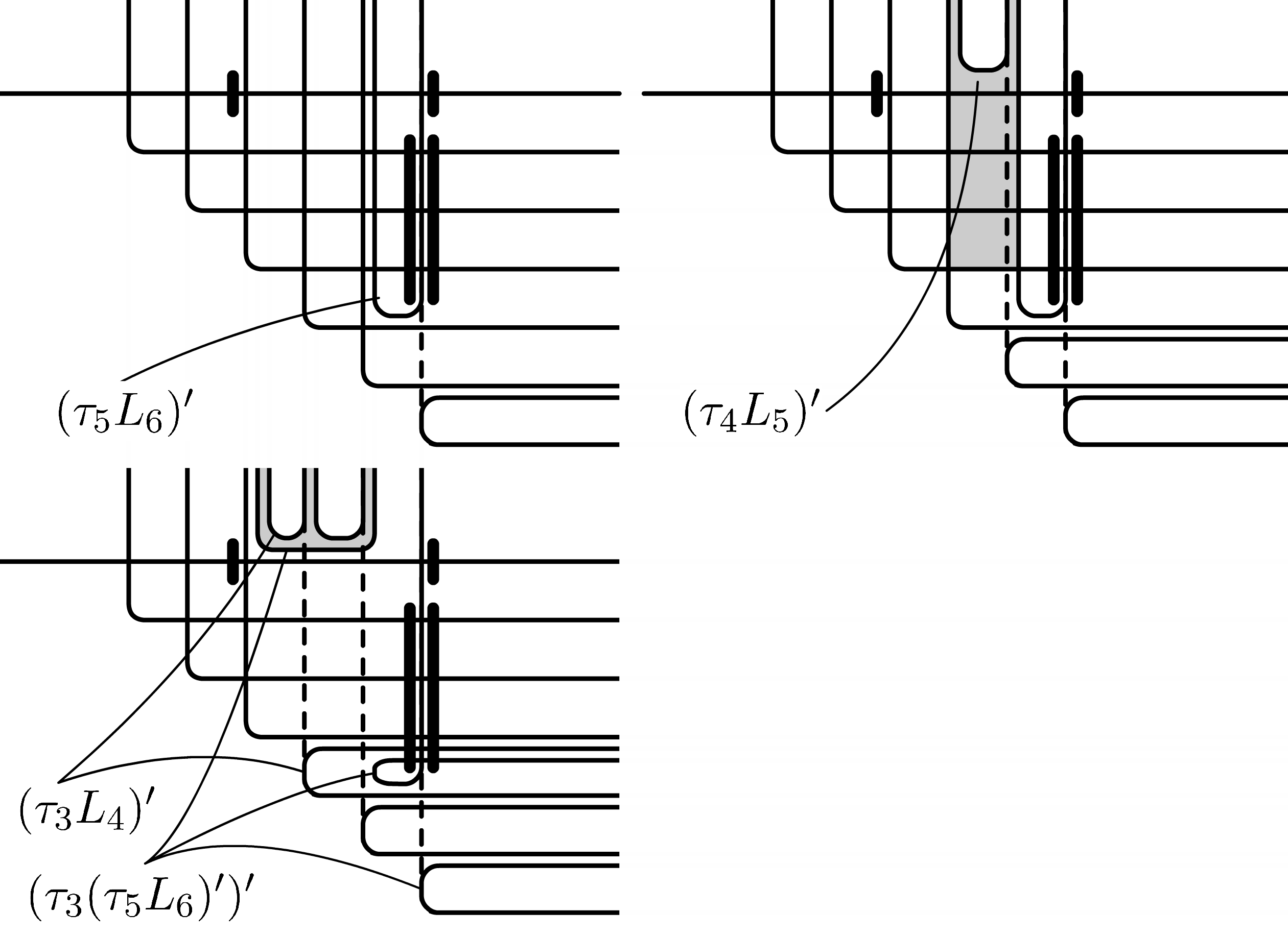}
\caption{$\boldsymbol{L}_2^{(3)}$}\label{fig:S2_1}
\end{figure}

Now, we see the core of $\boldsymbol{S}_2$, especially we study that
which curve $S'_j$  intersects with $S'_6$.
As in Figure \ref{fig:S2_1}, the curve $(\tau_5L_6)'$, obtained by the Dehn
twist and isotopy like the case of $A_1$, does not intersect with $L_4$,
so the curve is twisted
by $L_3$. After an action of isotopy, we can see that $(\tau_3 (\tau_5 L_6)')'$
does not intersect with $L_0$, $L_1$, and $L_2$.
Hence, we can replace $S'_6$ by $(\tau_3 (\tau_5 L_6)')'$.
By the construction, $S'_6$ intersects with $S'_j$ with $j = 5, 3$ at $H_j$.

\begin{figure}[hbt]
\centering
\includegraphics[width=8cm]{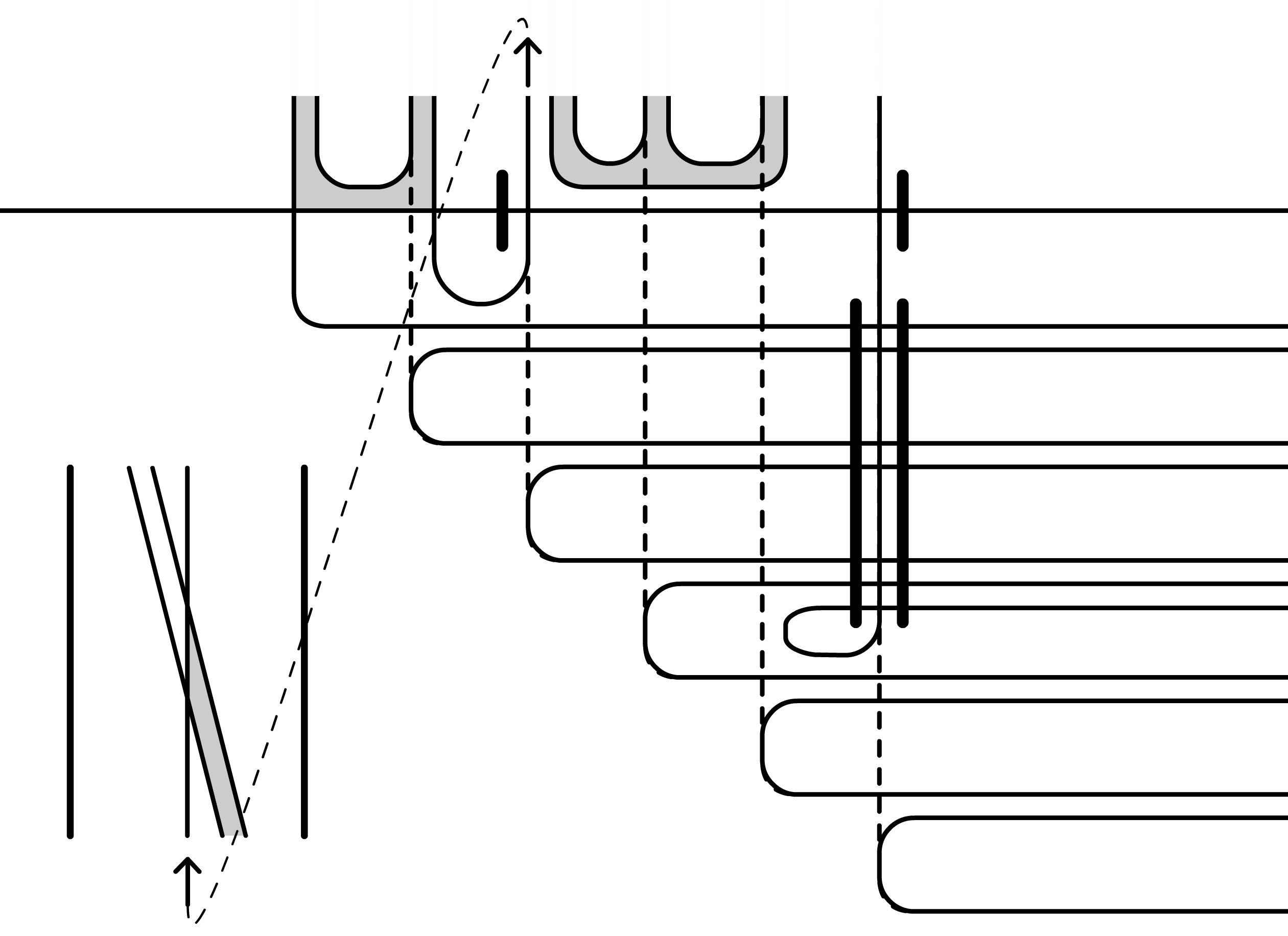}
\caption{$\boldsymbol{S}'_2$}\label{fig:S2_2}
\end{figure}

The upper right part of Figure \ref{fig:S2_1} teaches us that $(\tau_5 L_6)'$,
$(\tau_4 L_5)'$, $L_4$, and $L_3$ bounds
a quadrangle. 
This quadrangle do not be destroyed by the Dehn twists $\tau_0$, $\tau_1$,
$\tau_2$, and $\tau_3$ and isotopies.
By the same argument of the case of $A_1$, we can see that $S_3$, $S_2$,
$S_1$, and $S_0$ forms another quadrangle.
In fact, there are no polygons other than these two quadrangles and the core
of $\boldsymbol{S}_2'$ is as in Figure \ref{fig:S2_2}.

Together with the degree, a directed $A_\infty$-category
$\mathcal{B}_2 \coloneqq \mathcal{F}_2^\to(\boldsymbol{S}^\#_2) \cong \mathcal{F}_2^\to(\boldsymbol{S}'_2{}^\#)$
can be represented as follows: $\hom_{\mathcal{B}_2}^d(S^\#_i, S^\#_j) =
0$ except for $\hom_{\mathcal{B}_2}^0(S^\#_j, S^\#_j) =k \cdot 1_{S^\#_j}$,
$\hom_{\mathcal{B}_2}^1(S^\#_j, S^\#_{j-1}) =\ k \cdot q_{j, j-1}$, $\hom_{\mathcal{B}_2}^2(S^\#_6,
S^\#_3) = k \cdot q_{6, 3}$, and $\hom_{\mathcal{B}_2}^2(S^\#_3,
S^\#_0) = k \cdot q_{3, 0}$; $\mu$'s are all zero but $\mu^2$ with identity
morphisms, $\mu^3(q_{4, 3}, q_{5, 4}, q_{6, 5}) = q_{6, 3}$, and $\mu^3(q_{1,
0}, q_{2, 1}, q_{3, 2}) = q_{3, 0}$.
This formula again coincides with that in Theorem \ref{thm:MAIN1TheFormulaOfAInftyKoszulDual}.

\clearpage

\begin{figure}[hbt]
\centering
\includegraphics[width=8cm]{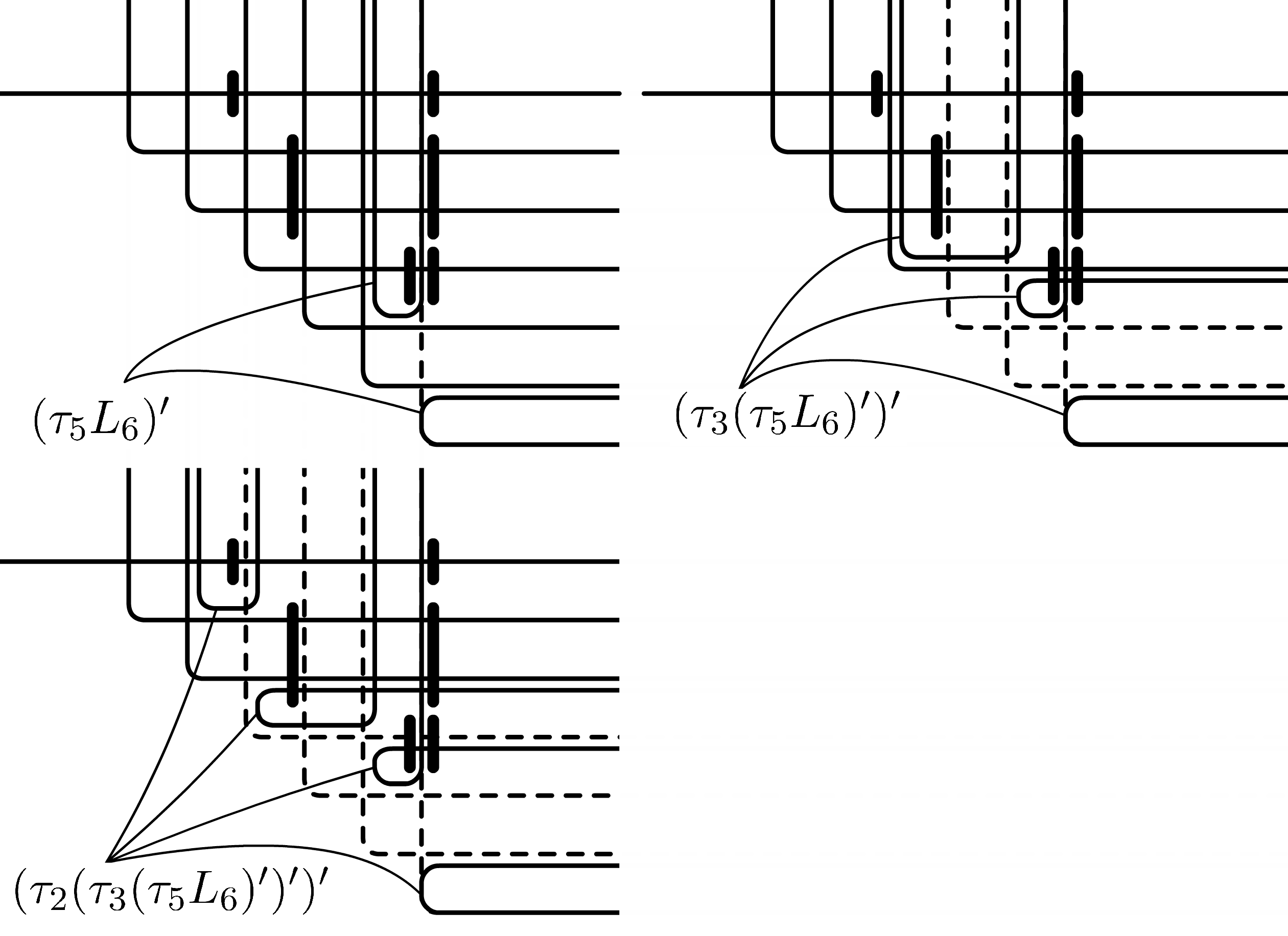}
\caption{from $L_6$ to $S_6$}\label{fig:S3_1}
\end{figure}

Finally, we study the case of $A_3$.
First, we focus on that which curve $S'_j$  does intersect with $S'_6$.
See Figure \ref{fig:S3_1}.
The curve $(\tau_3 (\tau_5 L_6)')'$ intersects with $L_2$ in contrast to
the case of $A_2$.
This difference comes from the bypass in $M_3$ corresponds with a relation
$\alpha_5 \alpha_4$.
This difference induce the difference $d_3(3) = 0 < 2 = d_3(5)
$ and $d_2(3) = 0 = d_2(5)$.
The subscript $2$ of $L_2$ comes from $d_3(5)$.
We can say that $(\tau_2(\tau_3 (\tau_5 L_6)')')'$ intersects
with $L_0 = L_{d_3(3)}$ ``since" $d_3(2) = - \infty < 0 = d_3(3)$.
This naive computation demonstrates the computation of the hom spaces of
an $A_\infty$-Koszul dual of $\mathcal{A}_3$.

Next, we study the order of intersection points of $S'_6$ and other $S'_j$'s.
When we go along $S_6 \simeq \tau_0 S'_6$, we pass through the handles $H_6$, $H_3$, $H_0$,
$H_2$, and $H_5$ in this
order.
These numbers are nothing but $6 = a^{(6)}_0$, $3 = a^{(6)}_2$, $0 = a^{(6)}_4$,
$2
= a^{(6)}_3$, $5 = a^{(6)}_1$.
When we read their subscripts from left to right, the subscripts are
$0, 2, 4, \dots$ and when we read them from right to left, the subscripts
are $1 ,3, \dots$.
As we see, this pattern holds in the general case. 

\begin{figure}[hbt]
\centering
\includegraphics[width=8cm]{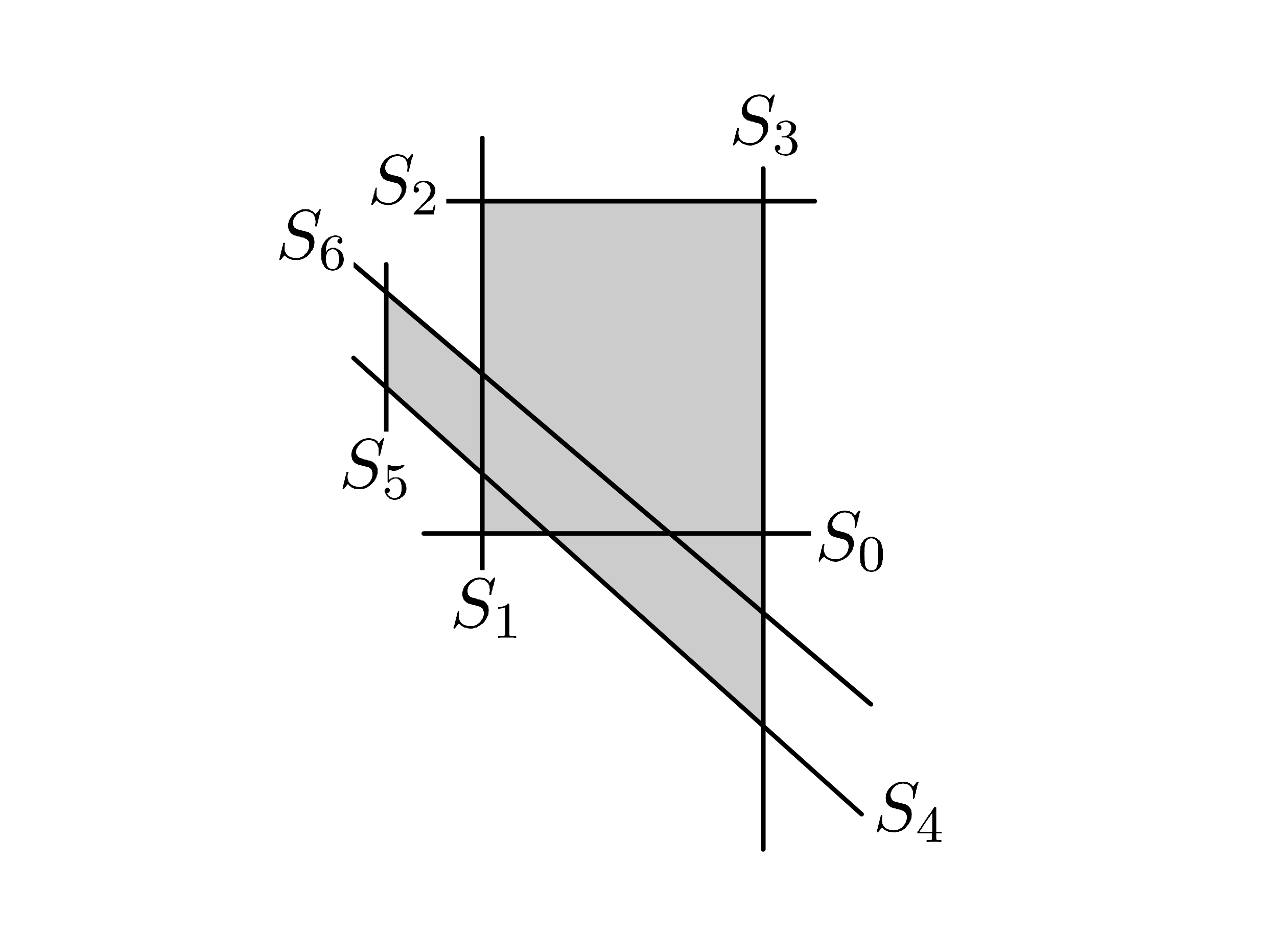}
\caption{$c\boldsymbol{S}_3$}\label{fig:S3_2}
\end{figure}

As in the general case, $\boldsymbol{S}_3$ has the
following properties.
\begin{itemize}
\item There are two quadrangles encircled by $(S_3, S_2, S_1, S_0)$ and $(S_6,
S_5, S_4, S_3)$ and there is a triangle encircled by $(S_4, S_3, S_2)$. \item
$S_j$ intersects with $S_i$ for $i < j$ only when $i$ emerges in the
sequence $\{ a^{(j)}_p \}_{0 \leq p \leq l_j}$.
\item The order of subscripts $i$ of intersection points $q_{j, i}$ of $S_i\cap
S_j$ for $i
< j$ in $S_j$ is $a^{(j)}_2$, $a^{(j)}_4, \dots , a^{(j)}_{l_j},
\dots , a^{(j)}_3$, $a^{(j)}_1$.
\end{itemize}
These properties uniquely determine the core of $\boldsymbol{S}_3$ as in
Figure \ref{fig:S3_2}.

(In this example, the relevant sequences are as follows: $\{ a^{(0)}_i \}
= \{ 0 \}$, $\{ a^{(1)}_i \} = \{ 1, 0\}$, $\{ a^{(2)}_i \} = \{ 2, 1 \}$,
$\{ a^{(3)}_i \} = \{ 3, 2, 0\}$, $\{ a^{(4)}_i \} = \{ 4, 3, 2, 0 \}$, $\{
a^{(5)}_i \} = \{ 5, 4 \}$, and $\{ a^{(6)}_i \} = \{ 6, 5, 3, 2, 0 \}$.
Thus the order of intersection in $S_j$ is as follows: $\varnothing$ for
$S_0$,
$(q_{1, 0})$ for $S_1$, $(q_{2, 1})$ for $S_2$, $(q_{3, 0}, q_{3, 2})$ for
$S_3$, $(q_{4, 2}, q_{4, 0}, q_{4, 3})$ for $S_4$, $(q_{5, 4})$ for $S_5$,
and
$(q_{6, 3}, q_{6, 0}, q_{6, 2}, q_{6, 5})$ for $S_6$.)

By the diagram in Figure \ref{fig:S3_2}, we can check that there is a unique desired polygon if the higher
composition can be non-zero in the sense of the definition of $\mathcal{B}_{S,
T}$ in subsection \ref{subsec:FormulaOfAInfKoszulDual}.
Of course this example supports the Theorem \ref{thm:MAIN1TheFormulaOfAInftyKoszulDual}.

\begin{rem}

A reader who just wants to compute by a picture, one should write the diagram
as
in Figure \ref{fig:S3_2}.
The drawing procedure is as follows:
\begin{enumerate}
\item Compute $\{ a^{(j)}_i \}_{0 \leq i \leq l_j}$.
\item Draw $S_j$ from $j = 0$ to $n$ which intersects only with $S_{a^{(j)}_i}$
and the order
of (subscripts of subscripts) is $2, 4, 6, \dots , l_j, \dots , 5, 3, 1$.
\item Verify that $S_j$ creats a $(t_l - s_l + 1)$-gon encircled by $S_{s_l},
S_{s_l + 1}, \dots , S_{t_l}$ if $j = t_l$ for some $1 \leq l \leq m$.
\item Verify that every desired polygon do exists as in the sense of the
definition of $\mathcal{B}_{S, T}$ in subsection \ref{subsec:FormulaOfAInfKoszulDual},
i.e. if there exist a collection of $S_{j_0}, S_{j_1}, \dots , S_{j_d}$ such
that $S_{j_l}$ and $S_{j_{l+1}}$ intersect, $S_{j_0}$ and $S_{j_d}$ intersect,
and the degree of the intersection points satisfy the degree condition of
higher composition maps, then there exists just one $(d+1)$-gon which contributes
the relevant $\mu^d$.
\end{enumerate}

\end{rem}

\clearpage




\begin{thebibliography}{99}

\bibitem[Ab08]{Ab08}
M. Abouzaid,
\textit{On the Fukaya categories of higher genus surfaces},
Advances in Mathematics 217.3 (2008): 1192-1235.

\bibitem[AKO08]{AKO08}
D. Auroux, L. Katzarkov, D. Orlov,
\textit{Mirror symmetry for weighted projective planes and their noncommutative
deformations}
Annals of Mathematics 167.3 (2008): 867-943.

\bibitem[BGS96]{BGS96}
A. Beilinson, V. Ginzburg, W. Soergel,
\textit{Koszul duality patterns in representation theory},
Journal of the American Mathematical Society 9.2 (1996): 473-527.

\bibitem[BCT09]{BCT09}
A. J. Blumberg, R. L. Cohen, and C. Teleman,
\textit{Open-closed field theories, string topology, and Hochschild homology},
Contemporary Mathematics 504 (2009): 53.

\bibitem[BK90]{BK90}
A. I. Bondal, M. M. Kapranov,
\textit{Enhanced triangulated categories}
Matematicheskii Sbornik 181.5 (1990): 669-683.

\bibitem[EL16]{EL16}
T. Etg\"{u}, Y. Lekili,
{\it Koszul duality patterns in Floer theory},
arXiv preprint arXiv:1502.07922 (2015).

\bibitem[FOOO10]{FOOO10}
K. Fukaya, Y. G. Oh, H. Ohta, K. Ono,
\textit{Lagrangian intersection Floer theory: anomaly and obstruction}
(two volumes), Vol. 46.1 and Vol 46.2. American Mathematical
Soc., 2010.

\bibitem[GK94]{GK94}
V. Ginzburg, M. Kapranov,
\textit{Koszul duality for operads},
Duke mathematical journal 76.1 (1994): 203-272.

\bibitem[HV00]{HV00}
K. Hori, C. Vafa,
{\it Mirror symmetry},
arXiv preprint hep-th/0002222 (2000).

\bibitem[Ka80]{Ka80}
A. Kas,
{\it On the handlebody decomposition associated to a Lefschetz fibration},
Pacific Journal of Mathematics 89.1 (1980): 89-104.

\bibitem[Ko94]{Ko94}
M. Kontsevich, 
\textit{Homological algebra of mirror symmetry},
In Proceedings of the International
Congress of Mathematicians (Zurich, 1994), pages 120-139. Birkhauser, 1995.

\bibitem[LV12]{LV12}
J-L. Loday, and Bruno Vallette,
\textit{Algebraic operads},
Vol. 346. Springer Science \& Business Media, 2012.

\bibitem[L\"o86]{Lo86}
C. L\"ofwall,
\textit{On the subalgebra generated by the one-dimensional elements in the
Yoneda Ext-algebra},
Algebra, algebraic topology and their interactions. Springer Berlin Heidelberg,
1986. 291-338.

\bibitem[LPWZ04]{LPWZ04}
D. M. Lu, J. H. Palmieri, Q. S. Wu, J. J. Zhang,
\textit{$A_{\infty}$-algebras for ring theorists},
Algebra Colloquium 11 (2004), 91-128.

\bibitem[Pr70]{Pr70}
S. B. Priddy,
\textit{Koszul resolutions},
Transactions of the American Mathematical Society 152.1 (1970): 39-60.

\bibitem[Se00]{Se00}
P. Seidel,
\textit{More about vanishing cycles and mutation},
Symplectic geometry and mirror symmetry (Seoul, 2000) (2001): 429-465.

\bibitem[Se01]{Se01}
P. Seidel,
\textit{Vanishing cycles and mutation},
European Congress of Mathematics. Birkh\"{a}user Basel, 2001.

\bibitem[Se08]{Se08}
P. Seidel, 
\textit{Fukaya categories and Picard-Lefschetz theory}, 
European Math. Soc., 2008.

\bibitem[Su16]{Su16}
S. Sugiyama,
\textit{On the Fukaya-Seidel categories of surface Lefschetz fibrations},
arXiv preprint arXiv:1607.02263 (2016).

\bibitem[Va07]{Va07}
B. Vallette,
\textit{A Koszul duality for props},
Transactions of the American Mathematical Society 359.10 (2007): 4865-4943.

\end{thebibliography}
\end{document}